\documentclass[10pt]{article}

\usepackage{amsfonts}
\usepackage{amssymb}
\usepackage{amsmath,appendix}
\usepackage{amscd,array,dsfont,graphicx}%,mathtools}
\usepackage{epsfig,array}
\usepackage{amstext,latexsym,xspace}
\usepackage[margin=2cm]{geometry}
\usepackage{comment}

\def\beq{\begin{equation}}
\def\eqn#1{\beq\label{#1}}
\def\eeq{\end{equation}}

\def\bb {\begin {eqnarray}}
\def\eqnn#1{\bb\label{#1}}
\def\ee {\end {eqnarray}}

\newcommand{\nn}{\nonumber }
\def\tablerule{\noalign{\hrule}}

\def\dim{{\rm dim}}
\def\nt{\noindent}

\def\np{\vfill\eject}
\def\nd{\end{document}}

\def\md{\bigskip}

\def\rank{{\rm rank}}
\def\downcirc#1{\mathop{\circ}\limits_{#1}}
\def\riga{-\kern-4pt - \kern-4pt -}
\font\fat=cmsy10 scaled\magstep5

\def\Bbullet{\raise-3pt\hbox{\fat\char"0F}}

\def\supdot{\cdot\kern-9pt\supset}

\def\subdot{\cdot\kern-10pt\subset}

\def\cupdot{\cdot\kern-6pt\cup}
\def\cupdota{\cdot\kern-6pt\cup}

\def\capdot{\cdot\kern-6pt\cap}

\font\verysmall=cmr5
\def\cupm{\raise2pt\hbox{\verysmall m}\kern-7pt\cup}
\def\cupmm{\raise2pt\hbox{\verysmall m}\kern-10pt\cup}

\def\supm{\raise1.5pt\hbox{\verysmall m}\kern-10pt\supset}
\def\subm{\raise1.5pt\hbox{\verysmall m}\kern-12pt\subset}

\def\black#1{\mathop{\bullet}\limits_{#1}}

\font\tfont=cmbx12 scaled\magstep1 %large
 \font\verysmall=cmr5

\def\Box{
\vbox{ \halign to5pt{\strut##& \hfil ## \hfil \cr &$\kern -0.5pt
\sqcap$ \cr \noalign{\kern -5pt \hrule} }}~}

\def\down{\raise1.5pt\hbox{$\phantom{a}_2$}\downarrow}

\def\downa{\raise1.5pt\hbox{$\phantom{a}_{2\atop m_2}$}\downarrow}

\def\({\left(}
\def\){\right)}

\def\lra{\longrightarrow}

\def\bbz{\mathbb{Z}}%{Z\!\!\!Z}
\def\bbc{\mathbb{C}}%{I\!\!\!\!C}}
\def\bbr{\mathbb{R}}%{{I\!\!R}}
\def\bac{\mathbb{C}}

\def\a{\alpha}

\def\vr{|}

\def\l{\lambda}

\def\D{{\Delta}}

\def\ca{{\cal A}} \def\cb{{\cal B}} 
  
\def\cg{{\cal G}} \def\ch{{\cal H}} 
 \def\ck{{\cal K}} 
\def\cm{{\cal M}} \def\cn{{\cal N}} 
\def\cp{{\cal P}}

\def\th{\theta} \def\Th{\Theta}

\begin{document}

%\rightline{\bf Preliminary text}

%\vskip 2cm

 \centerline{{\tfont Jordan Algebraic Interpretation of }} \vskip
2truemm \centerline{{\tfont Maximal Parabolic Subalgebras :}}\vskip
2truemm \centerline{{\tfont Exceptional Lie Algebras}}

\begin{center}

\vskip 1.5cm

{\bf Vladimir Dobrev}$^1$,~ {\bf Alessio Marrani}$^{2,3}$ \vskip
0.3cm

$^1$ Institute for Nuclear Research and Nuclear Energy,\\ Bulgarian
Academy of Sciences, 72 Tsarigradsko Chaussee,\\ 1784 Sofia,
Bulgaria,  dobrev@inrne.bas.bg

$^2$ Museo Storico della Fisica e Centro Studi e Ricerche Enrico
Fermi,\\ Via Panisperna 89A, I-00184, Roma, Italy

$^3$ Dipartimento di Fisica e Astronomia Galileo Galilei,
Universit`a di Padova,\\ and INFN, sezione di Padova, Via Marzolo 8,
I-35131 Padova, Italy\\ alessio.marrani@pd.infn.it

 \end{center}

\vskip 1.5cm

 \centerline{{\bf Abstract}}

With this paper we start a programme aiming at connecting two vast scientific areas:
Jordan algebras and representation theory. Within representation theory, we
focus on non-compact, real forms of semisimple Lie algebras and groups as well as on the modern
theory of their induced representations, in which a central role is played
by the parabolic subalgebras and subgroups. The aim of the present paper and
its sequels is to present a Jordan algebraic interpretations of maximal
parabolic subalgebras. In this
first paper, we confine ourselves to maximal parabolic subalgebras of the
non-compact real forms of finite-dimensional exceptional Lie algebras, in particular focusing on Jordan algebras of rank 2 and 3.

%\vskip 1.5cm

\newpage

\tableofcontents

%\newpage

\setcounter{equation}{0}
\section{Introduction}

The aim of this paper is to relate two vast scientific areas: Jordan
algebras and representation theory. Jordan algebras are beautiful
mathematical structures that arose together with the advent of quantum
mechanics \cite{Jordan:1933a}--\cite{Ferrar:1972}. They are widely used in
theoretical and mathematical physics, in particular in supergravity and
superstring theory; cf., e.g., \cite{F-Gimon-K}--\cite{Mkrtchyan-Nersessian}%
; we will review their role in such frameworks in Sec. \ref{sec:JA}.
Intriguingly, recent developments related the standard model of particle
physics with the Albert algebra, namely the exceptional rank-3 simple Jordan
algebra over the octonions \cite{DV-1}--\cite{Todorov-3}.

Within the area of representation theory, we focus on non-compact semisimple
Lie algebras and groups starting from the treatment of Gelfand \cite{GeNa},%
\cite{GeV} and Harish-Chandra \cite{HC},\cite{Har} (see also \cite{Bru}), up
to the most advanced approach of Langlands \cite{Lan}, later refined in \cite%
{KnZu}; see also \cite{War},\cite{Vog}. The main building block of the
modern theory of induced representations of non-compact semisimple Lie
algebras and groups are the parabolic subalgebras and subgroups, cf., e.g.,
\cite{Sat}--\cite{Bou}.

We aim at presenting Jordan algebraic interpretations of maximal parabolic
subalgebras of non-compact real forms of semisimple Lie algebras, focussing
on rank-3 and rank-2 (simple and semisimple) Jordan algebras. This is a
rather lengthy and far reaching project, whose proper treatment deserves to
be developed in a series of papers. In the present paper, which inits the
series, we confine ourselves to maximal parabolic subalgebras of the
non-compact real forms of finite-dimensional exceptional Lie algebras,
building over the results and treatment of \cite{Dobpar,Dobparab,Dobk}%
.\bigskip

The paper is organized as follows.

In Section 2 we present a concise treatment of the general theory of
parabolic subalgebras \cite{Dobpar, Dobk}. Then, Section 3 recalls the
classification of maximal parabolic subalgebras of all non-compact real
forms of finite-dimensional exceptional simple Lie algebras; moreover, in
Sec. \ref{parab-rel} we reconsider the definition of parabolically related
non-compact semisimple Lie algebras, and in Table A we present the
classification of maximally parabolically related exceptional Lie algebras,
slightly extending the results of \cite{Dobparab}. An overview of the
general theory of rank-3 and rank-2 Jordan algebras and their symmetries and
related structures, along with a \textit{r\'{e}sum\'{e}} of their relevance
for Maxwell-Einstein (super)gravity theories in various space-time
dimensions, is provided in Sec. \ref{sec:JA}. Such Section also contains
Table B, in which the Jordan algebraic interpretation of Table A is
provided. The content of Table B is derived in a detailed way in the long
Sec. 5, in which the maximal parabolic subalgebras of non-compact real forms
of exceptional Lie algebras are analyzed and obtained through sequences of
maximal embeddings of Lie algebras, which in turn allow for a natural
interpretation in terms of symmetries of Jordan algebras. Finally, a brief
summary and outlook to future developments is given in Sec. 6.

%###

%\np

\bigskip

\setcounter{equation}{0}

\section{Parabolic Subalgebras}

This section follows \cite{Dobpar}.
Let $G$ be a noncompact semisimple Lie group.   Let $\cg$ be the Lie
algebra of $G$, $\th$ be a Cartan involution in $\cg$, and $\cg =
\ck \oplus \cp$ be a Cartan decomposition of $\cg$, so that $\th X
= X, ~X \in \ck$, $\th X = -X , ~X \in \cp$~; $\ck$ is a maximal
compact subalgebra of $\cg$; in general, $\cp$ fits in a (reducible) representation of the algebra $\ck$.
\\ Let $\ca_0$ be a maximal subspace of
$\cp$ which is an abelian subalgebra of $\cg$~; $r =\ $dim $\ca_0$ is
the ~{\it split}~ (or {\it real}) rank of $\cg$, $1 \leq r \leq \ell
= \rank\,\cg$. The subalgebra ~$\ca_0$~ is called a Cartan subspace of
~$\cp$.

 Let $\D_{\ca_0}$ be the root system of the pair $(\cg , \ca_0)$:
\eqn{rroots} \D_{\ca_0} ~\doteq~ \{ \l \in \ca_0^* ~\vr ~ \l \neq 0 ,\
\cg^\l_{\ca_0} \neq 0 \} ~, ~~~\cg^\l_{\ca_0} \doteq \{ X \in \cg ~\vr ~
[Y,X] = \l(Y) X ~, ~~\forall Y\in \ca_0 \} ~. \eeq  The elements of
$\D_{\ca_0}$  are called $\ca_0\,$-{\it restricted roots}.
 For $\l \in \D_{\ca_0}\,$, ~$\cg_{\ca_0}^\l$ are called $\ca_0\,$-{\it\ restricted root
spaces}, ~dim$_R~\cg_{\ca_0}^\l \geq 1$. In a standard way, the $\ca_0\,$-{\it\ restricted roots} are split into positive and
negative restricted roots: ~$\D_{\ca_0} = \D_{\ca_0}^{+} \cup \D_{\ca_0}^{-}\,$.
Then we introduce the corresponding nilpotent subalgebras:
\eqn{nilpot} \cn^\pm ~\doteq~ \mathop{\oplus}\limits_{\l \in \D_{\ca_0}^{\pm} }
~~\cg^\l_{\ca_0} \ . \eeq

Next let $\cm_0$ be the centralizer of $\ca_0$ in $\ck$, i.e.,
 $\cm_0\ \doteq\ \{ X \in \ck ~\vr ~ [X,Y] = 0 , ~\forall Y \in \ca_0 \}$. In
general $\cm_0$ is a compact reductive Lie algebra.

For the Bruhat decomposition, it holds that \cite{Bru}:
\eqn{bruhat} \cg ~=~ \cn^+ \oplus \cm_0 \oplus \ca_0 \oplus \cn^-\ , \eeq
and the subalgebra  ~$\cp_0 \doteq \cm_0 \oplus \ca_0 \oplus \cn^-$~
called a ~{\it minimal parabolic subalgebra}~
of $\cg$. (Note that we may take equivalently
$\cn^+$ instead of $\cn^-\,$.)

   A {\it\ standard parabolic subalgebra}~ is any subalgebra
~$\cp'$~ of $\cg$ containing ~$\cp_0\,$. The number of standard parabolic subalgebras,
including $\cp_0$ and $\cg$, is $2^r$.

Thus, if $r=1$ the only nontrivial parabolic subalgebra is $\cp_0\,$.

Thus, further in this section $r>1$.

Any standard parabolic subalgebra is of the form:
\eqn{parag} \cp' = \cm' \oplus \ca' \oplus \cn'^-\ ,\eeq
so that $\cm' \supseteq \cm_0$, $\ca' \subseteq \ca_0$, $\cn'^-
\subseteq \cn^-$~; $\cm'$ is the centralizer of $\ca'$ in $\cg$
(mod~$\ca')$~; $\cn'^-$~   is comprised
from the negative   root spaces of the restricted
root system $\D_{\ca'}$ of $(\cg, \ca')$.
The decomposition \eqref{parag} is called the Langlands decomposition of ~$\cp'\,$.
One also has the  Bruhat decomposition  \eqref{bruhat} for this general situation:
\eqn{bruhatg}\cg ~=~ \cn'^+ \oplus \ca' \oplus \cm'
\oplus \cn'^- ~, \eeq
where $\cn'^+ = \th \cn'^-$.

The standard parabolic subalgebras may be described explicitly using the restricted simple root system
~$\D^S_{\ca_0}\,$, such that if $\l\in\D_{\ca_0}^+$, (resp. $\l\in\D_{\ca_0}^-$), one has:
\eqn{decomp} \l = \sum_{i=1}^r \,n_i\, \l_i\ , ~~~~ \l_i \in\D^S_{\ca_0}\ ,
\quad {\rm all}~ n_i\geq 0, ~~~ ({\rm resp.~ all} ~ n_i\leq 0)\ .\eeq

We shall follow Warner \cite{War}, where one may find all references to the original
mathematical work on parabolic subalgebras. For a short
formulation one may say that the parabolic subalgebras correspond
to the various subsets of $\D^S_{\ca_0}$ - hence their number $2^r$. To
formalize this let us denote: ~$S_r = \{ 1,2,\ldots,r\}$, and  let
~$\Th$~ denote any subset of ~$S_r\,$. Let ~$\D^\pm_\Th\in \D_{\ca_0}$~ denote all positive/negative restricted
roots which are linear combinations of the simple restricted roots $\l_i\,, \forall\,i\in \Th$.
Then a standard parabolic subalgebra corresponding to ~$\Th$~ will be denoted by  ~$\cp_\Th$~
and is given explicitly as:
\eqn{stpar} \cp_\Th ~=~ \cp_0\ \oplus\ \cn^+(\Th) \ , \qquad \cn^+(\Th) ~\doteq~
\mathop{\oplus}\limits_{\l\in \Delta^+_\Th}  \ \cg^\l_{\ca_0} \ .\eeq
Clearly, ~$\cp_\emptyset ~=~ \cp_0\,$, ~~$\cp_{S_r} ~=~ \cg\,$, ~since
~$\cn^+(\emptyset) ~=~ 0\,$, ~~$\cn^+(S_r) ~=~ \cn^+\,$.
Further, we need to bring \eqref{stpar} in the  form \eqref{parag}. First, define ~$\cg(\Th)$~ as the algebra generated
by $\cn^+(\Th)$ and $\cn^-(\Th) \doteq \th \cn^+(\Th)\,$.
Next, define $\ca(\Th) \doteq \cg(\Th) \cap \ca_0$,
and ~$\ca_\Th$~ as the orthogonal complement (relative to the Euclidean structure of $\ca_0$)
of $\ca(\Th)$ in $\ca_0$. Then ~$\ca_0 = \ca(\Th) \oplus \ca_\Th\,$.
Note that ~$\dim\, \ca(\Th) = |\Th|$,  ~$\dim\, \ca_\Th = r-|\Th|$.
Next, define:
\eqn{nilp}  \cn^+_\Th ~\doteq~ \mathop{\oplus}\limits_{\l\in \D^+_{\ca_0}-\Delta^+_\Th}  \ \cg^\l_{\ca_0}\ ,
\quad \cn^-_\Th ~\doteq~ \th\cn^+_\Th \ . \eeq
Then ~$\cn^\pm = \cn^\pm(\Th) \oplus \cn^\pm_\Th\,$.
Next, define ~$\cm_\Th \doteq \cm_0 \oplus \ca(\Th)\oplus \cn^+(\Th) \oplus \cn^-(\Th) $.
Then $\cm_\Th$ is the centralizer of $\ca_\Th$ in $\cg$ (mod $\ca_\Th$).
Finally, we can derive:
\eqnn{paragg}  \cp_\Th\ &=&\ \cp_0 \oplus \cn^+(\Th) = \cm_0 \oplus \ca_0 \oplus \cn^- \oplus \cn^+(\Th)\ =\cr
&=&\ \cm_0 \oplus \ca(\Th) \oplus \ca_\Th \oplus \cn^-(\Th) \oplus \cn^-_\Th \oplus \cn^+(\Th)\ =\cr
&=&\ \Big( \cm_0 \oplus \ca(\Th) \oplus \cn^-(\Th) \oplus \cn^+(\Th) \Big) \oplus \ca_\Th \oplus \cn^-_\Th =\cr
&=&\ \cm_\Th \oplus \ca_\Th \oplus \cn^-_\Th \ .\ee
Thus, we have rewritten explicitly the standard parabolic $\cp_\Th$ in the desired form \eqref{parag}.
The associated (generalized) Bruhat decomposition \eqref{bruhatg} is given now explicitly as:
\eqnn{bruhag}   \cg\ &=&\  \cn^+ \oplus \cp_0 = \cn^+_\Th \oplus \cn^+(\Th) \oplus \cp_0
= \cn^+_\Th \oplus   \cp_\Th \ =\cr &=&\  \cn^+_\Th \oplus \cm_\Th \oplus \ca_\Th \oplus \cn^-_\Th \ .
\ee

In this paper we concentrate on
the {\it maximal\ } parabolic subalgebras which
correspond to $\Th$ of the form:
\eqn{thitb} \Th^{\rm max}_j ~=~ S_r\backslash \{j\} \ , ~~~1 \leq j\leq r \ . \eeq
~$\dim\, \ca(\Th^{\rm max}_j) = r-1$,  ~$\dim\, \ca_{\Th^{\rm max}_j} = 1$.

%\np

 \setcounter{equation}{0}
\section{\label{sec:MP}Maximal Parabolic Subalgebras of the Exceptional Lie Algebras}
%\label{Max-parab}

The material of this Section is taken from \cite{Dobpar}. Here we list the real forms of the exceptional simple Lie
algebras and their maximal parabolic subalgebras.

First we mention a class of real forms, the maximally non-compact, or ~{\it split}, real forms, which exist for every complex simple algebra.
For the split real forms we can use the same basis as for the
corresponding complex simple Lie algebra $\cg^\bbc$, but over $\bbr$.
Restricting $\bbc\lra\bbr$ one obtains the
Bruhat decomposition of $\cg$ (with $\cm_0=0$) from  the triangular decomposition of $\cg^\bbc
= \cg^+ \oplus \ch^\bbc \oplus \cg^-$, and obtains the
   minimal parabolic subalgebras $\cp_0$ from the Borel subalgebra
~$\cb = \ch^\bbc \oplus \cg^+$, (or $\cg^-$ instead of  $\cg^+$).
Furthermore, in this case ~$\dim_\bbr\,\ck ~=~ \dim_\bbr\,\cn^\pm$.

For the real forms in general we need to use the Satake diagrams
\cite{Sat,War}. A Satake diagram has as a starting point
the Dynkin diagram of the corresponding complex form. For a split
real form it remains the same \cite{War}. In the other cases some dots of the Dynkin diagram are
painted in black - these considered by themselves are Dynkin
diagrams of the compact semisimple factors ~$\cm_0$~ of the minimal
parabolic subalgebras. Further, there are arrows connecting some
nodes which use the $\bbz_2$ symmetry of some Dynkin diagrams.  Then
the reduced root systems are described by Dynkin-Satake diagrams
which are obtained from the Satake diagrams by dropping the black
nodes, identifying the arrow-related nodes, and adjoining all nodes
in a connected Dynkin-like diagram, but in addition noting the
multiplicity of the reduced roots (which is in general different
from 1).

Below all dimensions of vector spaces are real.

Note that from now on we shall omit the superscript ~$^{\max}$~
in order to simplify the notation. Instead we include specification
for each maximal parabolics.

\subsection{EI :~   $E_{6(6)}$}

\nt
The split real form of ~$E_6$~ is denoted  as  ~$E_{6(6)}$~ or ~$E'_6\,$.
The maximal compact subgroup is ~$\ck \cong sp(4)$, $\cp =\wedge _{0}^{4}\mathbf{8}$, ~$\dim\,\cp = 42$,
~$\dim\,\cn^\pm = \dim\,\ck = 36$.
This real form does not have discrete series representations.

Since ~$E_{6(6)}$~ is split the Satake diagram coincides with  the  Dynkin diagram:
\eqn{satsixa}  \downcirc{{\a_1}} \riga\downcirc{\a_2} \riga
\downcirc{{\a_3}}\kern-8pt\raise11pt\hbox{$\vr$}
\kern-3.5pt\raise22pt\hbox{$\circ {{\scriptstyle{\a_6}}}$}
\riga\downcirc{{\a_4}} \riga\downcirc{{\a_5}} \eeq

Taking into account the above enumeration of simple roots and \eqref{thitb}
the maximal parabolic subalgebras are determined by \cite{Dobpar}:
\eqnn{maxsixa}
 \cp^{ 6(6)}_i ~&=&~ \cm^{ 6(6)}_i\ \oplus\ so(1,1)\ \oplus\ \cn^{\, 6(6)}_i\ , \quad i=1,2,3,4\ ;\\
 \cm^{ 6(6)}_1  \ &\cong&\ \cm^{ 6(6)}_5 ~=~ so(5,5)
\ , \qquad \dim\,\cn_1^{\, 6(6)\pm} ~=~ 16
\cr
\cm^{ 6(6)}_2 &=&\ sl(6,\bbr)
\ , \qquad \dim\,\cn_2^{\, 6(6)\pm} ~=~ 21\cr
 \cm^{ 6(6)}_3\ &\cong&  \ \cm^{ 6(6)}_6 ~=~
sl(5,\bbr) \oplus sl(2,\bbr)
\ , \qquad \dim\,\cn_3^{\, 6(6)\pm} ~=~ 25
\cr
 \cm^{ 6(6)}_4\ & = &\
sl(3,\bbr) \oplus sl(3,\bbr) \oplus sl(2,\bbr)
\ , \qquad \dim\,\cn_4^{\, 6(6)\pm} ~=~ 29
\nn \ee
There are essentially four maximal cases (instead of six) due to the symmetry of the Satake
diagram. The two isomorphic occurrences are distinguished by the fact which node we delete to obtain a maximal parabolic.\\
 Clearly, no maximal parabolic subalgebra is cuspidal.

\subsection{EII :~   $E_{6(2)}$}

\nt
Another real form of ~$E_6$~ is denoted as ~$E_{6(2)}\,$ or ~$E''_6\,$.
The maximal compact subgroup is ~$\ck \cong su(6)\oplus su(2)$, $\cp =\left( \mathbf{20},\mathbf{2}\right) $, where $\mathbf{20}=\wedge ^{3}%
\mathbf{6}$, ~$\dim\,\cp = 40$,
~$\dim\,\cn^\pm = 36$, ~$\cm_0 \cong u(1)\oplus u(1)$.
This real form has discrete series representations.

The Satake diagram is:
\eqn{satsixaz} \underbrace{ \downcirc{{\a_1}} \riga \underbrace{ \downcirc{\a_2} \riga
\downcirc{{\a_3}}\kern-8pt\raise11pt\hbox{$\vr$}
\kern-3.5pt\raise22pt\hbox{$\circ {{\scriptstyle{\a_6}}}$}
\riga\downcirc{{\a_4}}} \riga\downcirc{{\a_5}} } \eeq

The split rank is four and thus there are four
maximal parabolic subalgebras given by \cite{Dobpar}:
\eqnn{maxsixb}
 \cp^{ 6(2)}_i ~&=&~ \cm^{ 6(2)}_i\ \oplus\ so(1,1)\ \oplus\ \cn^{\, 6(2)}_i\ , \quad i=1,2,3,4\ ;\\
\cm^{ 6(2)}_1 &=&\ so(5,3) \oplus u(1)
\ , \qquad \dim\,\cn_1^{\, 6(2)\pm} ~=~ 24\cr
 \cm^{ 6(2)}_2\ &=&\
sl(3,\bbr) \oplus sl(2,\bbc)_\bbr \oplus u(1)
\ , \qquad \dim\,\cn_2^{\, 6(2)\pm} ~=~ 31 \cr
\cm^{ 6(2)}_3\ &=&  \
sl(3,\bbc)_\bbr \oplus sl(2,\bbr)
\ , \qquad \dim\,\cn_3^{\, 6(2)\pm} ~=~ 29
\cr
 \cm^{ 6(2)}_4  \ &=&\ su(3,3)
\ , \qquad \dim\, \cn_4^{\, 6(2)\pm} ~=~ 21
\nn \ee
Only the last case in \eqref{maxsixb} is cuspidal, and it has
  highest/lowest weight representations.

%\eqn\cm^{\rm max} ~=~ \begin{cases}
%su(3,3) \ ,~&~ \dim\,(\cn^\pm)^{\rm max}= 21, \cr
%so(5,3) \oplus u(1) \ ,~&~ \dim\,(\cn^\pm)^{\rm max}=24, \cr
%sl(3,\bbc) \oplus sl(2,\bbr)   \ ,~&~ \dim\,(\cn^\pm)^{\rm max}=29, \cr
%so(3,1) \oplus u(1)\oplus sl(3,\bbr)   \ ,~&~ \dim\,(\cn^\pm)^{\rm max}=31,
%\end{cases} \eeq

\subsection{EIII :~   $E_{6(-14)}$}

\nt
Another real form of ~$E_6$~ is denoted as ~$E_{6(-14)}\,$ or ~$E'''_6\,$.
The maximal compact subgroup is ~$\ck \cong so(10)\oplus so(2)$, ~$\cp = \mathbf{16}_{-3}\oplus \overline{\mathbf{16}}_{3}$, $\dim\,\cp = 32$,
~$\dim\,\cn^\pm = 30$, ~$\cm_0 \cong so(6)\oplus so(2)$.
This real form has discrete series representations (and highest/lowest weight
representations).

The Satake diagram is:
\eqn{satsixay} \underbrace{ \downcirc{{\a_1}} \riga  \black{\a_2} \riga
\black{{\a_3}}\kern-8pt\raise11pt\hbox{$\vr$}
\kern-3.5pt\raise22pt\hbox{$\circ {{\scriptstyle{\a_6}}}$}
\riga\black{{\a_4}} \riga\downcirc{{\a_5}} } \eeq

The split rank is two and thus there are two
maximal parabolic subalgebras given by \cite{Dobpar}:
\eqnn{maxsixc}
 \cp^{6(-14)}_i ~&=&~ \cm^{6(-14)}_i\ \oplus\ so(1,1)\ \oplus\ \cn^{\,6(-14)}_i\ , \quad i=1,2\ ;\\
\cm^{6(-14)}_1 &=&\ so(7,1) \oplus so(2)
\ , \qquad \dim\,\cn_1^{\,6(-14)\pm} ~=~ 24\cr
\cm^{6(-14)}_2 &=&\ su(5,1)
\ , \qquad \dim\,\cn_2^{\,6(-14)\pm} ~=~ 21
\nn \ee     The 2nd is cuspidal and it has
  highest/lowest weight representations.

\subsection{EIV :~   $E_{6(-26)}$}

\nt
Another real form of ~$E_6$~ is denoted as ~$E_{6(-26)}\,$ or ~$E^{iv}_6\,$. This the minimally non-compact real form of $E_6$.
The maximal compact subgroup is ~$\ck \cong f_4$, $\cp = \mathbf{26}$, ~$\dim\,\cp = 26$,
~$\dim\,\cn^\pm = 24$, ~$\cm \cong so(8)$.
This real form does not have discrete series representations.

The Satake diagram is:
\eqn{satsixd}
\downcirc{ {\a_1}} \riga \black{\a_2} \riga
\black{{\a_3}}\kern-8pt\raise11pt\hbox{$\vr$}
\kern-3.5pt\raise22pt\hbox{$\bullet
{{\scriptstyle{\a_6}}}$} \riga\black{{\a_4}}
\riga\downcirc{{\a_5}}  \eeq

The split rank is equal to 2, thus we have \cite{Dobpar}:
\eqnn{maxsixd}
 \cp^{6(-26)}_i ~&=&~ \cm^{6(-26)}_i\ \oplus\ so(1,1)\ \oplus\ \cn^{\,6(-26)}_i\,  \quad i=1,2\ ;\\
\cm^{6(-26)}_1 ~&\cong&~ \cm^{6(-26)}_2 ~=~ so(9,1)  \   , \qquad
\dim\,\cn^{\,6(-26)\pm}_i ~=~ 16 \, \quad i=1,2.\ \nn\ee
We distinguish the two isomorphic maximal parabolic subalgebras by the fact which noncompact node on
the Satake diagram we delete - the first - for ~$\cm^{6(-26)}_1\,$, or the last,
for ~$\cm^{6(-26)}_2\,$. This case is not cuspidal.

\subsection{EV :~   $E_{7(7)}$}

\nt
The split real form of ~$E_7$~ is denoted as ~$E_{7(7)}\,$ or ~$E'_7\,$.
The maximal compact subgroup is ~$\ck \cong su(8)$, $\cp = \wedge ^{4}\mathbf{8}$, ~$\dim\,\cp = 70$,
~$\dim\,\cn^\pm = \dim\,\ck = 63$.
This real form has discrete series representations.

The Dynkin-Satake diagram is:
\eqn{satseva}
\downcirc{{\a_1}} \riga\downcirc{\a_2} \riga
\downcirc{{\a_3}}\kern-8pt\raise11pt\hbox{$\vr$}
\kern-3.5pt\raise22pt\hbox{$\circ
{{\scriptstyle{\a_7}}}$} \riga\downcirc{{\a_4}}
\riga\downcirc{{\a_5}} \riga\downcirc{{\a_6}} \eeq

Taking into account the above enumeration of simple roots and
  \eqref{thitb} the maximal parabolic subalgebras are determined by (only the first case is cuspidal) :
\eqnn{maxseva}
 \cp^{7(7)}_i ~&=&~ \cm^{7(7)}_i\ \oplus\ so(1,1)\ \oplus\ \cn^{\,7(7)}_i\ , \quad i=1,\ldots,7\ ;\\
  \cm^{7(7)}_1   &=& so(6,6)
\ , \qquad \dim\,\cn_1^{\,7(7)\pm} ~=~ 33
\cr
  \cm^{7(7)}_2 &=& sl(6,\bbr) \oplus sl(2,\bbr)
\ , \qquad \dim\,\cn_2^{\,7(7)\pm} ~=~ 47
\cr
  \cm^{7(7)}_3 &=& sl(4,\bbr) \oplus sl(3,\bbr) \oplus sl(2,\bbr)
\ , \qquad \dim\,\cn_3^{\,7(7)\pm} ~=~ 53
\cr
 \cm^{7(7)}_4 &=& sl(5,\bbr) \oplus sl(3,\bbr)
\ , \qquad \dim\,\cn_4^{\,7(7)\pm} ~=~ 50
\cr
  \cm^{7(7)}_5   &=& so(5,5) \oplus sl(2,\bbr)
\ , \qquad \dim\,\cn_5^{\,7(7)\pm} ~=~ 42
\cr
  \cm^{7(7)}_6 &=& E_{6(6)}
\ , \qquad \dim\,\cn_6^{\,7(7)\pm} ~=~ 27
\cr
 \cm^{7(7)}_7 &=& sl(7,\bbr)
\ , \qquad \dim\,\cn_7^{\,7(7)\pm} ~=~ 42
\nn\ee

\subsection{EVI :~   $E_{7(-5)}$}

\nt
Another real form of ~$E_7$~ is denoted as ~$E_{7(-5)}\,$ or $E''_7\,$.
The maximal compact subgroup is ~$\ck \cong so(12)\oplus su(2)$, $\cp = \left( \mathbf{32},\mathbf{2}\right) $, where $\mathbf{32}$ is the semispinor irrepr. of $so(12)$, ~$\dim\,\cp = 64$,
~$\dim\,\cn^\pm = 60$, ~$\cm \cong su(2)\oplus su(2) \oplus su(2)$.
This real form has discrete series representations.

The Satake diagram is:
\eqn{satsevb}
\downcirc{{\a_1}} \riga\downcirc{\a_2} \riga
\downcirc{{\a_3}}\kern-8pt\raise11pt\hbox{$\vr$}
\kern-3.5pt\raise22pt\hbox{$\bullet
{{\scriptstyle{\a_7}}}$} \riga\black{{\a_4}}
\riga\downcirc{{\a_5}} \riga\black{{\a_6}} \eeq

The split rank is equal to 4, thus, there are four
maximal parabolic subalgebras given by \cite{Dobpar}:
\eqnn{maxsevb}
\cp^{7(-5)}_i ~&=&~ \cm^{7(-5)}_i\ \oplus\ so(1,1)\ \oplus\ \cn^{\,7(-5)}_i\ , \quad i=1,\ldots,4\ ;\\
\cm^{7(-5)}_1 ~&=&~
so^*(12)\ ,\qquad \dim\,\cn_1^{\,7(-5)\pm} ~=~ 33 \cr
\cm^{7(-5)}_2 ~&=&~
so(7,3) \oplus su(2)  \ ,
\qquad \dim\,\cn_2^{\,7(-5)\pm} ~=~ 42 \cr
\cm^{7(-5)}_3 ~&=&~
su^*(6) \oplus   sl(2,\bbr) \ ,
\qquad \dim\,\cn_3^{\,7(-5)\pm} ~=~ 47 \cr
\cm^{7(-5)}_4 ~&=&~
so(5,1) \oplus  sl(3,\bbr)  \oplus su(2) \ ,
\qquad \dim\,\cn_4^{7(-5)} ~=~ 53
\nn\ee
the first case being  cuspidal and having highest/lowest weight representations.

\subsection{EVII :~   $E_{7(-25)}$}

\nt
Another real form, the minimally non-compact one, of ~$E_7$~ is denoted as ~$E_{7(-25)}\,$ or ~$E'''_7\,$.
The maximal compact subgroup is ~$\ck \cong e_6\oplus so(2)$, $\cp = \mathbf{27}_{2}\oplus \overline{\mathbf{27}}_{-2}$, ~$\dim\,\cp = 54$,
~$\dim\,\cn^\pm = 51$, ~$\cm \cong so(8)$.
This real form has discrete series representations (and highest/lowest
weight representations).

The Satake diagram is:
\eqn{satsevbb}
\downcirc{{\a_1}} \riga\black{\a_2} \riga
\black{{\a_3}}\kern-8pt\raise11pt\hbox{$\vr$}
\kern-3.5pt\raise22pt\hbox{$\bullet
{{\scriptstyle{\a_7}}}$} \riga\black{{\a_4}}
\riga\downcirc{{\a_5}} \riga\downcirc{{\a_6}} \eeq

The split rank is equal to 3, thus, there are three
maximal parabolic subalgebras given by \cite{Dobpar}:
\eqnn{maxsevc}
\cp^{7(-25)}_i ~&=&~ \cm^{7(-25)}_i\ \oplus\ so(1,1)\ \oplus\ \cn^{\,7(-25)}_i\ , \quad i=1,2,3\ ;\\
\cm^{7(-25)}_1 ~&=&~ so(10,2)\ ,
\qquad \dim\,\cn_3^{\,7(-25)\pm} ~=~ 33 \cr
\cm^{7(-25)}_2 ~&=&~
so(9,1) \oplus sl(2,\bbr)    \ ,
\qquad \dim\,\cn_2^{\,7(-25)\pm} ~=~ 42 \cr
\cm^{7(-25)}_3 ~&=&~   E_{6(-26)}\   ,
\qquad \dim\,\cn_1^{\,7(-25)\pm} ~=~ 27
 \nn\ee
the first case being  cuspidal and having highest/lowest weight representations.

\subsection{EVIII :~   $E_{8(8)}$}

\nt
The split real form of ~$E_8$~ is denoted as ~$E_{8(8)}\,$ or ~ $E'_8\,$.
The maximal compact subgroup ~$\ck \cong so(16)$, $\cp = \mathbf{128}$ is the semispinor irrepr. of $so(16)$, ~$\dim\,\cp = 128$,
~$\dim\,\cn^\pm = \dim\,\ck = 120$.
This real form has discrete series representations.

The Dynkin-Satake diagram is:
\eqn{sateighta}
\downcirc{{\a_1}} \riga\downcirc{\a_2} \riga
\downcirc{{\a_3}}\kern-8pt\raise11pt\hbox{$\vr$}
\kern-3.5pt\raise22pt\hbox{$\circ
{{\scriptstyle{\a_8}}}$} \riga\downcirc{{\a_4}}
\riga\downcirc{{\a_5}} \riga\downcirc{{\a_6}}  \riga\downcirc{{\a_7}}  \eeq

Taking into account the above enumeration of simple roots and \eqref{thitb}
the maximal  parabolic subalgebras   are determined by:
\eqnn{maxeighta}
\cp^{8(8)}_i ~&=&~ \cm^{8(8)}_i\ \oplus\ so(1,1)\ \oplus\ \cn^{\,8(8)}_i\ , \quad i=1,\ldots,8\ ;\\
  \cm^{8(8)}_1   ~&=&~ so(7,7)
\ , \qquad \dim\,\cn_1^{\,8(8)\pm} ~=~ 78
\cr
  \cm^{8(8)}_2 ~&=&~ sl(7,\bbr) \oplus sl(2,\bbr)
\ , \qquad \dim\,\cn_2^{\,8(8)\pm} ~=~ 98
\cr
  \cm^{8(8)}_3 ~&=&~ sl(5,\bbr) \oplus sl(3,\bbr) \oplus sl(2,\bbr)
\ , \qquad \dim\,\cn_3^{\,8(8)\pm} ~=~ 106
\cr
 \cm^{8(8)}_4 ~&=&~ sl(5,\bbr) \oplus sl(4,\bbr)
\ , \qquad \dim\,\cn_4^{\,8(8)\pm} ~=~ 104
\cr
  \cm^{8(8)}_5   ~&=&~ so(5,5) \oplus sl(3,\bbr)
\ , \qquad \dim\,\cn_5^{\,8(8)\pm} ~=~ 97
\cr
  \cm^{8(8)}_6 ~&=&~ E_{6(6)} \oplus sl(2,\bbr)
\ , \qquad \dim\,\cn_6^{\,8(8)\pm} ~=~ 83
\cr
 \cm^{8(8)}_7 ~&=&~  E_{7(7)}
\ , \qquad \dim\,\cn_7^{\,8(8)\pm} ~=~ 57
\cr
\cm^{8(8)}_8 ~&=&~ sl(8,\bbr)
\ , \qquad \dim\,\cn_8^{\,8(8)\pm} ~=~ 92
\nn\ee
Only the seventh case is cuspidal.

\subsection{EIX :~   $E_{8(-24)}$}

\nt
Another real form of ~$E_8$~ is denoted as ~$E_{8(-24)}\,$ or ~$E''_8\,$.
The maximal compact subgroup is ~$\ck \cong e_7\oplus su(2)$, $\cp = \left( \mathbf{56},\mathbf{2}\right) $, ~$\dim\,\cp = 112$,
~$\dim\,\cn^\pm = 108$, ~$\cm \cong so(8)$.
 This real form has discrete series representations.

The Satake diagram
\eqn{sateightb}
\downcirc{{\a_1}} \riga\black{\a_2} \riga
\black{{\a_3}}\kern-8pt\raise11pt\hbox{$\vr$}
\kern-3.5pt\raise22pt\hbox{$\bullet
{{\scriptstyle{\a_8}}}$} \riga\black{{\a_4}}
\riga\downcirc{{\a_5}} \riga\downcirc{{\a_6}}
\riga\downcirc{{\a_7}}\eeq

The split rank is equal to 4, thus, there are four
maximal parabolic subalgebras given by \cite{Dobpar}:
\eqnn{maxeightb}
\cp^{8(-24)}_i ~&=&~ \cm^{8(-24)}_i\ \oplus\ so(1,1)\ \oplus\ \cn^{\,8(-24)}_i\ , \quad i=1,\ldots,4\ ;\\
\cm^{8(-24)}_1 ~&=&~  so(11,3)  \ ,
\qquad \dim\,\cn_1^{\,8(-24)\pm} ~=~ 78 \cr
\cm^{8(-24)}_2 ~&=&~  so(9,1) \oplus  sl(3,\bbr)  \ ,
\qquad \dim\,\cn_2^{\,8(-24)\pm} ~=~ 97 \cr
\cm^{8(-24)}_3 ~&=&~ E_{6(-26)} \oplus sl(2,\bbr)   \ ,
\qquad \dim\,\cn_3^{\,8(-24)\pm} ~=~ 83 \cr
\cm^{8(-24)}_4 ~&=&~ E_{7(-25)} \ ,
\qquad \dim\,\cn_4^{\,8(-24)\pm} ~=~ 57
\nn  \ee
the last case being  cuspidal and having highest/lowest weight representations.

\subsection{FI :~   $F_{4(4)}$}

\nt
The split real form of ~$F_4$~ is denoted as ~$F_{4(4)}\,$ or ~$F'_4\,$.
The maximal compact subgroup ~$\ck \cong sp(3) \oplus su(2)$, $\cp = \left( \mathbf{14'},\mathbf{2}\right) $, where $\mathbf{14'}=\wedge_0 ^{3}%
\mathbf{6}$, ~$\dim\,\cp = 28$,
~$\dim\,\cn^\pm = \dim\,\ck = 24$.
This real form has discrete series representations.

The Dynkin-Satake diagram is:
\eqn{satfour}  \downcirc{{\a_1}}
\riga\downcirc{{\a_2}} \Longrightarrow\downcirc{{\a_3}}
\riga\downcirc{{\a_4}} \eeq

 The maximal parabolic subalgebras are given by \cite{Dobpar} (being $F_4$ non-simply laced, we denote the short/long nature of the roots) :
\eqnn{maxfoura}
\cp^{4(4)}_i ~&=&~ \cm^{4(4)}_i\ \oplus\ so(1,1)\ \oplus\ \cn^{\,4(4)}_i\ , \quad i=1,2,3,4\ ;\\
\cm^{4(4)}_1 ~&=&~ sl(3,\bbr)_S \oplus sl(2,\bbr)_L \ ,
\qquad \dim\,\cn_1^{\,4(4)\pm} ~=~ 20 \ , \cr
&& (11 ~{\rm long~roots}, 9 ~{\rm short~roots}), \cr
\cm^{4(4)}_2 ~&=&~ sl(3,\bbr)_L \oplus sl(2,\bbr)_S \ ,
\qquad \dim\,\cn_2^{\,4(4)\pm} ~=~ 20 \ ,\cr
&& (9 ~{\rm long~roots}, 11 ~{\rm short~roots}), \cr
\cm^{4(4)}_3 ~&=&~ sp(3,\bbr) \ ,
\qquad \dim\,\cn_3^{\,4(4)\pm} ~=~ 15 \cr
&& (9 ~{\rm long~roots}, 6 ~{\rm short~roots}), \cr
\cm^{4(4)}_4 ~&=&~ so(4,3) \ ,
\qquad \dim\,\cn_4^{\,4(4)\pm} ~=~ 15 \cr
&& (6 ~{\rm long~roots}, 9 ~{\rm short~roots})
\nn \ee
Note  that the last two cases of \eqref{maxfoura}
are cuspidal and   have highest/lowest weight representations.

Not that ~$\cp^{4(4)}_1$~ and ~$\cp^{4(4)}_2$~ seem isomorphic, however,
they are distinguished   by the fact what roots
generate the ~$sl(3,\bbr)$ and  $sl(2,\bbr)$ subalgebras of ~$\cm^{4(4)}_{1,2}$~ - long or short.
Also the algebras  ~$\cn^{4(4)\pm}_1$~ and ~$\cn^{4(4)\pm}_2$~ differ by the number of long and short
roots, as indicated.
%Thus, we can write two versions of $\cm^{4(4)}_1$~:
%\eqnn{maxfouraz} \cm^{4(4)}_{1,1} ~&\cong&~  sl(3,\bbr)_S \oplus sl(2,\bbr)_L    \\
%\cm^{4(4)}_{1,2} ~&\cong&~ sl(3,\bbr)_L \oplus sl(2,\bbr)_S \nn\ee

%\np

%\np

\md

\subsection{FII :~   $F_{4(-20)}$}

\nt
The other real form of ~$F_4$~ is denoted as ~$F_{4(-20)}\,$ or ~$F''_4\,$.
The maximal compact subgroup ~$\ck \cong so(9)$, $\cp = \mathbf{16}$, ~$\dim\,\cp = 16$,
~$\dim\,\cn^\pm = 15$, ~$\cm_0  ~\cong~  so(7)$.
This real form has discrete series representations.

The Satake diagram is:
\eqn{satsour}
\black{{\a_1}} \riga\black{{\a_2}}
\Longrightarrow\black{{\a_3}} \riga\downcirc{{\a_4}} \eeq

The split rank is equal to 1, thus, the minimal and maximal parabolics coincide,  the
~$\cm$-factor  and ~$\cn$-factor are the same as in the Bruhat decomposition:
\eqnn{maxfourb}
\cp^{4(-20)} ~&=&~ \cm^{4(-20)}\ \oplus\ so(1,1)\ \oplus\ \cn^{\,4(-20)}\  ;\\
\cm^{4(-20)} ~&=&~  so(7)  \ ,
\qquad \dim\,\cn^{\,4(-20)\pm} ~=~ 15
\nn\ee

\subsection{G :~   $G_{2(2)}$}

\nt
The unique non-compact real form of ~$G_2$~ is the split form, denoted as ~$G_{2(2)}$ or ~$G'_2\,$.
The maximal compact subgroup ~$\ck \cong su(2) \oplus su(2)$, $\cp = \left( \mathbf{4},\mathbf{2}\right) $, where $\mathbf{4}=S^{3}%
\mathbf{2}$, ~$\dim\,\cp = 8$,
~$\dim\,\cn^\pm = \dim\,\ck = 6$.
This real form has discrete series representations.

The Dynkin-Satake diagram is:
\eqn{gsplit}
\downcirc{{\a_1}}\!\!\equiv\!\equiv\!\equiv\!\equiv\!\equiv\!\!>\!\!\downcirc{{\a_2}}
\eeq

The maximal parabolic subalgebras  are  given by (being $G_2$ non-simply laced, we denote the short/long nature of the roots):
\eqnn{maxge}
\cp^{2(2)}_{L\atop S} ~&=&~ \cm^{2(2)}_{L\atop S}\ \oplus\ so(1,1)\ \oplus\ \cn^{\,2(2)}_{L\atop S}\  ;\\
\cm^{2(2)}_L ~&=&~ sl(2,\bbr)_L\ ,  \qquad \dim\,\cn^{\,2(2)\pm}_L ~=~ 5 \nn \\
&& (2 ~{\rm long~roots}, 3 ~{\rm short~roots}), \nn\\
\cm^{2(2)}_S ~&=&~ sl(2,\bbr)_S\ ,  \qquad \dim\,\cn^{\,2(2)\pm}_S ~=~ 5 \nn \\
&& (3 ~{\rm long~roots}, 2 ~{\rm short~roots}), \nn
\ee
They are cuspidal   and have highest/lowest weight representations.

Not that ~$\cp^{2(2)}_L$~ and ~$\cp^{2(2)}_S$~ seem isomorphic, however,
they are distinguished   by the fact what root
generate the  ~$\cm^{2(2)} = sl(2,\bbr)$ subalgebras -- long or short.
Also the algebras  ~$\cn^{2(2)\pm}_L$~ and ~$\cn^{2(2)\pm}_S$~ differ by the number of long and short
roots, as indicated.

%\np

\subsection{Parabolically related non-compact semisimple Lie algebras}
\label{parab-rel}

Next, let us introduce the   notion of `parabolically related
non-compact semisimple Lie algebras' \cite{Dobparab} which
is also very useful in the study of the structure of real forms.

{\it Definition:}  ~~~Let ~$\cg,\cg'$~ be two non-compact
semisimple Lie algebras with the same complexification ~$\cg^\bac
\cong \cg'^\bac$. We call them ~{\it parabolically related}~ if they
have parabolic subalgebras ~$\cp = \cm \oplus \ca \oplus \cn$,
~$\cp' = \cm' \oplus \ca' \oplus \cn'$, such that: ~$\cm^\bac
~\cong~ \cm'^\bac$~  ($\Rightarrow \cp^\bac ~\cong~ \cp'^\bac$).

Certainly, there may be more than one parabolic relationship for an
algebra ~$\cg$. Furthermore, two algebras ~$\cg,\cg'$~ may be
parabolically related with different parabolic subalgebras. In the case the parabolic subalgebras determining the relation are maximal, the two corresponding algebras are said to be maximally parabolically related.

We summarize the maximally parabolically related exceptional Lie algebras
 in the following table, which slightly extends the results of \cite{Dobparab} :

\begin{center}

%\bigskip
\vbox{\offinterlineskip {\bf Table A} : maximally parabolically related non-compact real forms of finite-dimensional exceptional Lie
algebras; the corresponding interpretation in terms of symmetries of rank-2 and rank-3 Jordan algebras is given in {\bf Table B} \medskip \halign{\baselineskip12pt \strut
\vrule#\hskip0.1truecm & #\hfil&
\vrule#\hskip0.1truecm & #\hfil&
\vrule#\hskip0.1truecm & #\hfil&
\vrule#\hskip0.1truecm & #\hfil&
%\vrule#\hskip0.1truecm & #\hfil&
\vrule #\cr
\tablerule
&&& && &&&\cr
%%% HEADER ROW
&~ $\cg$ %&&~ $\ck'$
&&~  $\cm $ &&~$\cg'$&& ~$\cm'$&\cr
&&&&&&&&\cr
&&&~~$\dim \cn$&&&&&\cr
\tablerule

 &~1~:~$E_{6(-14)}$   %&& ~$so(10)$
 && ~$su(5,1)$ && ~$E_{6(6)}$
&&~$sl(6,\bbr)$&\cr
&&&&& &&&\cr
&&& ~~$21$
&&~$E_{6(2)}$&&~$su(3,3)$&\cr
%&&&&&&&&\cr
\tablerule

&~2~:~$E_{6(-14)}$
 && ~$so(7,1)\oplus so(2)$ && ~$E_{6(2)}$
&&~$so(5,3)\oplus u(1)$&\cr
&&&&& &&&\cr
&&& ~~$24$
&& && &\cr
%&&&&&&&&\cr
\tablerule

&~3~:~$E_{6(-26)}$
 && ~$so(9,1)$ && ~$E_{6(6)}$
&&~$so(5,5)$&\cr
&&&&& &&&\cr
&&& ~~$16$
&& && &\cr
%&&&&&&&&\cr
\tablerule

&~4~:~$E_{6(6)}$
 && ~$sl(3,\bbr) \oplus sl(3,\bbr)\oplus sl(2,\bbr)$ && ~$E_{6(2)}$
&&~$sl(3,\bbc)_\bbr \oplus sl(2,\bbr)$&\cr
&&&&& &&&\cr
&&& ~~$29$
&& && &\cr
%&&&&&&&&\cr
\tablerule

&~5~:~$E_{7(-25)}$    && ~$E_{6(-26)}$ && ~$E_{7(7)}$
&&~$E_{6(6)}$&\cr
&&&&&&&&\cr &&&~~$27$&&&&&\cr
\tablerule

&~6~:~$E_{7(-25)}$    && ~$so(9,1)\oplus sl(2,\bbr)$  && ~$E_{7(7)}$
&&~$so(5,5)\oplus sl(2,\bbr)$ &\cr
&&&&&&&&\cr
&&&~~$42$&&~$E_{7(-5)}$&& ~$so(7,3)\oplus su(2)$&\cr
\tablerule

&~7~:~$E_{7(-25)}$    && ~$so(10,2)$  && ~$E_{7(7)}$
&&~$so(6,6)$ &\cr
&&&&&&&&\cr
&&&~~$33$&&~$E_{7(-5)}$&& ~$so^*(12)$&\cr
\tablerule

&~8~:~$E_{7(7)}$
 && ~$sl(4,\bbr) \oplus sl(2,\bbr) \oplus sl(3,\bbr)$ && ~$E_{7(-5)}$
&&~$su^*(4) \oplus su(2) \oplus  sl(3,\bbr)$&\cr
&&&&& &&&\cr
&&& ~~$53$ && && &\cr
%&&&&&&&&\cr
\tablerule

&~9~:~$E_{7(-5)}$
 && ~$su^*(6) \oplus  sl(2,\bbr)$ && ~$E_{7(7)}$
&&~$sl(6,\bbr) \oplus  sl(2,\bbr) $&\cr
&&&&& &&&\cr
&&&~~$47$&&~$E_{7(-25)}$&& ~$su^*(6) \oplus su(2) $&\cr
%&&&&&&&&\cr
\tablerule

&~10~:~$E_{8(-24)}$      && ~$E_{7(-25)}$  && ~$E_{8(8)}$
&& ~$E_{7(7)}$  &\cr
&&&&& &&&\cr
&&& ~~$57$ && && &\cr
%&&&&&&&&\cr
\tablerule

&~11~:~$E_{8(-24)}$      && ~$so(11,3)$
&& ~$E_{8(8)}$      && ~$so(7,7)$   &\cr
&&&&& &&&\cr
&&& ~~$78$ && && &\cr
%&&&&&&&&\cr
\tablerule

&~12~:~$E_{8(-24)}$      && ~ $E_{6(-26)}\oplus sl(2,\bbr)$
&& ~$E_{8(8)}$      && ~ $E_{6(6)}\oplus sl(2,\bbr)$   &\cr
&&&&& &&&\cr
&&& ~~$83$ && && &\cr
%&&&&&&&&\cr
\tablerule

&~13~:~$E_{8(-24)}$      && ~ $so(9,1)\oplus sl(3,\bbr)$
&& ~$E_{8(8)}$      && ~ $so(5,5)\oplus sl(3,\bbr)$   &\cr
&&&&& &&&\cr
&&& ~~$97$ && && &\cr
%&&&&&&&&\cr
\tablerule

&~14~:~$F_{4(-20)}$      && ~ $so(7)$
&& ~$F_{4(4)}$      && ~ $so(4,3)$   &\cr
&&&&& &&&\cr
&&& ~~$15$ && && &\cr
%&&&&&&&&\cr
\tablerule
}}

\end{center}

\np

%\documentclass[11pt]{article}
%%%%%%%%%%%%%%%%%%%%%%%%%%%%%%%%%%%%%%%%%%%%%%%%%%%%%%%%%%%%%%%%%%%%%%%%%%%%%%%%%%%%%%%%%%%%%%%%%%%%%%%%%%%%%%%%%%%%%%%%%%%%%%%%%%%%%%%%%%%%%%%%%%%%%%%%%%%%%%%%%%%%%%%%%%%%%%%%%%%%%%%%%%%%%%%%%%%%%%%%%%%%%%%%%%%%%%%%%%%%%%%%%%%%%%%%%%%%%%%%%%%%%%%%%%%%
%\usepackage{amscd,array,dsfont,graphicx,mathtools}
%\usepackage[bookmarksnumbered,linktocpage]{hyperref}
%\usepackage[all]{hypcap}
%\usepackage{amsmath,epsfig,array}
%\usepackage{graphicx}
%\usepackage{amssymb}
%\usepackage{amsfonts,amstext,latexsym,xspace}
%\usepackage{amsfonts}
%\usepackage[margin=2cm]{geometry}

\setcounter{MaxMatrixCols}{10}

\newcolumntype{M}[1]{>{$}{#1}<{$}}
\newcommand{\half}{\tfrac{1}{2}}
\newcommand{\rep}[1]{\mathbf{#1}}
\newcommand{\cP}{\mathcal{P}}
\newcommand{\cQ}{\mathcal{Q}}
\newcommand{\cJ}{\mathcal{J}}
\newcommand{\N}{{\mathcal N}}
\newcommand{\R}{{\bf R}}
\newcommand{\eq}{\begin{equation}}
\newcommand{\en}{\end{equation}}
\newcommand{\enn}{\end{eqnarray}}
\newcommand{\M}{\ensuremath{\mathcal{M}}}
\newcommand{\CN}{\ensuremath{\mathcal{N}}}
\newcommand{\tM}{\ensuremath{\tilde{M}}}
\newcommand{\tb}{\ensuremath{\tilde{b}}}
\newcommand{\tc}{\ensuremath{\tilde{c}}}
\newcommand{\tx}{\ensuremath{\tilde{x}}}
\newcommand{\ty}{\ensuremath{\tilde{y}}}
\newcommand{\tz}{\ensuremath{\tilde{z}}}
\newcommand{\tu}{\ensuremath{\tilde{u}}}
\newcommand{\tw}{\ensuremath{\tilde{w}}}
\newcommand{\ti}{\ensuremath{\tilde{I}}}
\newcommand{\tj}{\ensuremath{\tilde{J}}}
\newcommand{\tk}{\ensuremath{\tilde{K}}}
\newcommand{\tm}{\ensuremath{\tilde{M}}}
\newcommand{\tq}{\ensuremath{\tilde{Q}}}
\newcommand{\tA}{\ensuremath{\tilde{A}}}
\newcommand{\tB}{\ensuremath{\tilde{B}}}
\newcommand{\tC}{\ensuremath{\tilde{C}}}
\newcommand{\tD}{\ensuremath{\tilde{D}}}
\newcommand{\tE}{\ensuremath{\tilde{E}}}
\newcommand{\tF}{\ensuremath{\tilde{F}}}
\newcommand{\tG}{\ensuremath{\tilde{G}}}
\newcommand{\JTP}[4][{}]{\{ #2\,#3\,#4 \}^{#1}}
\newcommand{\bi}{\ensuremath{\bar{i}}}
\newcommand{\bj}{\ensuremath{\bar{j}}}
\newcommand{\bk}{\ensuremath{\bar{k}}}
\newcommand{\hmu}{\ensuremath{\hat{\mu}}}
\newcommand{\hnu}{\ensuremath{\hat{\nu}}}
\newcommand{\hrho}{\ensuremath{\hat{\rho}}}
\newcommand{\hsigma}{\ensuremath{\hat{\sigma}}}
\newcommand{\hlambda}{\ensuremath{\hat{\lambda}}}
\newcommand{\he}{\ensuremath{{\hat{e}}}}
\newcommand{\hA}{\ensuremath{{\hat{A}}}}
\newcommand{\hF}{\ensuremath{{\hat{\F}}}}
\newcommand{\hcF}{\ensuremath{{\hat{\mathcal{F}}}}}
\newcommand{\hB}{\ensuremath{{\hat{B}}}}
\newcommand{\hcD}{\ensuremath{{\hat{\mathcal{D}}}}}
\newcommand{\tI}{\ensuremath{\tilde{I}}}
\newcommand{\tJ}{\ensuremath{\tilde{J}}}
\newcommand{\tK}{\ensuremath{\tilde{K}}}
\newcommand{\fg}{\ensuremath{\mathfrak{g}}}
\newcommand  {\Rbar} {{\mbox{\rm$\mbox{I}\!\mbox{R}$}}}
\newcommand  {\Hbar} {{\mbox{\rm$\mbox{I}\!\mbox{H}$}}}
\newcommand {\Cbar}{\mathord{\setlength{\unitlength}{1em}
     \begin{picture}(0.6,0.7)(-0.1,0) \put(-0.1,0){\rm C}
        \thicklines \put(0.2,0.05){\line(0,1){0.55}}\end {picture}}}
\newcommand{\cD}{\ensuremath{\mathcal{D}}\xspace}
\newcommand{\Gtwo}[1][{}]{\ensuremath{G_{2 #1}}\xspace}
\newcommand{\4}[1][{}]{\ensuremath{F_{4 #1}}\xspace}
\newcommand{\6}[1][{}]{\ensuremath{E_{6 #1}}\xspace}
\newcommand{\7}[1][{}]{\ensuremath{E_{7 #1}}\xspace}
\newcommand{\8}[1][{}]{\ensuremath{E_{8 #1}}\xspace}
\newcommand{\E}{E_{8(8)}}
\newcommand{\EE}{E_{7(7)}}
\newcommand{\EEE}{E_{6(6)}}
\newcommand{\9}{\ensuremath{E_9}\xspace}
\newcommand{\0}{\ensuremath{E_{10}}\xspace}
\newcommand{\USp}[1][8]{\ensuremath{\mbox{USp$(#1)$}}\xspace}
\newcommand{\SU}[1][8]{\ensuremath{\mbox{SU$(#1)$}}\xspace}
\newcommand{\SO}[1][8]{\ensuremath{\mbox{SO$(#1)$}}\xspace}
\newcommand{\SL}[1][8]{\ensuremath{\mbox{SL$(#1)$}}\xspace}
\newcommand{\SLR}[1][8]{\ensuremath{\mbox{SL$(#1,\mathbb{R})$}}\xspace}
\newcommand{\SUS}[1][6]{\ensuremath{\mbox{SU${}^\ast(#1)$}}\xspace}
\newcommand{\SOS}[1][12]{\ensuremath{\mbox{SO${}^\ast(#1)$}}\xspace}

%\begin{document}

\setcounter{equation}{0}
\section{\label{sec:JA}Jordan Algebras, Freudenthal Triple Systems, Their Symmetries and
Embeddings}

\subsection{\label{sec:J}Jordan Algebras}

A Jordan algebra $\mathfrak{J}$ \cite%
{Jordan:1933a,Jordan:1933b,Jordan:1933vh,Jacobson:1961,Jacobson:1968} is
vector space defined over a ground field $\mathbb{F}$ equipped with a
bilinear product (named Jordan product) satisfying
\begin{equation}
\begin{split}
X\circ Y& =Y\circ X, \\
X^{2}\circ (X\circ Y)& =X\circ (X^{2}\circ Y),\quad \forall \ X,Y\in
\mathfrak{J}.
\end{split}
\label{eq:Jid}
\end{equation}%
For the treatment given in the present investigation, the relevant Jordan
algebras are examples of the class of \textit{cubic} Jordan algebras over $%
\mathbb{F}=\mathbb{R}$ \cite{Elkies:1996,Gross:1996,
Krutelevich:2002,Krutelevich:2004}. A cubic Jordan algebra is endowed with a
cubic form $N:\mathfrak{J}\rightarrow \mathbb{R}$, such that $N(\lambda
X)=\lambda ^{3}N(X),\quad \forall \ \lambda \in \mathbb{R},\ X\in \mathfrak{J%
}$. Moreover, an element $c\in \mathfrak{J}$ exists, satisfying $N(c)=1$
(usually named \textit{base point}). A general procedure for constructing
cubic Jordan algebras, due to Freudenthal, Springer and Tits \cite%
{Springer:1962, McCrimmon:1969, McCrimmon:2004}, exists, in which all
properties of the Jordan algebra are determined by the cubic form itself.

Let $V$ be a vector space equipped with a cubic norm $N:V\rightarrow \mathbb{%
R}$ such that $N(\lambda X)=\lambda ^{3}N(X),\ \forall \ \lambda \in \mathbb{%
R},\ X\in V$, and with a base point $c\in V$ satisfying $N(c)=1$. Then, if
the full linearization of the cubic norm, denoted by $N(X,Y,Z)$ and defined
as%
\begin{equation}
6N(X,Y,Z)\doteq N\left( X+Y+Z\right) -N(X+Y)-N(Y+Z)-N(X+Z)+N(X)+N(Y)+N(Z),
\end{equation}%
is trilinear, the following four maps can be introduced :

\begin{enumerate}
\item The trace
\begin{equation}
Tr(X)\doteq 3N(c,c,X);  \label{eq:cubicdefs}
\end{equation}

\item A quadratic map
\begin{equation}
S(X)\doteq 3N(X,X,c),
\end{equation}

\item A bilinear map
\begin{equation}
S(X,Y)\doteq 6N(X,Y,c),
\end{equation}

\item A trace bilinear form
\begin{equation}
Tr(X,Y)\doteq Tr(X)Tr(Y)-S(X,Y).  \label{eq:tracebilinearform}
\end{equation}
\end{enumerate}

A cubic Jordan algebra $\mathfrak{J}$ with multiplicative identity $Id=c$
can be obtained starting from the vector space $V$ above \textit{iff} $N$ is
\textit{Jordan cubic}, namely \textit{iff}: \textbf{I]} The trace bilinear
form \eqref{eq:tracebilinearform} is non-degenerate, and \textbf{II]} the
quadratic adjoint map, $\sharp \colon \mathfrak{J}\rightarrow \mathfrak{J}$,
uniquely defined by $Tr(X^{\sharp },Y):=3N(X,X,Y)$, satisfies

\begin{equation}
(X^{\sharp })^{\sharp }=N(X)X,\quad \forall X\in \mathfrak{J}.
\label{eq:Jcubic}
\end{equation}

In a cubic Jordan algebra, the so-called \textit{Jordan product} can be
introduced in the following way :
\begin{equation}
X\circ Y\doteq \tfrac{1}{2}\left( X\times Y+Tr(X)Y+Tr(Y)X-S(X,Y)Id\right) ,
\end{equation}%
where $X\times Y$ denotes the linearization of the quadratic adjoint map :
\begin{equation}
X\times Y\doteq (X+Y)^{\sharp }-X^{\sharp }-Y^{\sharp }.
\label{eq:FreuProduct}
\end{equation}%
Another related map is the \textit{Jordan triple product} :
\begin{equation}
\{X,Y,Z\}\doteq (X\circ Y)\circ Z+X\circ (Y\circ Z)-(X\circ Z)\circ Y.
\label{eq:Jtripleproduct}
\end{equation}

Jordan algebras were introduced and completely classified in \cite%
{Jordan:1933vh} in an attempt to generalize quantum mechanics beyond the
complex numbers $\mathbb{C}$. Below, we list all allowed possibilities of
cubic Jordan algebras \cite%
{Jacobson:1961,Jacobson:1968,Krutelevich:2004,McCrimmon:1969,Baez:2001dm}:

\begin{enumerate}
\item the simplest case: $\mathfrak{J}$ $=\mathbb{R}$, $N(X)\doteq X^{3}$;

\item the infinite sequence of semi-simple Jordan algebras given by $%
\mathfrak{J}=\mathbb{R}\oplus \Gamma_{m,n}$ (named \textit{pseudo-Euclidean}
\textit{spin factors}), where $\Gamma _{m,n}$ is an $\left( m+n\right) $%
-dimensional vector space over $\mathbb{R}$, namely the Clifford algebra of $%
O(m,n)$, with $N(X=\xi \oplus \gamma )\doteq \xi \gamma ^{a}\gamma ^{b}\eta
_{ab}$;

\item Four exceptional, simple and Euclidean\footnote{%
The Lorentzian version of the exceptional, simple, Lorentzian cubic Jordan
algebras can also be defined; for its definition and symmetries, see \cite%
{Squaring-Magic}, and for its use in supergravity see \cite%
{GZ}.} cases, given by $\mathfrak{J}=J_{3}^{\mathbb{A}}$ or $%
\mathfrak{J}=J_{3}^{\mathbb{A}_{s}}$, the algebra of $3\times 3$ Hermitian
matrices over the four division algebras $\mathbb{A=R}$ (real numbers)$,%
\mathbb{C}$ (complex numbers)$,\mathbb{H}$ (quaternions)$,\mathbb{O}$
(octonions), or their split versions $\mathbb{A}_{s}=\mathbb{C}_{s},\mathbb{H%
}_{s},\mathbb{O}_{s}$ :
\begin{equation}
X=\left(
\begin{array}{ccc}
\alpha _{1} & x_{3} & \overline{x}_{2} \\
\overline{x}_{3} & \alpha _{2} & x_{1} \\
x_{2} & \overline{x}_{1} & \alpha _{3}%
\end{array}%
\right) ,~\alpha _{1},\alpha _{2},\alpha _{3}\in \mathbb{R}%
,~x_{1},x_{2},x_{3}\in \mathbb{A}\text{ (or }\mathbb{A}_{s}\text{)},
\end{equation}%
with conjugation (denoted by bar) pertaining to the relevant (division or
split) algebra. In these cases, the cubic norm is given by%
\begin{equation}
N\left( X\right) \doteq \alpha _{1}\alpha _{2}\alpha _{3}-\alpha _{1}x_{1}%
\overline{x}_{1}-\alpha _{2}x_{2}\overline{x}_{2}-\alpha _{3}x_{3}\overline{x%
}_{3}+2\text{Re}\left( x_{1}x_{2}x_{3}\right) .
\end{equation}%
This reproduces the usual determinant\footnote{%
For explicit constructions of $N(X)$, see \textit{e.g.} \cite{F-Gimon-K} and
\cite{Wissanji}.} for $\mathbb{A}=\mathbb{R\ }$\ and $\mathbb{C}$. In these
cases, the Jordan product simply reads%
\begin{equation}
X\circ Y\doteq \tfrac{1}{2}(XY+YX),\label{J-product}
\end{equation}%
where $XY$ is just the conventional $3\times 3$ matrix product; see \textit{%
e.g.} \cite{Jacobson:1968} for a comprehensive account.
\end{enumerate}

\md

\subsection{\label{Jordan-Pairs}Jordan Pairs}

Jordan algebras have traveled a long journey, since their appearance in the
30's \cite{Jordan:1933vh}. The modern formulation \cite{jaco} involves a quadratic
map $U_{x}y$ (like $xyx$ for associative algebras) instead of the original
symmetric product (\ref{J-product}). The quadratic map and its linearization
$V_{x,y}z=(U_{x+z}-U_{x}-U_{z})y$ (like $xyz+zyx$ in the associative case)
reveal the mathematical structure of Jordan algebras much more clearly,
through the notion of inverse, inner ideal, generic norm, etc. The axioms
are:
\begin{equation}
U_{1}=Id\quad ,\qquad U_{x}V_{y,x}=V_{x,y}U_{x}\quad ,\qquad
U_{U_{x}y}=U_{x}U_{y}U_{x}  \label{qja}
\end{equation}%
The quadratic formulation led to the concept of \textit{Jordan triple systems%
} \cite{myb}, an example of which is a pair of modules represented by
rectangular matrices. There is no way of multiplying two matrices $x$ and $y$
, say $n\times m$ and $m\times n$ respectively, by means of a bilinear
product. But one can do it using a product like $xyx$, quadratic in $x$ and
linear in $y$. Notice that, like in the case of rectangular matrices, there
needs not be a unity in these structures. The axioms are in this case:
\begin{equation}
U_{x}V_{y,x}=V_{x,y}U_{x}\quad ,\qquad V_{U_{x}y,y}=V_{x,U_{y}x}\quad
,\qquad U_{U_{x}y}=U_{x}U_{y}U_{x}  \label{jts}
\end{equation}

Finally, a \textit{Jordan pair} \cite{loos1} is just a pair of modules $%
(V^{+},V^{-})$ acting on each other (but not on themselves) like a Jordan
triple:
\begin{equation}
\begin{array}{ll}
U_{x^{\sigma }}V_{y^{-\sigma },x^{\sigma }} & =V_{x^{\sigma },y^{-\sigma
}}U_{x^{\sigma }} \\
V_{U_{x^{\sigma }}y^{-\sigma },y^{-\sigma }} & =V_{x^{\sigma },U_{y^{-\sigma
}}x^{\sigma }} \\
U_{U_{x^{\sigma }}y^{-\sigma }} & =U_{x^{\sigma }}U_{y^{-\sigma
}}U_{x^{\sigma }}%
\end{array}
\label{jp}
\end{equation}%
where $\sigma =\pm $ and $x^{\sigma }\in V^{+\sigma }\,,\;y^{-\sigma }\in
V^{-\sigma }$.

\textit{Jordan pairs} \cite{loos1} $(V^{+},V^{-})\doteq \left( \mathfrak{J},%
\mathfrak{J}^{\prime }\right) $ (whose recent mathematical and physical
treatment can be found \textit{e.g.} in \cite{Truini-1,FMZ,Truini-2}; also
\textit{cfr.} \cite{McCrimmon:2004} for a review) are strictly related to
the Tits-Kantor-Koecher construction of 3-graded Lie Algebras $\mathfrak{L}$
\cite{tits1,kantor1,koecher1} (see also the interesting relation to Hopf
algebras \cite{faulk}):
\begin{equation}
\mathfrak{L}=\mathfrak{J}^{\prime }\oplus str(\mathfrak{J})\oplus \mathfrak{J%
},  \label{tkk}
\end{equation}%
where $\mathfrak{J}$ is a Jordan algebra, $str(\mathfrak{J})=\mathcal{L}(%
\mathfrak{J})\oplus der(\mathfrak{J})$ is the structure Lie algebra of $%
\mathfrak{J}$ \cite{McCrimmon:2004}, $\mathcal{L}\left( \mathfrak{J}\right) $
is the left multiplication in $\mathfrak{J}$, and $der(\mathfrak{J})=[%
\mathcal{L}(\mathfrak{J}),\mathcal{L}(\mathfrak{J})]$ is the algebra of
derivations of $\mathfrak{J}$ (\textit{i.e.}, the algebra of the
automorphism group of $\mathfrak{J}$) \cite{Schafer:1966,Schafer-2}. As we
will see in the next subsection (cfr. (\ref{conf})), $\mathfrak{L}$ can be
identified with the conformal Lie algebra $conf(\mathfrak{J})$ of $\mathfrak{%
J}$ itself.

\md

\subsection{\label{Symms-JAs}Symmetries of Jordan Algebras}

To each cubic Jordan algebra, a number of symmetry groups can be associated :

\begin{itemize}
\item $Aut(\mathfrak{J})$, the group of automorphisms of $\mathfrak{J}$,
which leaves invariant the structure constants of the Jordan product (the
Lie algebra of $Aut(\mathfrak{J})$ is given by the derivations $der(%
\mathfrak{J})$ of $\mathfrak{J}$).

\item $Str(\mathfrak{J})$, the \textit{structure} group, with Lie algebra $%
str(\mathfrak{J})$, which leaves the cubic norm $N$ \textit{invariant} up to
a rescaling:%
\begin{equation}
N(\mathbf{g}(X))=\lambda N(X),\text{ }\lambda \in \mathbb{R},\quad \forall \
\mathbf{g}\in Str(\mathfrak{J});
\end{equation}%
the \textit{reduced structure} group $Str_{0}(\mathfrak{J})$, with Lie
algebra $str_{0}(\mathfrak{J})$, is obtained from $Str(\mathfrak{J})$ by
modding it out by its center \cite{Schafer:1966, Jacobson:1968, Brown:1969}:%
\begin{equation}
N(\mathbf{g}(X))=N(X),\quad \forall \ \mathbf{g}\in Str_{0}(\mathfrak{J}).
\end{equation}%
It should be here remarked that the structure group of $\mathfrak{J}$ is the
automorphism group of the Jordan pair $(\mathfrak{J},\mathfrak{J}^{\prime })$
:%
\begin{equation}
Str(\mathfrak{J})\cong Aut\left( \mathfrak{J},\mathfrak{J}^{\prime }\right) .
\end{equation}

\item $Conf(\mathfrak{J})$, the \textit{conformal} group, whose Lie algebra $%
conf(\mathfrak{J})$ can be given a 3-graded structure with respect to $str(%
\mathfrak{J})$ :%
\begin{equation}
conf(\mathfrak{J})=g^{-1}\oplus str\left( \mathfrak{J}\right) \oplus g^{1}.%
\label{conf}
\end{equation}%
The Tits-Kantor-Koecher construction \cite{tits1,kantor1,koecher1} (\ref{tkk}) of $conf(%
\mathfrak{J})\equiv \mathfrak{L}$ establishes a one-to-one mapping between
the grade +1 subspace $g^{1}$ of $conf(\mathfrak{J})$ and the corresponding
Jordan algebra $\mathfrak{J}$ : $g^{1}\Leftrightarrow \mathfrak{J}$. Every
Lie algebra $\mathfrak{L}$ (\ref{tkk}) admits a conjugation (involutive
automorphism) under which the elements of the grade +1 subspace get mapped
into the elements of the grade -1 subspace, and vice versa (grade-reversing
nature of the involution); cfr. e.g. \cite{Gun-CQC}. Remarkably, $Conf(%
\mathfrak{J})$ is isomorphic to the automorphism group of the (reduced)
Freudenthal triple system defined over $\mathfrak{J}$; see below.

\item $QConf(\mathfrak{J})$, the \textit{quasi-conformal} group, with Lie
algebra $qconf(\mathfrak{J})$, which can be defined by introducing
Freudenthal triple systems and their further extension named \textit{extended%
} Freudenthal triple system \cite{Gunaydin:2000xr}; see below.
\end{itemize}

All symmetry groups of (simple and semi-simple) cubic Jordan algebras over $%
\mathbb{R}$ are listed\footnote{%
Besides cubic Jordan algebras and their symmetries, in the subsequent
treatment we will also consider another remarkable \textit{Hermitian Jordan
triple system}, namely given by $2$-dimensional octonionic vectors, and
denoted by $M_{2,1}\left( \mathbb{O}\right) $. Its relevant symmetries are $%
conf\left( M_{2,1}\left( \mathbb{O}\right) \right) \cong E_{6(-14)}$, and $%
str_{0}\left( M_{2,1}\left( \mathbb{O}\right) \right) \cong so(8,2)$ (%
\textit{cfr.} \cite{Koecher,Loos,Gun-Bars, Gun-CQC} and \cite{GST}).} in
Table 1.

Remarkably, the symmetries of cubic Jordan algebras over $J_{3}^{\mathbb{A}}$
and $J_{3}^{\mathbb{A}_{s}}$ respectively arrange as entries of the
single-split and doubly-split Magic Squares of order three (reported in
Tables 2 and 3), which are non-compact, real forms of the
Freudenthal-Rozenfeld-Tits Magic Square of Lie algebras itself \cite{FRT}
(reported in Table 4); generally, $aut$, $str_{0}$, $conf$ and $qconf$ Lie
algebras enter the first, second, third and fourth rows of the $4\times 4$
array of algebras constituting the Magic Square (\textit{cfr.} \cite%
{Squaring-Magic} for a comprehensive review).

\begin{table}[tbp]
\begin{equation*}
\begin{array}{|c||c|c|c|c|}
\hline
\mathfrak{J} & aut(\mathfrak{J}) & str_{0}(\mathfrak{J}) & conf(\mathfrak{J})
& qconf(\mathfrak{J}) \\ \hline\hline
\mathbb{R} & \varnothing & \varnothing & sl(2,\mathbb{R}) & G_{2(2)} \\
\hline
\mathbb{R\oplus }\Gamma _{m,n} & so(m)\oplus so(n) & so(m,n) & sl(2,\mathbb{R%
})\oplus so(m+1,n+1) & so(m+3,n+3) \\ \hline
J_{3}^{\mathbb{R}} & so(3) & sl(3,\mathbb{R}) & sp(3,\mathbb{R}) & F_{4(4)}
\\ \hline
J_{3}^{\mathbb{C}} & su(3) & sl(3,\mathbb{C})_{\mathbb{R}} & su(3,3) &
E_{6(2)} \\ \hline
J_{3}^{\mathbb{C}_{s}} & sl(3,\mathbb{R}) & sl(3,\mathbb{R})\oplus sl(3,%
\mathbb{R}) & sl(6,\mathbb{R}) & E_{6(6)} \\ \hline
J_{3}^{\mathbb{H}} & usp(3) & su^{\ast }(6) & so^{\ast }(12) & E_{7(-5)} \\
\hline
J_{3}^{\mathbb{H}_{s}} & sp(3,\mathbb{R}) & sl(6,\mathbb{R}) & so(6,6) &
E_{7(7)} \\ \hline
J_{3}^{\mathbb{O}} & F_{4(-52)} & E_{6(-26)} & E_{7(-25)} & E_{8(-24)} \\
\hline
J_{3}^{\mathbb{O}_{s}} & F_{4(4)} & E_{6(6)} & E_{7(7)} & E_{8(8)} \\ \hline
\end{array}%
\end{equation*}%
\caption{Lie algebras associated to cubic (i.e., rank-3) Euclidean Jordan
algebras. The notation $g(\mathbb{C})_{\mathbb{R}}$ means the algebra $g(%
\mathbb{C})$ seen as a real algebra}
\end{table}

\begin{table}[h]
\begin{center}
\begin{tabular}{|c|c|c|c|c|}
\hline
& $\mathbb{R}$ & $\mathbb{C}$ & $\mathbb{H}$ & $\mathbb{O}$ \\ \hline
$\mathbb{R}$ & $so(3)$ & $su(3)$ & $usp(3)$ & $F_{4(-52)}$ \\ \hline
$\mathbb{C}_{s}$ & $sl(3,\mathbb{R})$ & $sl(3,\mathbb{C})_{\mathbb{R}}$ & $%
su^{\ast }(6)$ & $E_{6(-26)}$ \\ \hline
$\mathbb{H}_{s}$ & $sp(3,\mathbb{R})$ & $su(3,3)$ & $so^{\ast }(12)$ & $%
E_{7(-25)}$ \\ \hline
$\mathbb{O}_{s}$ & $F_{4(4)}$ & $E_{6(2)}$ & $E_{7(-5)}$ & $E_{8(-24)}$ \\
\hline
\end{tabular}%
\end{center}
\caption{The \textit{single-split} (non-symmetric) real form of the Magic
Square ${\mathcal{L}}_{3}(\mathbb{A}_{s},\mathbb{B})$ \protect\cite{GST}}
\label{magicgst}
\end{table}

\begin{table}[h]
\begin{center}
\begin{tabular}{|c|c|c|c|c|}
\hline
& $\mathbb{R}$ & $\mathbb{C}_{s}$ & $\mathbb{H}_{s}$ & $\mathbb{O}_{s}$ \\
\hline
$\mathbb{R}$ & $so(3)$ & $sl(3,\mathbb{R})$ & $sp(3,\mathbb{R})$ & $F_{4(4)}$
\\ \hline
$\mathbb{C}_{s}$ & $sl(3,\mathbb{R})$ & $sl(3,\mathbb{R})\oplus sl(3,\mathbb{%
R})$ & $sl(6,\mathbb{R})$ & $E_{6(6)}$ \\ \hline
$\mathbb{H}_{s}$ & $sp(3,\mathbb{R})$ & $sl(6,\mathbb{R})$ & $so(6,6)$ & $%
E_{7(7)}$ \\ \hline
$\mathbb{O}_{s}$ & $F_{4(4)}$ & $E_{6(6)}$ & $E_{7(7)}$ & $E_{8(8)}$ \\
\hline
\end{tabular}%
\end{center}
\caption{The \textit{double-split} (symmetric) real form of the Magic Square
${\mathcal{L}}_{3}(\mathbb{A}_{s},\mathbb{B}_{s})$ \protect\cite{BS}}
\label{magicBS}
\end{table}

\begin{table}[h]
\begin{center}
\begin{tabular}{|c|c|c|c|c|}
\hline
& $\mathbb{R}$ & $\mathbb{C}$ & $\mathbb{H}$ & $\mathbb{O}$ \\ \hline
$\mathbb{R}$ & $so(3)$ & $su(3)$ & $usp(3)$ & $F_{4(-52)}$ \\ \hline
$\mathbb{C}$ & $su(3)$ & $su(3)\oplus su(3)$ & $su(6)$ & $E_{6(-78)}$ \\
\hline
$\mathbb{H}$ & $usp(3)$ & $su(6)$ & $so(12)$ & $E_{7(-133)}$ \\ \hline
$\mathbb{O}$ & $F_{4(-52)}$ & $E_{6(-78)}$ & $E_{7(-133)}$ & $E_{8(-248)}$
\\ \hline
\end{tabular}%
\end{center}
\caption{{}The compact, real form of the \textit{Freudenthal-Rozenfeld-Tits}
symmetric Magic Square ${\mathcal{L}}_{3}(\mathbb{A},\mathbb{B})$
\protect\cite{FRT}}
\label{magicfrt}
\end{table}

\md

\subsection{\label{sec:F}Freudenthal Triple Systems}

Starting from a cubic Jordan algebra $\mathfrak{J}$, a (reduced \cite{Brown:1969}) \textit{Freudenthal
triple system} (FTS) is defined as the vector space
\begin{equation}
\mathbf{F}\mathfrak{(J)}\doteq \mathbb{R}\oplus \mathbb{R}\oplus \mathfrak{%
J\oplus J}.
\end{equation}%
An element $\mathbf{x}\in \mathbf{F}\mathfrak{(J)}$ can thus formally be
written as a \textquotedblleft $2\times 2$ matrix\textquotedblright\ :
\begin{equation}
\mathbf{x}=%
\begin{pmatrix}
x & X \\
Y & y%
\end{pmatrix}%
,\text{\ }x,y\in \mathbb{R},~X,Y\in \mathfrak{J}.
\end{equation}%
An FTS is endowed\footnote{%
It is worth remarking that all the other necessary definitions, such as the
cubic and trace bilinear forms, are inherited from the underlying Jordan
algebra $\mathfrak{J}$.} with a non-degenerate bilinear antisymmetric
quadratic form, a quartic form and a trilinear triple product \cite%
{Freudenthal:1954,Brown:1969,Faulkner:1971, Ferrar:1972, Krutelevich:2004} :

\begin{enumerate}
\item Quadratic form $\{\bullet ,\mathbf{\bullet }\}$: $\mathbf{F}\mathfrak{%
(J)}\times \mathbf{F}\mathfrak{(J)}\rightarrow \mathbb{R}$, defined as
\begin{subequations}
\begin{equation}
\{\mathbf{x},\mathbf{y}\}\doteq \alpha \delta -\beta \gamma +Tr\left(
A,D\right) -Tr\left( B,C\right) ,  \label{eq:bilinearform}
\end{equation}%
where
\end{subequations}
\begin{equation}
\mathbf{x}=%
\begin{pmatrix}
\alpha  & A \\
B & \beta
\end{pmatrix}%
,~\mathbf{y}=%
\begin{pmatrix}
\gamma  & C \\
D & \delta
\end{pmatrix}%
.
\end{equation}

\item Quartic form $\Delta :\mathbf{F}\mathfrak{(J)}\rightarrow \mathbb{R}$,
defined as
\begin{subequations}
\begin{equation}
\Delta (\mathbf{x})\doteq -4\left( \alpha N(A)+\beta N(B)+\kappa (\mathbf{x}%
)^{2}-Tr(A^{\sharp },B^{\sharp })\right) ,  \label{eq:quarticnorm}
\end{equation}%
where
\end{subequations}
\begin{equation}
\kappa (\mathbf{x})\doteq \tfrac{1}{2}(\alpha \beta -Tr(A,B)).
\end{equation}

\item Triple product $T:\mathbf{F}\mathfrak{(J)}\times \mathbf{F}\mathfrak{%
(J)}\times \mathbf{F}\mathfrak{(J)}\rightarrow \mathbf{F}\mathfrak{(J)}$,
defined as

\begin{equation}
\{T(\mathbf{x},\mathbf{y},\mathbf{w}),\mathbf{z}\}=2\Delta (\mathbf{x},%
\mathbf{y},\mathbf{w},\mathbf{z}),
\end{equation}%
where $\Delta (\mathbf{x},\mathbf{y},\mathbf{w},\mathbf{z})$ is the full
linearization of $\Delta (\mathbf{x})$, such that $\Delta (\mathbf{x},%
\mathbf{x},\mathbf{x},\mathbf{x})=\Delta (\mathbf{x})$.
\end{enumerate}

The \textit{automorphism} group $Aut(\mathbf{F}\mathfrak{(J)})$, with Lie
algebra $aut(\mathbf{F}\mathfrak{(J)})$, is defined as the set of all
invertible $\mathbb{R}$-linear transformations which leave both $\{\mathbf{x}%
,\mathbf{y}\}$ and $\Delta (\mathbf{x},\mathbf{y},\mathbf{w},\mathbf{z})$
invariant \cite{Brown:1969}.

It can be proved \cite{Gunaydin:1975mp,Gunaydin:1989dq,Gunaydin:1992zh}
that, as anticipated above :

\begin{equation}
Aut(\mathbf{F}\mathfrak{(J)})\cong Conf\left( \mathfrak{J}\right) .
\end{equation}

\md

\subsection{\label{EFTS}Extended Freudenthal Triple Systems}

Every simple Lie algebra $g$ can be endowed with a $5$-grading, determined
by one of its generators $\mathcal{G}\equiv so(1,1)$, with one-dimensional $%
\pm 2$-graded subspaces :
\begin{equation}
g=g^{-2}\oplus g^{-1}\oplus g^{0}\oplus g^{+1}\oplus g^{+2}\,.
\end{equation}%
where
\begin{eqnarray}
g^{0} &=&conf\left( \mathfrak{J}\right) \oplus \mathcal{G}; \\
\lbrack \mathcal{G},\mathfrak{t}] &=&m\mathfrak{t}\;\;\;\;\forall \mathfrak{t%
}\in g^{m}\;\;,\;m=0,\pm 1,\pm 2
\end{eqnarray}%
As firstly discussed in \cite{Gunaydin:2000xr}, a $5$-graded\footnote{%
For $sl(2)$, the $5$-grading degenerates into a $3$-grading.} Lie algebra $g$
can geometrically be constructed as the \textit{quasi-conformal} Lie algebra
$qconf\left( \mathfrak{J}\right) $ over a vector space $\mathbf{EF}\left(
\mathfrak{J}\right) $, named \textit{extended Freudenthal triple system} (EFTS), which is coordinatized by $\mathcal{X}:=\left( \mathbf{x},\Phi
\right) \in \mathbf{EF}\left( \mathfrak{J}\right) $, where $\mathbf{x}\in
\mathbf{F}\left( \mathfrak{J}\right) $, and $\Phi $ is an extra real
variable \cite{Gunaydin:2000xr,Gunaydin:2005zz}:%
\begin{equation*}
\mathbf{EF}\left( \mathfrak{J}\right) :=\mathbf{F}\left( \mathfrak{J}\right)
\oplus \mathbb{R}.
\end{equation*}

Remarkably, a norm $\mathcal{N}:\mathbf{EF}\mathfrak{(J)}\rightarrow \mathbb{%
R}$ can be defined by using the quartic form $\Delta $ previously introduced
in $\mathbf{F}\left( \mathfrak{J}\right) $, as follows\footnote{%
Since the image of $\mathcal{\Delta }$ in $\mathbf{F}\left( \mathfrak{J}%
\right) $ extends over the whole $\mathbb{R}$, for $\mathcal{\Delta }(%
\mathbf{x})<0$ the light-like condition $\mathcal{N}(\mathcal{X})=0$ in $%
\mathbf{EF}\left( \mathfrak{J}\right) $ does not yield real solutions for $%
\Phi $. However, as discussed in \cite{Gunaydin:2000xr}, this problem can be
solved by complexifying the whole $\mathbf{EF}\left( \mathfrak{J}\right) $ (%
\textit{i.e.}, by considering $\mathbb{F}=\mathbb{C}$ as ground field), thus
obtaining a realization of the \textit{complexified} Lie algebra $\mathfrak{g%
}\left( \mathbb{C}\right) $ over $\left[ \mathbf{EF}\left( \mathfrak{J}%
\right) \right] _{\mathbb{C}}$.} :
\begin{equation}
\mathcal{N}(\mathcal{X}):=\mathcal{\Delta }(\mathbf{x})-\Phi ^{2}.
\label{N-call}
\end{equation}%
Furthermore, a \textquotedblleft quartic distance" $d_{4}:\mathbf{EF}%
\mathfrak{(J)\times }\mathbf{EF}\mathfrak{(J)}\rightarrow \mathbb{R}$
between any two points $\mathcal{X}=\left( \mathbf{x},\Phi \right) $ and $%
\mathcal{Y}:=\left( \mathbf{y},\Psi \right) $ in $\mathbf{EF}\left(
\mathfrak{J}\right) $ can be defined as
\begin{equation}
d_{4}(\mathcal{X},\mathcal{Y}):=\mathcal{\Delta }(\mathbf{x}-\mathbf{y}%
)-\left( \Phi -\Psi +\left\{ \mathbf{x},\mathbf{y}\right\} \right) ^{2},
\label{d4}
\end{equation}%
such that $\mathcal{N}(\mathcal{X})=d_{4}(\mathcal{X},\mathcal{Y}=0)$.

Then, the \textit{quasi-conformal} group $QConf\left( \mathfrak{J}\right) $
over $\mathbf{EF}\left( \mathfrak{J}\right) $, with Lie algebra $qconf\left(
\mathfrak{J}\right) $, is defined as the set of all invertible $\mathbb{R}$%
-linear transformations which leave invariant the \textquotedblleft quartic
light-cone", namely, the geometrical locus defined by \cite{Gunaydin:2000xr}
\begin{equation}
d_{4}(\mathcal{X},\mathcal{Y})=0,~\forall \left( \mathcal{X},\mathcal{Y}%
\right) \in \left( \mathbf{EF}\left( \mathfrak{J}\right) \right) ^{2},
\end{equation}

Thus, every $5$-graded Lie algebra $g$, geometrically realized as the
\textit{quasi-conformal} Lie algebra $qconf\left( \mathfrak{J}\right) $ over
a vector space $\mathbf{EF}\left( \mathfrak{J}\right) $, admits a \textit{%
conformal invariant} given by the norm $\mathcal{N}$ (\ref{N-call}); see
also \cite{Gunaydin:2005zz}.

\md

\subsection{\label{Embs}Embeddings}

In this Subsection we will present some general structures of embeddings
involving the symmetry algebras of cubic Jordan algebras introduced above,
along with some comments on their physical meaning and relevance.

We start and remark that $str_{0}$, $conf$ and $qconf$ Lie algebras can be
interpreted as global, electric-magnetic duality ($U$-duality\footnote{%
Here $U$-duality is referred to as the \textquotedblleft continuous"
symmetries of \cite{j}. Their discrete versions are the $U$-duality
non-perturbative string theory symmetries introduced in \cite{HT}.})
symmetries of suitable theories of gravity (possibly with local
supersymmetry) coupled to scalar fields and Abelian vectors, respectively in
$D=5,4,3$ (Lorenzian-signed) space-time dimensions (\textit{cfr.} \textit{%
e.g.} \cite{Gun-2}, and Refs. therein). Such symmetries are non-linearly
realized on the scalars, while vectors do sit in some linear representations
of them (in $D=3$, they are completely dualized into scalar fields, as well).

As mentioned above, $aut$, $str_{0}$, $conf$ and $qconf$ Lie algebras of
cubic (i.e. rank-3) Jordan algebras over $J_{3}^{\mathbb{A}}$ and $J_{3}^{%
\mathbb{A}_{s}}$ respectively enter the first, second, third and fourth rows
of the single-split and doubly-split Magic Squares of order three (see
related Tables). In general, such algebras are embedded maximally into the
algebras of the row above, along the same column, possibly with a commuting
summand; namely, the following maximal embeddings hold (\textit{cfr. e.g.}
\cite{Squaring-Magic}, and Refs. therein):%
\begin{eqnarray}
str_{0} &\supset &aut;  \label{s1} \\
conf &\supset &str_{0}\oplus so(1,1);  \label{s2} \\
qconf &\supset &conf\oplus sl(2,\mathbb{R}),  \label{s3}
\end{eqnarray}%
and these actually hold for any of the (simple and semi-simple) cubic Jordan
algebras introduced above.

Within the physical interpretation of $U$-dualities, the $so(1,1)$ in (\ref%
{s2}) can be regarded as corresponding to the Kaluza-Klein (KK)
compactification radius of the $S^{1}$-reduction from $D=5$ to $D=4$;
alternatively, such an $so(1,1)$ can also be conceived as the Lie algebra
associated to the pseudo-K\"{a}hler connection of the pseudo-special K\"{a}%
hler (and pseudo-Riemannian) symmetric coset $\frac{Conf}{Str_{0}\times
SO(1,1)}$, obtained from\footnote{%
\textquotedblleft $mcs$" denotes the maximal compact subalgebra/subgroup
throughout. Note that $mcs\left( str_{0}\right) =aut$ for $\mathfrak{J}%
=J_{3}^{\mathbb{A}}$, $\mathbb{R}$ and $\mathbb{R}\oplus \Gamma _{m,n}$.} $%
\frac{Str_{0}}{mcs\left( Str_{0}\right) }$ by applying the inverse $R^{\ast
} $-map pertaining to a timelike compactification from $D=5$ Lorentzian
dimensions to $D=4$ spacelike dimensions \cite{R-map,timelike-reduction}. On
the other hand, the $sl(2,\mathbb{R})$ in (\ref{s3}) can be identified as
corresponding to the Ehlers symmetry $sl(2,\mathbb{R})_{Ehlers}$ arising
from the $S^{1}$-reduction from $D=4$ to $D=3$; such an $sl(2,\mathbb{R})$
can also be regarded as the Lie algebra associated to the connection of the
para-quaternionic (and pseudo-Riemannian) symmetric coset $\frac{QConf}{%
Conf\times SL(2,\mathbb{R})}$, obtained from $\frac{Conf}{mcs\left(
Conf\right) }$ by applying the inverse $c^{\ast }$-map pertaining to a
timelike compactification from $D=4$ Lorentzian dimensions to $D=3$
spacelike dimensions\footnote{%
Note that the $\frac{Conf}{mcs\left( Conf\right) }$ is (special) K\"{a}hler
only in certain cases.} \cite{BGM,c-map,timelike-reduction}.

A generalization of the embedding (\ref{s3}) is provided by the so-called
\textit{super-Ehlers embeddings}, recently discussed in \cite{super-Ehlers} :%
\begin{equation}
g_{3}\supset g_{D}\oplus sl(D-2,\mathbb{R})_{Ehlers},  \label{sE}
\end{equation}%
where $g_{3}$ is the $D=3$ $U$-duality Lie algebra (namely, $qconf$ in the
cases treated above), $g_{D}$ is the $U$-duality Lie algebra in $%
3<D\leqslant 11$ dimensions, and $sl(D-2,\mathbb{R})_{Ehlers}$ is the Ehlers
algebra in $D$ dimensions. These embeddings, discussed and generally proven
in \cite{super-Ehlers} (see also \cite{Keurentjes}) have different features,
depending on $D$; they can be maximal or non-maximal, symmetric or
non-symmetric, \textit{etc.}. The (symmetric and maximal) case $D=4$ of (\ref%
{sE}), matching (\ref{s3}), is usually simply named \textit{Ehlers embedding}%
. Morever, the (non-symmetric but maximal) case $D=5$ of (\ref{sE}) pertains
to the \textit{Jordan pairs} introduced above, and it is therefore usually
dubbed \textit{Jordan-Pairs' embedding}.

As the embeddings (\ref{s1})-(\ref{s3}) are obtained by moving along the
columns of the relevant Magic Square (for a fixed row entry), another class
of embeddings can be obtained by moving along the rows of the relevant Magic
Square (for a fixed column entry). In symmetric Magic Squares, as the
non-split ${\mathcal{L}}_{3}(\mathbb{A},\mathbb{B})$ and the double-split ${%
\mathcal{L}}_{3}(\mathbb{A}_{s},\mathbb{B}_{s})$ (respectively given in
Table 4 and 3), these embeddings formally match (\ref{s1})-(\ref{s3}), but
their intepretation corresponds to the restriction from one (division $%
\mathbb{A}$ or split $\mathbb{A}_{s}$) algebra to the next smaller (division
$\mathbb{A}$ or split $\mathbb{A}_{s}$) algebra. For the single-split Magic
Square ${\mathcal{L}}_{3}(\mathbb{A}_{s},\mathbb{B})$ (reported in Table 2),
this class presents different embeddings (at the level of non-compact, real
forms) with respect to the ones given in (\ref{s1})-(\ref{s3}). In general,
it holds that (recall $aut_{\mathbb{B}}=der_{\mathbb{B}}$):%
\begin{equation}
\mathbb{B}\subset \mathbb{C}\Leftrightarrow \left\{
\begin{array}{l}
aut_{\mathbb{B}}\oplus \mathcal{A}_{\mathbb{B}}\subset Aut_{\mathbb{C}}; \\
str_{0}{}_{\mathbb{B}}\oplus \mathcal{A}_{\mathbb{B}}\subset Str_{0\mathbb{C}%
}; \\
conf_{\mathbb{B}}\oplus \mathcal{A}_{\mathbb{B}}\subset Conf_{\mathbb{C}};
\\
qconf_{\mathbb{B}}\oplus \mathcal{A}_{\mathbb{B}}\subset QConf_{\mathbb{C}},%
\end{array}%
\right.   \label{embb}
\end{equation}%
where%
\begin{equation}
\mathcal{A}_{\mathbb{B}}\doteq tri\left( \mathbb{B}\right) \ominus so(%
\mathbb{B}),
\end{equation}%
with $tri\mathbb{\ }$and $so$ respectively denote the triality and
orthogonal (norm-preserving) Lie algebras (see \textit{e.g. }\cite%
{CFMZ1-D=5,MCD}, and Refs. therein). More explicitly :%
\begin{eqnarray}
\mathcal{A}_{\mathbb{B}} &\equiv &\mathcal{A}_{q}\dot{=}tri\left( \mathbb{B}%
\right) \ominus so(\mathbb{B})=\varnothing ,so(3),so(2),\varnothing \text{
for~}q:=\text{dim}_{\mathbb{R}}\mathbb{B}=8,4,2,1~\text{(\textit{i.e.}~for~}%
\mathbb{B}=\mathbb{O},\mathbb{H},\mathbb{C},\mathbb{R}\text{)};  \notag \\
&& \\
\mathcal{A}_{\mathbb{B}_{s}} &\equiv &\widetilde{\mathcal{A}}_{q}\dot{=}%
tri\left( \mathbb{B}_{s}\right) \ominus so(\mathbb{B}_{s})=\varnothing ,sl(2,%
\mathbb{R}),so(1,1)\text{ for~}q:=\text{dim}_{\mathbb{R}}\mathbb{B}%
_{s}=8,4,2~\text{(\textit{i.e.}~for~}\mathbb{B}_{s}=\mathbb{O}_{s},\mathbb{H}%
_{s},\mathbb{C}_{s}\text{)}.  \notag \\
&&
\end{eqnarray}%
In \cite{CFMZ1-D=5}, the appearance of $\mathcal{A}_{q}$ was observed within
the study of the charge orbits of asymptotically flat $0$- (black holes) and
$1$- (black strings) branes in minimal \textquotedblleft magical"
Maxwell-Einstein supergravity theories in $D=5$ space-time dimensions.
Moreover, it is worth noticing that $\mathcal{A}_{q}$ also occurs in the
treatment of supergravity billiards and timelike Kaluza-Klein reductions
(for recent treatment and set of related Refs., see \textit{e.g.} \cite{MCD}%
).

All in all, the right-hand side of (\ref{embb}) expresses the consequences
of the algebraic embedding on its left-hand side at the level of $aut$, $%
str_{0}$, $conf$ and $qconf$ symmetries pertaining to the (division or
split) algebras (as denoted by the subscripts). Clearly, one may consider
non-maximal algebraic embeddings $\mathbb{B}\subset \mathbb{C}$, as well.

Another remarkable class of embedding involves the relation between simple
cubic Jordan algebras $J_{3}^{\mathbb{A}}$ or $J_{3}^{\mathbb{A}_{s}}$ and
some elements of the (bi-parametric) infinite sequence of semi-simple Jordan
algebras $\mathbb{R}\oplus \mathbf{\Gamma }_{m,n}$ introduced above,
exploiting the Jordan-algebraic isomorphisms%
\begin{eqnarray}
J_{2}^{\mathbb{A}} &\cong &\Gamma _{1,q+1}\left( \cong \Gamma
_{q+1,1}\right) ; \\
J_{2}^{\mathbb{A}_{s}} &\cong &\Gamma _{q+2+1,q/2+1},
\end{eqnarray}%
where $q:=$dim$_{R}\mathbb{A}=8,4,2,1$ for $\mathbb{A}=\mathbb{O},\mathbb{H},%
\mathbb{C},\mathbb{R}$, and $q:=$dim$_{R}\mathbb{A}_{s}=8,4,2$ for $\mathbb{A%
}_{s}=\mathbb{O}_{s},\mathbb{H}_{s},\mathbb{C}_{s}$ (see \textit{e.g.} App.
A of \cite{Gun-2} - and Refs. therein - for an introduction to division and
split algebras). Indeed, the following (maximal, rank-preserving)
Jordan-algebraic embeddings hold :%
\begin{eqnarray}
J_{3}^{\mathbb{A}} &\supset &\mathbb{R}\oplus J_{2}^{\mathbb{A}}\cong
\mathbb{R}\oplus \Gamma _{1,q+1}; \\
J_{3}^{\mathbb{A}_{s}} &\supset &\mathbb{R}\oplus J_{2}^{\mathbb{A}%
_{s}}\cong \mathbb{R}\oplus \Gamma _{q/2+1,q/2+1}.
\end{eqnarray}%
Thus, one can consider their consequences at the level of symmetries of
cubic Jordan algebras defined over the corresponding algebras, obtaining :%
\begin{eqnarray}
&&%
\begin{array}{lll}
aut\left( J_{3}^{\mathbb{A}}\right) \supset aut\left( \mathbb{R}\oplus
J_{2}^{\mathbb{A}}\right) \oplus \mathcal{A}_{q}; &  & aut\left( J_{3}^{%
\mathbb{A}_{s}}\right) \supset aut\left( \mathbb{R}\oplus J_{2}^{\mathbb{A}%
_{s}}\right) \oplus \widetilde{\mathcal{A}}_{q}; \\
str_{0}\left( J_{3}^{\mathbb{A}}\right) \supset str_{0}\left( \mathbb{R}%
\oplus J_{2}^{\mathbb{A}}\right) \oplus \mathcal{A}_{q}; &  & str_{0}\left(
J_{3}^{\mathbb{A}_{s}}\right) \supset str_{0}\left( \mathbb{R}\oplus J_{2}^{%
\mathbb{A}_{s}}\right) \oplus \widetilde{\mathcal{A}}_{q}; \\
conf\left( J_{3}^{\mathbb{A}}\right) \supset conf\left( \mathbb{R}\oplus
J_{2}^{\mathbb{A}}\right) \oplus \mathcal{A}_{q}; &  & conf\left( J_{3}^{%
\mathbb{A}_{s}}\right) \supset conf\left( \mathbb{R}\oplus J_{2}^{\mathbb{A}%
_{s}}\right) \oplus \widetilde{\mathcal{A}}_{q}; \\
qconf\left( J_{3}^{\mathbb{A}}\right) \supset qconf\left( \mathbb{R}\oplus
J_{2}^{\mathbb{A}}\right) \oplus \mathcal{A}_{q}; &  & qconf\left( J_{3}^{%
\mathbb{A}_{s}}\right) \supset qconf\left( \mathbb{R}\oplus J_{2}^{\mathbb{A}%
_{s}}\right) \oplus \widetilde{\mathcal{A}}_{q}.%
\end{array}
\notag \\
&&  \label{embs!}
\end{eqnarray}%
It is here worth noticing the maximal nature of the embeddings (\ref{embs!}%
), as well as the presence of the algebras $\mathcal{A}_{q}$ and $\widetilde{%
\mathcal{A}}_{q}$. Within the physical ($U$-duality) interpretation, $%
\mathcal{A}_{q}$ and $\widetilde{\mathcal{A}}_{q}$ are consistent with the
properties of spinors in $q+2$ dimensions, with Lorentzian signature $\left(
1,q+1\right) $ resp. Kleinian signature $\left( q/2+1,q/2+1\right) $;
indeed, the electric-magnetic ($U$-duality) symmetry algebra in $D=6$
(Lorentzian) space-time dimensions is $so\left( 1,q+1\right) \oplus \mathcal{%
A}_{q}$ for $\mathbb{A}$-based theories (which are endowed with minimal,
chiral $\left( 1,0\right) $ supersymmetry) and $so\left( q/2+1,q/2+1\right)
\oplus \widetilde{\mathcal{A}}_{q}$ for $\mathbb{A}_{s}$-based theories
(which are non-supersymmetric for $q=2,4$ - see \cite{Magic-Non-Susy} for a
recent treatment - and endowed with maximal, non-chiral $\left( 2,2\right) $
supersymmetry for $q=8$); \textit{cfr. e.g.} \cite{Kugo-Townsend}, \cite%
{magic-D=6} (and Refs. therein) and \cite{Mkrtchyan-Nersessian} for further
discussion.

In light of the above treatment, the Table of page 10 enjoys a rather
simple Jordan algebraic interpretation, given in the Table at page 18,
which thus characterizes the (maximal) parabolical relation among
non-compact real forms of exceptional Lie algebras in terms of relations
among Lie symmetries of rank-2 (i.e., quadratic) and rank-3 (i.e. cubic)
Jordan algebras.

\begin{center}

%\bigskip
\vbox{\offinterlineskip {\bf Table B} : Jordan algebraic interpretation of maximally parabolically related non-compact real forms of finite-dimensional exceptional Lie
algebras  \medskip \halign{\baselineskip12pt \strut
\vrule#\hskip0.1truecm & #\hfil&
\vrule#\hskip0.1truecm & #\hfil&
\vrule#\hskip0.1truecm & #\hfil&
\vrule#\hskip0.1truecm & #\hfil&
\vrule #\cr
\tablerule
&&& && &&&\cr
&~ $\cg$ &&~  $\cm $ &&~$\cg'$&& ~$\cm'$&\cr
&&&&&&&&\cr
&&&~~$\dim \cn$&&&&&\cr
\tablerule

 &~1~:~$conf\left(M_{2,1}\left( \mathbb{O}\right) \right)$    && ~$su(5,1)$ && ~$qconf\left(
J_{3}^{\mathbb{C}_{s}}\right)$
&&~$conf\left(
J_{3}^{\mathbb{C}_{s}}\right)$&\cr
&&&&& &&&\cr
&&& ~~$21$
&&~$qconf\left(
J_{3}^{\mathbb{C}}\right)$&&~$conf\left(
J_{3}^{\mathbb{C}}\right)$&\cr
\tablerule

&~2~:~$conf\left(M_{2,1}\left( \mathbb{O}\right) \right)$
 && ~$str_{0}\left( \Gamma_{7,1}\right)\oplus u(1)$ && ~$qconf\left(
J_{3}^{\mathbb{C}}\right)$
&&~$str_{0}\left( \Gamma_{5,3}\right)\oplus u(1)$&\cr
&&&&& &&&\cr
&&& ~~$24$
&& && &\cr
\tablerule

&~3~:~$str_{0}\left( J_{3}^{\mathbb{O}}\right)$
 && ~$str_{0}\left( J_{2}^{\mathbb{O}}\right)$ && ~$str_{0}\left( J_{3}^{\mathbb{O}_s}\right)$
&&~$str_{0}\left( J_{2}^{\mathbb{O}_s}\right)$&\cr
&&&&& &&&\cr
&&& ~~$16$
&& && &\cr
\tablerule

&~4~:~$qconf\left(
J_{3}^{\mathbb{C}_s}\right)$
 && ~$str_{0}\left( J_{3}^{\mathbb{C}_s}\right) \oplus sl(2,\bbr)$ && ~$qconf\left(
J_{3}^{\mathbb{C}}\right)$
&&~$str_{0}\left( J_{3}^{\mathbb{C}}\right) \oplus sl(2,\bbr)$&\cr
&&&&& &&&\cr
&&& ~~$29$
&& && &\cr
\tablerule

&~5~:~$conf\left( J_{3}^{\mathbb{O}}\right)$    && ~$str_{0}\left( J_{3}^{\mathbb{O}}\right)$ && ~$conf\left( J_{3}^{\mathbb{O}_s}\right)$
&&~$str_{0}\left( J_{3}^{\mathbb{O}_s}\right)$&\cr
&&&&&&&&\cr &&&~~$27$&&&&&\cr
\tablerule

&~6~:~$conf\left( J_{3}^{\mathbb{O}}\right)$    && ~$str_{0}\left( J_{2}^{\mathbb{O}}\right)\oplus sl(2,\bbr)$  && ~$conf\left( J_{3}^{\mathbb{O}_s}\right)= qconf\left( J_{3}^{\mathbb{H}_s}\right)$
&&~$str_{0}\left( J_{2}^{\mathbb{O}_s}\right) \oplus sl(2,\bbr)$ &\cr
&&&&&&&&\cr
&&&~~$42$&&~$qconf\left( J_{3}^{\mathbb{H}}\right)$&& ~$str_{0}\left( \Gamma_{7,3}\right) \oplus su(2)$&\cr
\tablerule

&~7~:~$conf\left( J_{3}^{\mathbb{O}}\right)$    && ~$conf\left( J_{2}^{\mathbb{O}}\right)$  && ~$conf\left( J_{3}^{\mathbb{O}_s}\right) = qconf\left( J_{3}^{\mathbb{H}_s}\right)$
&&~$conf\left( J_{2}^{\mathbb{O}_s}\right) = conf\left( J_{3}^{\mathbb{H}_s}\right)$ &\cr
&&&&&&&&\cr
&&&~~$33$&&~$qconf\left( J_{3}^{\mathbb{H}}\right)$&& ~$conf\left( J_{3}^{\mathbb{H}}\right)$&\cr
\tablerule

&~8~:~$qconf\left( J_{3}^{\mathbb{H}_s}\right)$
 && ~$sl(4,\mathbb{R})\oplus \widetilde{\mathcal{A}}_{4}\oplus sl(3,\mathbb{R})$ && ~$qconf\left( J_{3}^{\mathbb{H}}\right)$
&&~$su^{\ast }(4)\oplus \mathcal{A}_{4}\oplus sl(3,\mathbb{R})$&\cr
&&&&& &&&\cr
&&& ~~$53$ && && &\cr
\tablerule

&~9~:~$qconf\left( J_{3}^{\mathbb{H}}\right)$
 && ~$str_{0}\left( J_{3}^{\mathbb{H}}\right) \oplus  sl(2,\bbr)$ && ~$qconf\left( J_{3}^{\mathbb{H}_s}\right) = conf\left( J_{3}^{\mathbb{O}_s}\right)$
&&~$str_{0}\left( J_{3}^{\mathbb{H}_{s}}\right) \oplus sl(2,\mathbb{R})
~$&\cr
&&&&& &&&\cr
&&&~~$47$&&~$conf\left( J_{3}^{\mathbb{O}}\right)$&& ~$str_{0}\left( J_{3}^{\mathbb{H}}\right) \oplus su(2) $&\cr
\tablerule

&~10~:~$qconf\left( J_{3}^{\mathbb{O}}\right)$      && ~$conf\left( J_{3}^{\mathbb{O}}\right)$  && ~$qconf\left( J_{3}^{\mathbb{O}_s}\right)$
&& ~$conf\left( J_{3}^{\mathbb{O}_s}\right)$  &\cr
&&&&& &&&\cr
&&& ~~$57$ && && &\cr
\tablerule

&~11~:~$qconf\left( J_{3}^{\mathbb{O}}\right)$      && ~$qconf\left( \mathbb{R}\oplus \Gamma_{8,0}\right)$
&& ~$qconf\left( J_{3}^{\mathbb{O}_s}\right)$      && ~$qconf\left( \mathbb{R}\oplus \Gamma_{4,4}\right)$   &\cr
&&&&& &&&\cr
&&& ~~$78$ && && &\cr
\tablerule

&~12~:~$qconf\left( J_{3}^{\mathbb{O}}\right)$      && ~ $str_{0}\left( J_{3}^{\mathbb{O}}\right) \oplus sl(2,\bbr)$
&& ~$qconf\left( J_{3}^{\mathbb{O}_s}\right)$      && ~ $str_{0}\left( J_{3}^{\mathbb{O}_s}\right) \oplus sl(2,\bbr)$   &\cr
&&&&& &&&\cr
&&& ~~$83$ && && &\cr
\tablerule

&~13~:~$qconf\left( J_{3}^{\mathbb{O}}\right)$      && ~ $str_{0}\left( J_{2}^{\mathbb{O}}\right) \oplus sl(3,\bbr)$
&& ~$qconf\left( J_{3}^{\mathbb{O}_s}\right)$      && ~ $str_{0}\left( J_{2}^{\mathbb{O}_s}\right) \oplus sl(3,\bbr)$   &\cr
&&&&& &&&\cr
&&& ~~$97$ && && &\cr
\tablerule

&~14~:~$der\left( J_{1,2}^{\mathbb{O}}\right)$      && ~ $str_{0}\left( \Gamma_{7,0}\right)$
&& ~$der\left( J_{1,2}^{\mathbb{O}_s}\right)$      && ~ $str_{0}\left( \Gamma_{4,3}\right)$   &\cr
&&&&& &&&\cr
&&& ~~$15$ && && &\cr
\tablerule
}}

\np

\end{center}

Some comments are in order (for further details, see Sec. 5).

\begin{enumerate}
\item The maximal parabolical relation \textbf{1.} hints for a
quasi-conformal interpretation of $E_{6(-14)}$, despite it is characterized
as conformal symmetry algebra of the Hermitian Jordan triple system $M_{2,1}(%
\mathbb{O})$ \cite{Gun-CQC}. In fact, $E_{6(-14)}$ is the U-duality symmetry
of $\mathcal{N}=10$, $D=3$ supergravity (after complete dualization of
1-forms), obtained as the dimensional reduction of $\mathcal{N}=5$, $D=4$
supergravity, which does not admit matter coupling, and whose U-duality
algebra is $su(5,1)$. Usually, U-duality symmetries in $D=4$ and $D=3$ can
be characterized as conformal resp. quasi-conformal symmetries of Jordan
triple systems. Essentially, a quasi-conformal realization of $E_{6(-14)}$
would concern an EFTS of \textit{non-reduced} type (i.e., whose
corresponding FTS is \textit{not} constructed in terms of cubic Jordan
algebras). Work is in progress to investigate this possibility, which would
render the maximal parabolical relation 1. simply the restriction from $qconf
$ to $conf$ symmetries of $M_{2,1}(\mathbb{O})$, $J_{3}^{\mathbb{C}}$ and $%
J_{3}^{\mathbb{C}_{s}}$.

\item In the maximal parabolical relation \textbf{2.}, $str_{0}\left( \Gamma
_{5,3}\right) $ can be enhanced to $str_{0}\left( \Gamma _{6,4}\right)
=qconf\left( \mathbb{R}\oplus J_{2}^{\mathbb{C}}\right) $, thus allowing for
the commuting algebra $u(1)$ to be interpreted as $\mathcal{A}_{2}$; thus,
the parabolic relation between $E_{6(2)}$ and $so(5,3)\oplus u(1)$ can be
interpreted as a consequence of the restriction, at $qconf$ level, of the
simple cubic Jordan algebra $J_{3}^{\mathbb{C}}$ to its maximal semi-simple
Jordan subalgebra $\mathbb{R}\oplus J_{2}^{\mathbb{C}}\simeq \mathbb{R}%
\oplus \Gamma _{3,1}$. Again, this would hint for a $qconf$ interpretation
of $E_{6(-14)}$, despite the fact that $str_{0}\left( \Gamma _{7,1}\right) $
cannot be characterized as $qconf$ symmetry.

\item The maximal parabolical relation \textbf{3.} can be simply explained
as the realization, at $str_{0}$ level, of the Jordan algebraic restriction $%
J_{3}^{\mathbb{O}}\supset \mathbb{R}\oplus J_{2}^{\mathbb{O}}$ and $J_{3}^{%
\mathbb{O}_{s}}\supset \mathbb{R}\oplus J_{2}^{\mathbb{O}_{s}}$, noting that
$str_{0}\left( \mathbb{R}\oplus \Gamma _{m,n}\right) =so(1,1)\oplus
str_{0}\left( \Gamma _{m,n}\right) $.

\item The maximal parabolical relation \textbf{4.} can be simply explained
as the realization, for $J_{3}^{\mathbb{C}_{s}}$ and $J_{3}^{\mathbb{C}}$,
of the embedding $qconf\supset str_{0}\oplus sl(3,\mathbb{R})$ (Jordan
pairs' embedding, or $D=5$ Ehlers embedding), with the $sl(3,\mathbb{R})$
Ehlers symmetry further branched to $gl(2,\mathbb{R})$ in order to give rise
to the grading $so(1,1)$ algebra.

\item The maximal parabolical relation \textbf{5.} can be simply explained
as the realization, for $J_{3}^{\mathbb{O}}$ and $J_{3}^{\mathbb{O}_{s}}$,
of the embedding $conf\supset str_{0}\oplus so(1,1)$.

\item The maximal parabolical relation \textbf{6.} is based on the two-fold
characterization of $E_{7(7)}$ as $conf\left( J_{3}^{\mathbb{O}_{s}}\right) $
(thus allowing a parabolical relation to $conf\left( J_{3}^{\mathbb{O}%
}\right) $), as well as $qconf\left( J_{3}^{\mathbb{H}_{s}}\right) $ (thus
allowing a relation to $qconf\left( J_{3}^{\mathbb{H}}\right) $). Concerning
the $conf$-part of the parabolical relation, we note that $str_{0}\left(
J_{2}^{\mathbb{O}}\right) \oplus sl(2,\mathbb{R})$ can be enhanced to $%
conf\left( J_{2}^{\mathbb{O}}\right) \oplus sl(2,\mathbb{R})\simeq
conf\left( \mathbb{R}\oplus J_{2}^{\mathbb{O}}\right) $, and thus it can be
traced back to the Jordan algebraic restriction $J_{3}^{\mathbb{O}}\supset
\mathbb{R}\oplus J_{2}^{\mathbb{O}}$ at $conf$ level. On the other hand,
concerning the $qconf$-part of the parabolical relation, we note that $%
str_{0}\left( J_{2}^{\mathbb{O}_{s}}\right) \oplus sl(2,\mathbb{R})$ can be
enhanced to $qconf\left( J_{2}^{\mathbb{H}_{s}}\right) \oplus sl(2,\mathbb{R}%
)$, and thus it can be traced back to the non-maximal Jordan algebraic
restriction $J_{3}^{\mathbb{H}_{s}}\supset J_{2}^{\mathbb{H}_{s}}$ at $qconf$
level, with the algebra $sl(2,\mathbb{R})$ interpreted as $\widetilde{%
\mathcal{A}}_{4}$. Analogously, $str_{0}\left( \Gamma _{7,3}\right) \oplus
su(2)$ can be enhanced to $qconf\left( J_{2}^{\mathbb{H}}\right) \oplus su(2)
$, and thus it can be traced back to the non-maximal Jordan algebraic
restriction $J_{3}^{\mathbb{H}}\supset J_{2}^{\mathbb{H}}$ at $qconf$ level,
with the algebra $su(2)$ interpreted as $\mathcal{A}_{4}$.

\item Also the maximal parabolical relation \textbf{7.} is based on the
two-fold characterization of $E_{7(7)}$ as $conf\left( J_{3}^{\mathbb{O}%
_{s}}\right) $ (thus allowing a parabolical relation to $conf\left( J_{3}^{%
\mathbb{O}}\right) $), as well as $qconf\left( J_{3}^{\mathbb{H}_{s}}\right)
$ (thus allowing a relation to $qconf\left( J_{3}^{\mathbb{H}}\right) $).
Concerning the $conf$-part of the parabolical relation, we note that $%
conf\left( J_{2}^{\mathbb{O}}\right) $ can be enhanced to $conf\left( J_{2}^{%
\mathbb{O}}\right) \oplus sl(2,\mathbb{R})\simeq conf\left( \mathbb{R}\oplus
J_{2}^{\mathbb{O}}\right) $, and thus it can be traced back to the Jordan
algebraic restriction $J_{3}^{\mathbb{O}}\supset \mathbb{R}\oplus J_{2}^{%
\mathbb{O}}$ at $conf$ level. Analogously, we note that $conf\left( J_{2}^{%
\mathbb{O}_{s}}\right) $ can be enhanced to $conf\left( J_{2}^{\mathbb{O}%
_{s}}\right) \oplus sl(2,\mathbb{R})\simeq conf\left( \mathbb{R}\oplus
J_{2}^{\mathbb{O}_{s}}\right) $, and thus it can be traced back to the
Jordan algebraic restriction $J_{3}^{\mathbb{O}_{s}}\supset \mathbb{R}\oplus
J_{2}^{\mathbb{O}_{s}}$ at $conf$ level. On the other hand, concerning the $%
qconf$-part of the parabolical relation, it can be interpreted as the
realization, for $J_{3}^{\mathbb{H}_{s}}$ and $J_{3}^{\mathbb{H}}$, of the
embedding restriction $qconf\supset conf\oplus sl(2,\mathbb{R})$ ($D=4$
Ehlers embedding), with the $sl(2,\mathbb{R})$ Ehlers symmetry further
branched to the grading $so(1,1)$ algebra.

\item The maximal parabolical relation \textbf{8.} can be explained as the
realization, for $J_{3}^{\mathbb{H}_{s}}$ and $J_{3}^{\mathbb{H}}$, of the
embedding $qconf\supset str_{0}\oplus sl(3,\mathbb{R})$ (Jordan pairs'
embedding, or $D=5$ Ehlers embedding). Indeed, $sl(4,\mathbb{R})\oplus
\widetilde{\mathcal{A}}_{4}$ and $su^{\ast }(4)\oplus \mathcal{A}_{4}$ can
respectively be enhanced to $str_{0}\left( J_{3}^{\mathbb{H}_{s}}\right) $
and $str_{0}\left( J_{3}^{\mathbb{H}_{s}}\right) $.

\item The maximal parabolical relation \textbf{9.} is based on the two-fold
characterization of $E_{7(7)}$ as $conf\left( J_{3}^{\mathbb{O}_{s}}\right) $
(thus allowing a parabolical relation to $conf\left( J_{3}^{\mathbb{O}%
}\right) $), as well as $qconf\left( J_{3}^{\mathbb{H}_{s}}\right) $ (thus
allowing a relation to $qconf\left( J_{3}^{\mathbb{H}}\right) $). Concerning
the $qconf$-part of the parabolical relation, we note that $str_{0}\left(
J_{3}^{\mathbb{H}}\right) \oplus sl(2,\mathbb{R})$ can be enhanced to $%
str_{0}\left( J_{3}^{\mathbb{H}}\right) \oplus sl(3,\mathbb{R})$ (Jordan
pairs' embedding). Analogously, $str_{0}\left( J_{3}^{\mathbb{H}_{s}}\right)
\oplus sl(2,\mathbb{R})$ can be enhanced to $str_{0}\left( J_{3}^{\mathbb{H}%
_{s}}\right) \oplus sl(3,\mathbb{R})$. On the other hand, concerning the $%
conf$-part of the parabolical relation, we note that the relation between $%
E_{7(7)}$ and $sl(6,\mathbb{R})\oplus sl(2,\mathbb{R})$ can be interpreted
as the realization, at $conf$ level, of the non-maximal algebraic
restriction $\mathbb{O}_{s}\supset \mathbb{C}_{s}$; indeed, $sl(6,\mathbb{R})
$ enjoys a twofold characterization : as $str_{0}\left( J_{3}^{\mathbb{H}%
_{s}}\right) $ and as $conf\left( J_{3}^{\mathbb{C}_{s}}\right) $. Finally,
the relation between $E_{7(-25)}$ and $su^{\ast }(6)\oplus su(2)$ is less
direct, in the sense that it moves before horizontally along the 3rd row of
the single-split (non-symmetric) Magic Square $\mathcal{L}_{3}(\mathbb{A}%
_{s},\mathbb{B})$, from the slot $\mathbb{O}$ to the slot $\mathbb{H}$, and
then it moves vertically along the 3rd column, from the slot $\mathbb{H}_{s}$
to $\mathbb{C}_{s}$ (note that $E_{7(-25)}$ does not have a $qconf$
interpretation).

\item The maximal parabolical relation \textbf{10.} can be simply explained
as the realization, for $J_{3}^{\mathbb{O}}$ and $J_{3}^{\mathbb{O}_{s}}$,
of the embedding $qconf\supset conf\oplus sl(2,\mathbb{R})$ ($D=4$ Ehlers
embedding).

\item The maximal parabolical relation \textbf{11.} can be explained as the
realization, at $qconf$ level, of the maximal Jordan algebraic embeddings $%
J_{3}^{\mathbb{O}}\supset \mathbb{R}\oplus J_{2}^{\mathbb{O}}$ and $J_{3}^{%
\mathbb{O}_{s}}\supset \mathbb{R}\oplus J_{2}^{\mathbb{O}_{s}}$. Indeed, $%
qconf\left( \mathbb{R}\oplus \Gamma _{8,0}\right) $ and $qconf\left( \mathbb{%
R}\oplus \Gamma _{4,4}\right) $ can respectively be enhanced to $qconf\left(
\mathbb{R}\oplus \Gamma _{9,1}\simeq \mathbb{R}\oplus J_{2}^{\mathbb{O}%
}\right) =$ and $qconf\left( \mathbb{R}\oplus \Gamma _{5,5}\simeq \mathbb{R}%
\oplus J_{2}^{\mathbb{O}_{s}}\right) $.

\item The maximal parabolical relation \textbf{12.} can be explained as the
realization, for $J_{3}^{\mathbb{O}}$ and $J_{3}^{\mathbb{O}_{s}}$, of the
embedding $qconf\supset str_{0}\oplus sl(3,\mathbb{R})$ (Jordan pairs'
embedding, or $D=5$ Ehlers embedding). Indeed, $str_{0}\left( J_{3}\right)
\oplus sl(2,\mathbb{R})$ can be enhanced to $str_{0}\left( J_{3}\right)
\oplus gl(2,\mathbb{R})$, which is maximal in\ $str_{0}(J_{3})\oplus sl(3,%
\mathbb{R})$.

\item The maximal parabolical relation \textbf{13.} can also be explained as
the realization, for $J_{3}^{\mathbb{O}}$ and $J_{3}^{\mathbb{O}_{s}}$, of
the embedding $qconf\supset str_{0}\oplus sl(3,\mathbb{R})$ (Jordan pairs'
embedding, or $D=5$ Ehlers embedding). Indeed, $str_{0}(J_{2})\oplus sl(3,%
\mathbb{R})$ can be enhanced to $str_{0}(J_{3})\oplus sl(3,\mathbb{R})$.

\item The maximal parabolical relation \textbf{14.} can be explained as the
realization, for $J_{1,2}^{\mathbb{O}}$ and $J_{1,2}^{\mathbb{O}_{s}}$ (i.e.
for $q=8$), respectively of the maximal (symmetric) embeddings%
\begin{eqnarray}
der\left( J_{1,2}^{\mathbb{A}}\right)  &\supset &str_{0}\left( \Gamma
_{q,1}\right) ; \\
der\left( J_{1,2}^{\mathbb{A}_{s}}\right)  &\supset &str_{0}\left( \Gamma
_{q/2+1,q/2}\right) .
\end{eqnarray}
Indeed, $str_{0}\left( \Gamma _{q-1,0}\right) $ and $str_{0}\left( \Gamma
_{q/2,q/2-1}\right) $ can trivially be enhanced to $str_{0}\left( \Gamma
_{q,1}\right) $ resp. $str_{0}\left( \Gamma _{q/2+1,q/2}\right) $, which in
turn are maximal in $str_{0}\left( \Gamma _{q+1,1}\simeq J_{2}^{\mathbb{A}%
}\right) $ and in $str_{0}\left( \Gamma _{q/2+1,q/2+1}\simeq J_{2}^{\mathbb{A%
}_{s}}\right) $, respectively.
\end{enumerate}

\setcounter{equation}{0}

\section{Jordan Structures in Maximal Parabolics of Exceptional Lie Algebras
: Analysis}

In this Section, we will analyze all maximal parabolic subalgebras
(shortened as \textit{parabolics}) of all non-compact, real forms of
finite-dimensional exceptional Lie algebras, determining them by means of
(chains of) maximal embeddings, and providing Jordan algebraic
interpretations for them, in light of the treatment given in previous \
Section. We will also provide an $\mathcal{M}^{\max }$-covariant
decomposition of the (generally reducible) vector spaces $\mathcal{N}^{\max
} $'s, occurring in the Bruhat branchings giving rise to the maximal
parabolics. The numbering of maximal parabolics refers to the listing of
Sec. \ref{sec:MP}.

\md

\subsection{$E_{6(6)}$}

This is the maximally non-compact (\textit{split}) real form of $E_{6}$. Its
Jordan interpretation is essentially twofold (due to the symmetry of the
double-split Magic Square ${\mathcal{L}}_{3}(\mathbb{A}_{s},\mathbb{B}_{s})$
\cite{BS}, reported in Table 3) :%
\begin{eqnarray}
E_{6(6)} &\cong &qconf\left( J_{3}^{\mathbb{C}_{s}}\right)  \label{uno-3} \\
&\cong &str_{0}\left( J_{3}^{\mathbb{O}_{s}}\right) \cong der\left( J_{3}^{%
\mathbb{O}_{s}},J_{3}^{\mathbb{O}_{s}\prime }\right) \ominus so(1,1).
\label{uno-4}
\end{eqnarray}

\md

\subsubsection{$\mathcal{P}_{1}^{6(6)}\cong \mathcal{P}_{5}^{6(6)}$\label{611}}

The maximal parabolics ~ $\mathcal{P}_{1}^{6(6)}\cong \mathcal{P}_{5}^{6(6)}$
~ from \eqref{maxsixa} corresponds to the Bruhat decomposition %
\eqref{bruhatg}:
\begin{eqnarray}
E_{6(6)} &=&\mathcal{N}_{1}^{6(6)-}\oplus so(5,5)\oplus so(1,1)\oplus
\mathcal{N}_{1}^{6(6)+}\ ,  \label{3} \\
\mathbf{78} &=&\mathbf{16}_{-3}\oplus \mathbf{45}_{0}\oplus \mathbf{1}%
_{0}\oplus \mathbf{16}_{3}^{\prime }\ ,  \label{4}
\end{eqnarray}%
yielding a $3$-grading\footnote{%
The subscripts denote $so(1,1)$-weights throughout. Unless otherwise
indicated, all embeddings are symmetric; non-symmetric embeddings will be
denoted by a \textquotedblleft $ns$" upperscript. Only (chains of) maximal
embeddings are considered throughout.}. Since $so(5,5)\cong str_{0}\left(
J_{2}^{\mathbb{O}_{s}}\right) $, \textit{at least} two Jordan algebraic
interpretations (denoted by $I$ and $II$) of (\ref{3})-(\ref{4}) can be
given:

\begin{enumerate}
\item the first one is%
\begin{equation}
str_{0}\left( J_{3}^{\mathbb{O}_{s}}\right) \supset str_{0}\left( J_{2}^{%
\mathbb{O}_{s}}\right) \oplus so(1,1),  \label{611I}
\end{equation}%
where the $so(1,1)$ generating the $3$-grading is the Kaluza-Klein (KK) $%
so(1,1)$ of the $S^{1}$-reduction $D=6\rightarrow 5$.

\item the second one stems from the Jordan algebraic embedding $J_{3}^{%
\mathbb{O}_{s}}\supset \mathbb{R}\oplus J_{2}^{\mathbb{O}_{s}}$, at the
level of $str_{0}$ :%
\begin{equation}
str_{0}\left( J_{3}^{\mathbb{O}_{s}}\right) \supset str_{0}\left( \mathbb{R}%
\oplus J_{2}^{\mathbb{O}_{s}}\right) \times \widetilde{\mathcal{A}}_{8},
\label{611II}
\end{equation}%
where the $so(1,1)$ generating the $3$-grading is the dilatonic $so(1,1)$ in
$D=5$.
\end{enumerate}

\md

\subsubsection{$\cp_{2}^{6(6)}$~\label{612}}

The maximal parabolics ~$\cp_{2}^{6(6)}$~ from \eqref{maxsixa} corresponds
to the Bruhat decomposition \eqref{bruhatg}:
\begin{equation}
E_{6(6)}=\cn_{2}^{6(6)-}\oplus sl(6,\mathbb{R})\oplus so(1,1)\oplus \cn%
_{2}^{6(6)+}\ ,  \label{11}
\end{equation}%
which can be obtained by the following chain of embeddings :
\begin{eqnarray}
E_{6(6)} &\supset &sl(6,\bbr)\oplus sl(2,\mathbb{R})\supset sl(6,\bbr)\oplus
so(1,1);  \label{pre-12} \\
\mathbf{78} &=&\left( \mathbf{35},\mathbf{1}\right) \oplus \left( \mathbf{1},%
\mathbf{3}\right) \oplus \left( \mathbf{20},\mathbf{2}\right) =\mathbf{1}%
_{-2}\oplus \mathbf{20}_{-1}\oplus \mathbf{35}_{0}\oplus \mathbf{1}%
_{0}\oplus \mathbf{20}_{1}\ \oplus \mathbf{1}_{2},  \label{12}
\end{eqnarray}%
exhibiting a $5$-grading of contact, with ~$\mathcal{N}_{2}^{6(6)\pm }=%
\mathbf{1}_{\pm 2}+\mathbf{20}_{\pm 1}\,$.

Since $sl(6,\mathbb{R})$ has a two-fold Jordan algebraic interpretation :
\begin{eqnarray}
sl(6,\mathbb{R}) &\cong &str_{0}\left( J_{3}^{\mathbb{H}_{s}}\right) \cong
der\left( J_{3}^{\mathbb{H}_{s}},J_{3}^{\mathbb{H}_{s}\prime }\right)
\ominus so(1,1)  \label{SL6-1} \\
&\cong &conf\left( J_{3}^{\mathbb{C}_{s}}\right) ,  \label{SL6-3}
\end{eqnarray}%
the first step of the above chain can be interpreted in two ways :

\begin{enumerate}
\item
\begin{equation*}
qconf\left( J_{3}^{\mathbb{C}_{s}}\right) \supset conf\left( J_{3}^{\mathbb{C%
}_{s}}\right) \oplus sl(2,\mathbb{R}),
\end{equation*}%
namely the $D=4$ Ehlers embedding for $J_{3}^{\mathbb{C}_{s}}$.

\item as a consequence of the split algebraic embedding $\mathbb{O}%
_{s}\supset \mathbb{H}_{s}\Rightarrow \mathfrak{J}_{3}^{\mathbb{O}%
_{s}}\supset \mathfrak{J}_{3}^{\mathbb{H}_{s}}$, at the level of $str_{0}$
it holds that :%
\begin{equation}
str_{0}\left( J_{3}^{\mathbb{O}_{s}}\right) \supset Str_{0}\left( J_{3}^{%
\mathbb{H}_{s}}\right) \times \widetilde{\mathcal{A}}_{4}.  \label{j612}
\end{equation}
\end{enumerate}

\md

\subsubsection{$\cp_{3}^{6(6)}\cong \cp_{6}^{6(6)}$\label{613}}

The maximal parabolics ~ $\cp_{3}^{6(6)}\cong \cp_{6}^{6(6)}$ ~ from %
\eqref{maxsixa} corresponds to the Bruhat decomposition \eqref{bruhatg}:
\begin{eqnarray}
E_{6(6)} &=&(\cn^{-})_{3}^{6(6)}\oplus (sl(5,\mathbb{R})\oplus sl(2,\mathbb{R%
})\oplus so(1,1)\oplus (\cn^{+})_{3}^{6(6)}\ ,  \label{5} \\
\mathbf{78} &=&(\mathbf{5}^{\prime },\mathbf{1})_{-6}\oplus (\mathbf{10},%
\mathbf{2})_{-3}\oplus (\mathbf{24},\mathbf{1})_{0}\oplus (\mathbf{1},%
\mathbf{3})_{0}\oplus \left( \mathbf{1},\mathbf{1}\right) _{0}\ \oplus (%
\mathbf{10}^{\prime },\mathbf{2})_{3}\oplus (\mathbf{5},\mathbf{1})_{6},
\end{eqnarray}%
thus yielding a $5$-grading (recall that $\dim \mathcal{N}_{3}^{6(6)\pm }=25$%
), which can be obtained by the following chain of embeddings:%
\begin{equation}
E_{6(6)}\supset sl(6,\mathbb{R})\oplus sl(2,\mathbb{R})\supset sl(5,\mathbb{R%
})\oplus sl(2,\mathbb{R})\oplus so(1,1).
\end{equation}%
The first embedding is the very same of the chain (\ref{pre-12}); here the $%
sl(6,\mathbb{R})$ is branched to generate the $so(1,1)$ producing the
parabolic $5$-grading, whereas in (\ref{pre-12}) it is the branching of the $%
sl(2,\mathbb{R})$ algebra to produce the relevant $so(1,1)$.

\md

\subsubsection{$\cp_{4}^{6(6)}$\label{614}}

The maximal parabolics ~$\cp_{4}^{6(6)}$ ~ from \eqref{maxsixa} corresponds
to the Bruhat decomposition \eqref{bruhatg}:
\begin{equation}
E_{6(6)}=\cn_{4}^{6(6)-}\oplus sl(3,\bbr)\oplus sl(3,\bbr)\oplus sl(2,\bbr%
)\oplus so(1,1)\oplus \cn_{4}^{6(6)+}\ ,
\end{equation}%
which can be obtained at least in two different ways, associated to two
embedding chains\footnote{%
For a recent quantum informational interpretation of the first step of the
chain, see \textit{e.g.} \cite{DF-E6}.}, respectively denoted by $1$ and $2$
:%
\begin{eqnarray}
1 &:&E_{6(6)}\supset ^{ns}sl(3,\mathbb{R})_{\mathbf{1}}\oplus sl(3,\mathbb{R}%
)_{\mathbf{2}}\oplus sl(3,\mathbb{R})_{\mathbf{3}}\supset sl(3,\mathbb{R}%
)\oplus sl(3,\mathbb{R})\oplus sl(2,\mathbb{R})\oplus so(1,1);  \notag \\
&&  \label{7} \\
\mathbf{78} &=&\left( \mathbf{8,1,1}\right) \oplus \left( \mathbf{1,8,1}%
\right) \oplus \left( \mathbf{1,1,8}\right) \oplus \left( \mathbf{3,3,3}%
^{\prime }\right) \oplus \left( \mathbf{3}^{\prime }\mathbf{,3}^{\prime }%
\mathbf{,3}\right)  \notag \\
&=&\left( \mathbf{8,1,1}\right) _{0}\oplus \left( \mathbf{1,8,1}\right)
_{0}\oplus \left( \mathbf{1,1,3}\right) _{0}\oplus \left( \mathbf{1,1,1}%
\right) _{0}\oplus \left( \mathbf{1,1,2}\right) _{3}\oplus \left( \mathbf{%
1,1,2}\right) _{-3}  \notag \\
&&\oplus \left( \mathbf{3,3,1}\right) _{2}\oplus \left( \mathbf{3}^{\prime }%
\mathbf{,3}^{\prime }\mathbf{,1}\right) _{-2}\oplus \left( \mathbf{3,3,2}%
\right) _{-1}\oplus \left( \mathbf{3}^{\prime }\mathbf{,3}^{\prime }\mathbf{%
,2}\right) _{1};  \label{8}
\end{eqnarray}%
\begin{eqnarray}
2 &:&E_{6(6)}\supset sl(6,\mathbb{R})\oplus sl(2,\mathbb{R})\supset sl(3,%
\mathbb{R})\oplus sl(3,\mathbb{R})\oplus sl(2,\mathbb{R})\oplus so(1,1);
\label{9} \\
\mathbf{78} &=&\left( \mathbf{35,1}\right) \oplus \left( \mathbf{1,3}\right)
\oplus \left( \mathbf{20,2}\right)  \notag \\
&=&\left( \mathbf{8,1,1}\right) _{0}\oplus \left( \mathbf{1,8,1}\right)
_{0}\oplus \left( \mathbf{1,1,3}\right) _{0}\oplus \left( \mathbf{1,1,1}%
\right) _{0}\oplus \left( \mathbf{1,1,2}\right) _{3}\oplus \left( \mathbf{%
1,1,2}\right) _{-3}  \notag \\
&&\oplus \left( \mathbf{3,3}^{\prime }\mathbf{,1}\right) _{2}\oplus \left(
\mathbf{3}^{\prime }\mathbf{,3,1}\right) _{-2}\oplus \left( \mathbf{3,3}%
^{\prime }\mathbf{,2}\right) _{-1}\oplus \left( \mathbf{3}^{\prime }\mathbf{%
,3,2}\right) _{1}.  \label{10}
\end{eqnarray}%
(\ref{7})-(\ref{8}) and (\ref{9})-(\ref{10}) both yield a $7$-grading with
the same $\mathcal{M}_{4}^{6(6)}$ but with different $\mathcal{N}$'s :

\begin{enumerate}
\item (\ref{7})-(\ref{8}) implies $\mathcal{N}_{4}^{6(6)+}=\left( \mathbf{%
1,1,2}\right) _{3}+\left( \mathbf{3,3,1}\right) _{2}+\left( \mathbf{3}%
^{\prime }\mathbf{,3}^{\prime }\mathbf{,2}\right) _{1}$;

\item (\ref{9})-(\ref{10}) yields $\mathcal{N}_{4}^{6(6)+}=\left( \mathbf{%
1,1,2}\right) _{3}+\left( \mathbf{3,3}^{\prime }\mathbf{,1}\right)
_{2}+\left( \mathbf{3}^{\prime }\mathbf{,3,2}\right) _{1}$.
\end{enumerate}

Namely, $1\leftrightarrow 2$ \textit{iff} $\mathbf{3}^{\prime
}\leftrightarrow \mathbf{3}$ in the second $sl(3,\mathbb{R})$ algebra of $%
\mathcal{M}_{4}^{6(6)}$.

The Jordan-algebraic interpretation of the first step of chains (\ref{7})-(%
\ref{8}) and (\ref{9})-(\ref{10}) is \textit{at least} twofold. Indeed, it
holds that :%
\begin{equation}
sl(3,\mathbb{R})\oplus sl(3,\mathbb{R})\cong str_{0}\left( J_{3}^{\mathbb{C}%
_{s}}\right) \cong der\left( J_{3}^{\mathbb{C}_{s}},J_{3}^{\mathbb{C}%
_{s}\prime }\right) \ominus so(1,1).
\end{equation}%
A first interpretation (denoted by $I$) is provided by the Jordan pairs'
(JP) embedding for the $\mathbb{C}_{s}$-based gravity theories \cite%
{Magic-Non-Susy} :%
\begin{equation}
I:qconf\left( J_{3}^{\mathbb{O}_{s}}\right) \supset str_{0}\left( J_{3}^{%
\mathbb{C}_{s}}\right) \oplus sl(3,\mathbb{R}),
\end{equation}%
where the $sl(3,\mathbb{R})$ commuting factor is the \textit{Ehlers group}
in $D=5$ (Lorentzian-signed) space-time dimensions. A\ second interpretation
(denoted by $II$) is based on the non-maximal split algebraic embedding $%
\mathbb{O}_{s}\supset \mathbb{C}_{s}\Rightarrow J_{3}^{\mathbb{O}%
_{s}}\supset J_{3}^{\mathbb{C}_{s}}$, at the level of $str_{0}$ :
\begin{equation}
II:str_{0}\left( J_{3}^{\mathbb{O}_{s}}\right) \supset str_{0}\left( J_{3}^{%
\mathbb{C}_{s}}\right) \oplus sl(3,\mathbb{R}).
\end{equation}

\md

\subsection{$E_{6(2)}$}

The Jordan algebraic interpretation of\ $E_{6(2)}$ is%
\begin{equation}
E_{6(2)}\cong qconf\left( J_{3}^{\mathbb{C}}\right) .
\end{equation}

\md

\subsubsection{$\cp_{1}^{6(2)}$\label{621}}

The maximal parabolics ~$\cp_{1}^{6(2)}$ ~ from \eqref{maxsixb} corresponds
to the Bruhat decomposition:
\begin{equation}
E_{6(2)}~=\cn_{1}^{6(2)-}\oplus ~so(5,3)\oplus so(2)\oplus so(1,1)\oplus \cn%
_{1}^{6(2)+}\ ,
\end{equation}%
which can be obtained through the following embedding chain:
\begin{eqnarray}
E_{6(2)} &\supset &so(6,4)\oplus so(2)\supset so(5,3)\oplus so(2)\oplus
so(1,1);  \label{13} \\
\mathbf{78} &=&\mathbf{28}_{0,0}\oplus \mathbf{1}_{0,0}\oplus \mathbf{1}%
_{0,0}\oplus \mathbf{8}_{v,0,2}\oplus \mathbf{8}_{v,0,-2}\oplus \mathbf{8}%
_{c,-3,1}\oplus \mathbf{8}_{s,-3,-1}\oplus \mathbf{8}_{s,3,1}\oplus \mathbf{8%
}_{c,3,-1},  \label{14}
\end{eqnarray}%
and accounting the action of ~$\cm_{1}^{6(2)}$ ~ we obtain a $5$-grading%
\footnote{%
The first subscripts denote $so(2)$ charges. Note that we use Slansky's
conventions on the triality of $so(8)$ \cite{Slansky}.}, with $\mathcal{N}%
_{1}^{6(2)\pm }=\mathbf{8}_{v,0,\pm 2}+\mathbf{8}_{c,-3,\pm 1}+\mathbf{8}%
_{s,3,\pm 1}$.\newline
The Jordan algebraic interpretation of the first embedding of (\ref{13})-(%
\ref{14}) stems from the Jordan algebraic embedding $\mathfrak{J}_{3}^{%
\mathbb{C}}\supset \mathbb{R}\oplus \mathfrak{J}_{2}^{\mathbb{C}}$, at the $%
qconf$ level :%
\begin{equation}
qconf\left( \mathfrak{J}_{3}^{\mathbb{C}}\right) \supset qconf\left( \mathbb{%
R}\oplus \mathfrak{J}_{2}^{\mathbb{C}}\right) \oplus \mathcal{A}_{2}.
\label{uno}
\end{equation}

\md

\subsubsection{$\cp_{2}^{6(2)}$\label{622}}

The maximal parabolics ~$\cp_{2}^{6(2)}$ ~ from \eqref{maxsixb} corresponds
to the Bruhat decomposition:
\begin{equation}
E_{6(2)}=\cn_{2}^{6(2)-}\oplus sl(3,\bbr)\oplus sl(2,\bbc)_{\bbr}\oplus
so(2)\oplus so(1,1)\oplus \cn_{2}^{6(2)}\ ,
\end{equation}%
which can be obtained by the following embedding chain:
\begin{eqnarray}
E_{6(2)} &\supset &^{ns}sl(3,\mathbb{R})\oplus sl(3,\mathbb{C})_{\mathbb{R}%
}\supset sl(3,\mathbb{R})\oplus sl(2,\mathbb{C})_{\mathbb{R}}\oplus
so(2)\oplus so(1,1);  \label{15} \\
\mathbf{78} &=&\left( \mathbf{8,1,1}\right) \oplus \left( \mathbf{1,8,1}%
\right) \oplus \left( \mathbf{1,1,8}\right) \oplus \left( \mathbf{3,3,}%
\overline{\mathbf{3}}\right) \oplus \left( \mathbf{3}^{\prime }\mathbf{,}%
\overline{\mathbf{3}}\mathbf{,3}\right)  \notag \\
&=&\left( \mathbf{8,1,1}\right) _{0,0}\oplus \left( \mathbf{1,3,1}\right)
_{0,0}\oplus \left( \mathbf{1,1,1}\right) _{0,0}\oplus \left( \mathbf{1,1,3}%
\right) _{0,0}\oplus \left( \mathbf{1,1,1}\right) _{0,0}  \notag \\
&&\oplus \left( \mathbf{1,2,1}\right) _{3,-3}\oplus \left( \mathbf{1,2,1}%
\right) _{-3,3}\oplus \left( \mathbf{1,1,2}\right) _{3,-3}\oplus \left(
\mathbf{1,1,2}\right) _{-3,3}  \notag \\
&&\oplus \left( \mathbf{3,2,2}\right) _{0,-2}\oplus \left( \mathbf{3,2,1}%
\right) _{3,1}\oplus \left( \mathbf{3,1,2}\right) _{-3,1}\oplus \left(
\mathbf{3,1,1}\right) _{0,4}  \notag \\
&&\oplus \left( \mathbf{3}^{\prime }\mathbf{,2,2}\right) _{0,2}\oplus \left(
\mathbf{3}^{\prime }\mathbf{,2,1}\right) _{-3,-1}\oplus \left( \mathbf{3}%
^{\prime }\mathbf{,1,2}\right) _{3,-1}\oplus \left( \mathbf{3}^{\prime }%
\mathbf{,1,1}\right) _{0,-4},  \label{16}
\end{eqnarray}%
thus yielding $9$-grading, noting that $\mathcal{N}_{2}^{6(2)\oplus }=\left(
\mathbf{1,2,1}\right) _{-3,3}\oplus \left( \mathbf{1,1,2}\right)
_{-3,3}\oplus \left( \mathbf{3,2,1}\right) _{3,1}\oplus \left( \mathbf{3,1,2}%
\right) _{-3,1}\oplus \left( \mathbf{3,1,1}\right) _{0,4}\oplus \left(
\mathbf{3}^{\prime }\mathbf{,2,2}\right) _{0,2}\,$. It holds that :%
\begin{equation}
sl(3,\mathbb{C})_{\mathbb{R}}\cong Str_{0}\left( \mathfrak{J}_{3}^{\mathbb{C}%
}\right) \cong Aut\left( \mathfrak{J}_{3}^{\mathbb{C}},\overline{\mathfrak{J}%
_{3}^{\mathbb{C}}}\right) /so(1,1),
\end{equation}%
and a possible interpretation of the first step of the chain (\ref{15})-(\ref%
{16}) is provided by the JP embedding for the $\mathbb{C}$-based gravity
theories :%
\begin{equation}
QConf\left( \mathfrak{J}_{3}^{\mathbb{C}}\right) \supset Str_{0}\left(
\mathfrak{J}_{3}^{\mathbb{C}}\right) \oplus sl(3,\mathbb{R}),
\end{equation}%
where, once again, the $sl(3,\mathbb{R})$ commuting factor is the \textit{%
Ehlers group} in $D=5$ (Lorentzian-signed) space-time dimensions. Also the
second step of the chain (\ref{15})-(\ref{16}) can be given a
Jordan-algebraic interpretation; this latter stems from the same embedding
of (\ref{uno}), but here considered at the level of the reduced structure
symmetries :%
\begin{equation}
Str_{0}\left( \mathfrak{J}_{3}^{\mathbb{C}}\right) \supset Str_{0}\left(
\mathbb{R}\oplus \mathfrak{J}_{2}^{\mathbb{C}}\right) \oplus \mathcal{A}_{2}.
\label{due}
\end{equation}

\md

\subsubsection{$\cp_{3}^{6(2)}$\label{623}}

The maximal parabolics ~$\cp_{3}^{6(2)}$ ~ from \eqref{maxsixb} corresponds
to the Bruhat decomposition:
\begin{equation}
E_{6(2)}=\cn_{3}^{6(2)-}\oplus sl(3,\bbc)_{\bbr}\oplus sl(2,\bbr)\oplus
so(1,1)\oplus \cn_{3}^{6(2)+}\ ,
\end{equation}%
which can be realized by \textit{at least} two chains of embeddings,
respectively denoted by $1$ and $2$ :%
\begin{eqnarray}
1 &:&E_{6(2)}\supset su(3,3)\oplus sl(2,\mathbb{R})\supset sl(3,\mathbb{C})_{%
\mathbb{R}}\oplus sl(2,\mathbb{R})\oplus so(1,1);  \label{17} \\
\mathbf{78} &=&\left( \mathbf{35},\mathbf{1}\right) \oplus \left( \mathbf{1},%
\mathbf{3}\right) \oplus \left( \mathbf{20},\mathbf{2}\right) =\left\{
\begin{array}{l}
\left( \mathbf{8},\mathbf{1,1}\right) _{0}\oplus \left( \mathbf{1},\mathbf{%
8,1}\right) _{0}\oplus \left( \mathbf{1},\mathbf{1,3}\right) _{0}\oplus
\left( \mathbf{1},\mathbf{1,1}\right) _{0} \\
\oplus \left( \mathbf{3},\mathbf{3}^{\prime }\mathbf{,1}\right) _{2}\oplus
\left( \mathbf{3}^{\prime },\mathbf{3,1}\right) _{-2}\oplus \left( \mathbf{1}%
,\mathbf{1,2}\right) _{3} \\
\oplus \left( \mathbf{1},\mathbf{1,2}\right) _{-3}\oplus \left( \mathbf{3},%
\mathbf{3}^{\prime }\mathbf{,2}\right) _{-1}\oplus \left( \mathbf{3}^{\prime
},\mathbf{3,2}\right) _{1};%
\end{array}%
\right.  \label{18}
\end{eqnarray}%
\begin{eqnarray}
2 &:&E_{6(2)}\supset ^{ns}sl(3,\mathbb{C})_{\mathbb{R}}\oplus sl(3,\mathbb{R}%
)\supset sl(3,\mathbb{C})_{\mathbb{R}}\oplus sl(2,\mathbb{R})\oplus so(1,1);
\label{19} \\
\mathbf{78} &=&\left( \mathbf{8,1,1}\right) \oplus \left( \mathbf{1,8,1}%
\right) \oplus \left( \mathbf{1,1,8}\right) \oplus \left( \mathbf{3,3}%
^{\prime },\mathbf{3}\right) \oplus \left( \mathbf{3}^{\prime }\mathbf{,3,3}%
^{\prime }\right)  \notag \\
&=&\left\{
\begin{array}{l}
\left( \mathbf{8},\mathbf{1,1}\right) _{0}\oplus \left( \mathbf{1},\mathbf{%
8,1}\right) _{0}\oplus \left( \mathbf{1},\mathbf{1,3}\right) _{0}\oplus
\left( \mathbf{1},\mathbf{1,1}\right) _{0} \\
\oplus \left( \mathbf{3},\mathbf{3}^{\prime }\mathbf{,1}\right) _{-2}\oplus
\left( \mathbf{3}^{\prime },\mathbf{3,1}\right) _{2}\oplus \left( \mathbf{1},%
\mathbf{1,2}\right) _{3} \\
\oplus \left( \mathbf{1},\mathbf{1,2}\right) _{-3}\oplus \left( \mathbf{3},%
\mathbf{3}^{\prime }\mathbf{,2}\right) _{1}\oplus \left( \mathbf{3}^{\prime
},\mathbf{3,2}\right) _{-1}.%
\end{array}%
\right.  \label{20}
\end{eqnarray}%
(\ref{17})-(\ref{18}) and (\ref{19})-(\ref{20}) both yield a $7$-grading,
but with different $\mathcal{N}$'s (both of real dimension $29$) :

\begin{enumerate}
\item (\ref{17})-(\ref{18}) implies $\mathcal{N}_{3}^{6(2)+}=\left( \mathbf{1%
},\mathbf{1,2}\right) _{3}+\left( \mathbf{3},\mathbf{3}^{\prime }\mathbf{,1}%
\right) _{2}+\left( \mathbf{3}^{\prime },\mathbf{3,2}\right) _{1}$;

\item (\ref{19})-(\ref{20}) yields $\mathcal{N}_{3}^{6(2)+}=\left( \mathbf{%
1,1,2}\right) _{3}+\left( \mathbf{3}^{\prime },\mathbf{3,1}\right)
_{2}+\left( \mathbf{3},\mathbf{3}^{\prime }\mathbf{,2}\right) _{1}$.
\end{enumerate}

Namely, $1\leftrightarrow 2$ \textit{iff} $\mathbf{3}^{\prime
}\leftrightarrow \mathbf{3}$ in the $sl(3,\mathbb{C})_{\mathbb{R}}$ algebra
of $\mathcal{M}_{3}^{6(2)}$.

Both steps of each of the chains $1$ and $2$ admits \textit{at least} one
Jordan algebraic interpretation, as follows : the chain $1$ (\ref{17})-(\ref%
{18}) starts with the so-called \textit{Ehlers embedding} for $\mathbb{C}$%
-based magic supergravity theories, and then proceeds with an inverse $%
R^{\ast }$-map (generating an $so(1,1)_{KK}$, determining the $7$-grading):
\begin{equation}
1:qconf\left( J_{3}^{\mathbb{C}}\right) \supset conf\left( J_{3}^{\mathbb{C}%
}\right) \oplus sl(2,\mathbb{R})_{Ehlers}\supset str_{0}\left( J_{3}^{%
\mathbb{C}}\right) \oplus sl(2,\mathbb{R})_{Ehlers}\oplus so(1,1)_{KK}.
\label{tre}
\end{equation}%
On the other hand, the chain $2$ (\ref{19})-(\ref{20}) starts with the JP
embedding for the $\mathbb{C}$-based magic supergravity theories (generating
the\footnote{%
The $D=5$ \textit{Ehlers }$sl(3,\mathbb{R})_{Ehlers}$ can also be regarded
as the \textit{enhancement} of $sl(2,\mathbb{R})_{Ehlers}\oplus so(1,1)_{KK}$%
, obtained by a composition of an inverse $c^{\ast }$-map and of an inverse $%
R^{\ast }$-map.} $D=5$ \textit{Ehlers }$sl(3,\mathbb{R})_{Ehlers}$), and
then proceeds with a further branching of this latter symmetry into $sl(2,%
\mathbb{R})_{Ehlers}$ ($D=4$ Ehlers) $\oplus \ so(1,1)_{KK}$ :%
\begin{equation}
2:qconf\left( J_{3}^{\mathbb{C}}\right) \supset str_{0}\left( J_{3}^{\mathbb{%
C}}\right) \oplus sl(3,\mathbb{R})_{Ehlers}\supset str_{0}\left( J_{3}^{%
\mathbb{C}}\right) \oplus sl(2,\mathbb{R})_{Ehlers}\oplus so(1,1)_{KK}.
\label{quattro}
\end{equation}

\md

\subsubsection{$\cp_{4}^{6(2)}$\label{624}}

The maximal parabolics ~$\cp_{4}^{6(2)}$ ~ from \eqref{maxsixb} corresponds
to the Bruhat decomposition:
\begin{equation}
E_{6(2)}=\cn_{4}^{6(2)-}\oplus su(3,3)\oplus so(1,1)\oplus \cn{6(2)+}_{4}\ ,
\end{equation}%
which can be obtained by the following embedding chain:
\begin{eqnarray}
E_{6(2)} &\supset &su(3,3)\oplus sl(2,\mathbb{R})\supset su(3,3)\oplus
so(1,1);  \label{21} \\
\mathbf{78} &=&\left( \mathbf{35},\mathbf{1}\right) \oplus \left( \mathbf{1},%
\mathbf{3}\right) \oplus \left( \mathbf{20},\mathbf{2}\right) =\mathbf{35}%
_{0}\oplus \mathbf{1}_{0}\oplus \mathbf{1}_{2}\oplus \mathbf{1}_{-2}\oplus
\mathbf{20}_{1}\oplus \mathbf{20}_{-1},  \label{22}
\end{eqnarray}%
exhibiting a $5$-grading with ~$\mathcal{N}_{4}^{6(2)\pm }=\mathbf{1}_{\pm
2}+\mathbf{20}_{\pm 1}\,$.

The Jordan algebraic interpretation starts with the Ehlers embedding for $%
\mathbb{C}$-based theories, and then proceeds by further branching this
latter symmetry, in order to generate the $so(1,1)$ responsible for the $5$%
-grading :%
\begin{equation}
qconf\left( J_{3}^{\mathbb{C}}\right) \supset conf\left( J_{3}^{\mathbb{C}%
}\right) \oplus sl(2,\mathbb{R})_{Ehlers}\supset conf\left( J_{3}^{\mathbb{C}%
}\right) \oplus so(1,1).  \label{cinque}
\end{equation}%
Notice that in this case the $so(1,1)$ has not a KK interpretation, but it
is rather the non-compact Cartan of $sl(2,\mathbb{R})_{Ehlers}$.

\md

\subsection{$E_{6(-14)}$}

The Jordan interpretation of\ $E_{6(-14)}$ is twofold :%
\begin{eqnarray}
E_{6(-14)} &\cong &conf\left( M_{1,2}\left( \mathbb{O}\right) \right) \\
&\cong &\mathcal{K}\left( J_{3}^{\mathbb{O}}\right) ,
\end{eqnarray}%
where, as mentioned in the previous Section, $M_{1,2}\left( \mathbb{O}%
\right) $ denotes an \textit{Hermitian} \textit{Jordan triple system} formed
by octonionic $2$-vectors (\textit{cfr.} \cite{Koecher,Loos,Gun-Bars} and
\cite{GST}). On the other hand, $\mathcal{K}\left( J_{3}^{\mathbb{O}}\right)
$ stands for the stabilizer of the rank-$4$ orbit of the action of $%
Conf\left( J_{3}^{\mathbb{O}}\right) \cong Aut\left( \mathbf{F}\left( J_{3}^{%
\mathbb{O}}\right) \right) $ on its (fundamental) irrep. $\mathbf{56}$ with
positive quartic invariant $I_{4}>0$ and representative \textquotedblleft $%
++--$" (for further detail, \textit{cfr. e.g.} \cite{BFGM1,CFM1,Small-Orbits}%
).

\md

\subsubsection{$\cp_{1}^{6(-14)}$\label{631}}

The maximal parabolics ~$\cp_{1}^{6(-14)}$ ~ from \eqref{maxsixc}
corresponds to the Bruhat decomposition:
\begin{equation}
E_{6(-14)}~=~\cn_{1}^{6(-14)-}\oplus so(7,1)\oplus so(2)\oplus so(1,1)\oplus %
\cn_{1}^{6(-14)+}\ ,\newline
\end{equation}%
which can be obtained by the following embedding chain:
\begin{eqnarray}
E_{6(-14)} &\supset &so(8,2)\oplus so(2)\supset so(7,1)\oplus so(2)\oplus
so(1,1);  \label{23} \\
\mathbf{78} &=&\mathbf{28}_{0,0}\oplus \mathbf{1}_{0,0}\oplus \mathbf{1}%
_{0,0}\oplus \mathbf{8}_{v,0,2}\oplus \mathbf{8}_{v,0,-2}\oplus \mathbf{8}%
_{c,-3,1}\oplus \mathbf{8}_{s,-3,-1}\oplus \mathbf{8}_{s,3,1}\oplus \mathbf{8%
}_{c,3,-1},  \label{24}
\end{eqnarray}%
thus yielding a $5$-grading with ~$\mathcal{N}_{1}^{6(-14)+}=\mathbf{8}%
_{v,0,2}+\mathbf{8}_{c,-3,1}+\mathbf{8}_{s,3,1}\,$.

A Jordan algebraic interpretation (of the first step) of the chain (\ref{23}%
)-(\ref{24}) is a consequence of the Jordan algebraic embedding $J_{3}^{%
\mathbb{O}}\supset \mathbb{R}\oplus J_{2}^{\mathbb{O}}$ at the level of the $%
\mathcal{K}$-symmetry, namely :%
\begin{equation}
\mathcal{K}\left( J_{3}^{\mathbb{O}}\right) \supset \mathcal{K}\left(
\mathbb{R}\oplus J_{2}^{\mathbb{O}}\right) \oplus \mathcal{A}_{8},
\label{j631}
\end{equation}%
where $\mathcal{K}\left( \mathbb{R}\oplus J_{2}^{\mathbb{O}}\right) $
denotes stabilizer of the rank-$4$ orbit of the action of \linebreak $%
Conf\left( \mathbb{R}\oplus J_{2}^{\mathbb{O}}\right) \cong Aut\left(
\mathbf{F}\left( \mathbb{R}\oplus J_{2}^{\mathbb{O}}\right) \right) $ on its
(bi-fundamental) rep. $\left( \mathbf{2},\mathbf{10}\right) $ with positive
quartic invariant $I_{4}>0$ and representative \textquotedblleft $++--$"
\cite{BFGM1,CFM1,Small-Orbits}). Indeed, it holds that (\textit{cfr.}
\textit{e.g.} case $4b$ with $n=10$ in Table VIII of \cite{Small-Orbits})%
\begin{equation}
so(8,2)\oplus so(2)\cong \mathcal{K}\left( \mathbb{R}\oplus J_{2}^{\mathbb{O}%
}\right) ,
\end{equation}%
which is nothing but the $q=8$ case of the general relation :%
\begin{equation}
so(q,2)\oplus so(2)\cong \mathcal{K}\left( \mathbb{R}\oplus J_{2}^{\mathbb{A}%
}\right) \ .
\end{equation}

\md

\subsubsection{$\cp_{2}^{6(-14)}$\label{632}}

The maximal parabolics ~$\cp_{2}^{6(-14)}$ ~ from \eqref{maxsixc}
corresponds to the Bruhat decomposition:
\begin{equation}
E_{6(-14)}~=~\cn_{2}^{6(-14)-}\oplus su(5,1)\oplus so(1,1)\oplus \cn%
_{2}^{6(-14)+}\ ,
\end{equation}%
which can be obtained through the following embedding chain:
\begin{eqnarray}
E_{6(-14)} &\supset &su(5,1)\oplus sl(2,\mathbb{R})\supset su(5,1)\oplus
so(1,1);  \label{25} \\
\mathbf{78} &=&\left( \mathbf{35},\mathbf{1}\right) \oplus \left( \mathbf{1},%
\mathbf{3}\right) \oplus \left( \mathbf{20},\mathbf{2}\right) =\mathbf{35}%
_{0}+\mathbf{1}_{0}+\mathbf{1}_{2}+\mathbf{1}_{-2}+\mathbf{20}_{1}+\mathbf{20%
}_{-1},  \label{26}
\end{eqnarray}%
exhibiting a $5$-grading with~ $\mathcal{N}_{2}^{6(-14)\pm }=\mathbf{1}_{\pm
2}+\mathbf{20}_{\pm 1}\,$.

A possible Jordan algebraic interpretation (of the first step) of the chain (%
\ref{25})-(\ref{26}) would be based on the characterization of $su(5,1)$,
which, at the moment, we can characterize as a maximal subalgebra (with
commutant $sl(2,\mathbb{R})$) of $E_{6(-14)}$, only. As mentioned above,work
is in progress for a sharper interpretation of $su(5,1)$.

\md

\subsection{$E_{6(-26)}$}

This is the minimally non-compact, real form of $E_{6}$. Its Jordan
interpretation reads:%
\begin{equation}
E_{6(-26)}\cong str_{0}\left( J_{3}^{\mathbb{O}}\right) \cong der\left(
J_{3}^{\mathbb{O}},J_{3}^{\mathbb{O}\prime }\right) \ominus so(1,1).
\end{equation}

%\subsubsection{\label{}}
There is only one maximal parabolics ~$\cp^{6(-26)}$ ~ from %
\eqref{maxsixd} corresponding to the Bruhat decomposition:
\begin{equation}
E_{6(-26)}~=~\cn^{6(-26)-}\oplus so(9,1)\oplus so(1,1)\oplus \cn^{6(-26)+}\ ,%
\newline
\label{28}
\end{equation}%
thus yielding a $3$-grading. Since%
\begin{equation}
so(9,1)\cong str_{0}\left( J_{2}^{\mathbb{O}}\right) ,
\end{equation}%
\textit{at least} two Jordan algebraic interpretations (denoted by $1$ and $%
2 $) of (\ref{28}) can be given :

\begin{enumerate}
\item
\begin{equation}
1:str_{0}\left( J_{3}^{\mathbb{O}}\right) \supset str_{0}\left( J_{2}^{%
\mathbb{O}}\right) \oplus so(1,1),  \label{64I}
\end{equation}%
where the $so(1,1)$ generating the $3$-grading is the KK $so(1,1)$ of the $%
S^{1}$-reduction $D=6\rightarrow 5$.

\item The second interpretation stems from the embedding $\mathfrak{J}_{3}^{%
\mathbb{O}}\supset \mathbb{R}\oplus \mathfrak{J}_{2}^{\mathbb{O}}$,
evaluated at $str_{0}$ level :%
\begin{equation}
2:str_{0}\left( J_{3}^{\mathbb{O}}\right) \supset str_{0}\left( \mathbb{R}%
\oplus J_{2}^{\mathbb{O}}\right) \oplus \mathcal{A}_{8},  \label{64II}
\end{equation}%
where the $so(1,1)$ generating the $3$-grading is the dilatonic $so(1,1)$ in
$D=5$.
\end{enumerate}

\md

\subsection{$E_{7(7)}$}

This is the \textit{split} real form of $E_{7}$. Its Jordan interpretation
is essentially twofold (due to the symmetry of the double-split Magic Square
${\mathcal{L}}_{3}(\mathbb{A}_{s},\mathbb{B}_{s})$ \cite{BS}, reported in
Table 3) :%
\begin{eqnarray}
E_{7(7)} &\cong &qconf\left( J_{3}^{\mathbb{H}_{s}}\right)  \label{E7(7)-1}
\\
&\cong &conf\left( J_{3}^{\mathbb{O}_{s}}\right) \cong der\left( \mathbf{F}%
\left( J_{3}^{\mathbb{O}_{s}}\right) \right) .  \label{E7(7)-2}
\end{eqnarray}

\subsubsection{$\cp_{1}^{7(7)}$\label{711}}

The maximal parabolics ~$\cp_{1}^{7(7)}$ ~ from \eqref{maxseva} corresponds
to the Bruhat decomposition:
\begin{equation}
E_{7(7)}~=~\cn_{1}^{7(7)-}\oplus so(6,6)\oplus so(1,1)\oplus \cn%
_{1}^{7(7)+}\ ,\newline
\end{equation}%
which can be obtained through the following embedding chain :
\begin{eqnarray}
E_{7(7)} &\supset &so(6,6)\oplus sl(2,\mathbb{R})\supset so(6,6)\oplus
so(1,1);  \label{29} \\
\mathbf{133} &=&\left( \mathbf{66,1}\right) \oplus \left( \mathbf{1},\mathbf{%
3}\right) \oplus \left( \mathbf{32}^{\prime }\mathbf{,2}\right) =\mathbf{66}%
_{0}\oplus \mathbf{1}_{0}\oplus \mathbf{1}_{2}\oplus \mathbf{1}_{-2}\oplus
\mathbf{32}_{1}^{\prime }\oplus \mathbf{32}_{-1}^{\prime },  \label{30}
\end{eqnarray}%
thus yielding a $5$-grading as $\mathcal{N}_{1}^{7(7)+}=\mathbf{1}_{2}+%
\mathbf{32}_{1}^{\prime }$.

The Jordan algebraic interpretation of (the first step of) (\ref{29})-(\ref%
{30}) is \textit{at least} threefold, and it is based on the following
identifications :%
\begin{eqnarray}
so(6,6)\oplus sl(2,\mathbb{R}) &\cong &conf\left( \mathbb{R}\oplus J_{2}^{%
\mathbb{O}_{s}}\right) \cong der\left( \mathbf{F}\left( \mathbb{R}\oplus
J_{2}^{\mathbb{O}_{s}}\right) \right) ;  \label{uno-1} \\
so(6,6) &\cong &conf\left( J_{3}^{\mathbb{H}_{s}}\right) \cong der\left(
\mathbf{F}\left( J_{3}^{\mathbb{H}_{s}}\right) \right)  \label{uno-pre-2} \\
&\cong &qconf\left( \mathbb{R}\oplus J_{2}^{\mathbb{H}_{s}}\right) .
\label{uno-2}
\end{eqnarray}

\begin{enumerate}
\item The first interpretation stems from the embedding $\mathbb{O}%
_{s}\supset \mathbb{H}_{s}$, at $conf$ level :%
\begin{equation}
conf\left( J_{3}^{\mathbb{O}_{s}}\right) \supset conf\left( J_{3}^{\mathbb{H}%
_{s}}\right) \oplus \widetilde{\mathcal{A}}_{4};  \label{j711I}
\end{equation}

\item the second one stems from the embedding $J_{3}^{\mathbb{O}_{s}}\supset
\mathbb{R}\oplus J_{2}^{\mathbb{O}_{s}}$ at the $conf$ level :%
\begin{equation}
conf\left( J_{3}^{\mathbb{O}_{s}}\right) \supset conf\left( \mathbb{R}\oplus
J_{2}^{\mathbb{O}_{s}}\right) \oplus \widetilde{\mathcal{A}}_{8};
\label{j711II}
\end{equation}

\item the third interpretation stems from $J_{3}^{\mathbb{H}_{s}}\supset
\mathbb{R}\oplus J_{2}^{\mathbb{H}_{s}}$, at the $qconf$ level :%
\begin{equation}
qconf\left( J_{3}^{\mathbb{H}_{s}}\right) \supset qconf\left( \mathbb{R}%
\oplus J_{2}^{\mathbb{H}_{s}}\right) \oplus \widetilde{\mathcal{A}}_{4}.
\label{j711III}
\end{equation}
\end{enumerate}

\subsubsection{$\cp_{2}^{7(7)}$\label{712}}

The maximal parabolics ~$\cp_{2}^{7(7)}$ ~ from \eqref{maxseva} corresponds
to the Bruhat decomposition:
\begin{equation}
E_{7(7)}~=\cn_{2}^{7(7)-}\oplus ~\left( sl(6,\mathbb{R})\oplus sl(2,\mathbb{R%
})\right) \oplus so(1,1)\oplus \cn_{2}^{7(7)+}\ ,\newline
\end{equation}%
which can be obtained through at least two embedding chains, respectively
denoted by $I$ and $II$ :%
\begin{eqnarray}
I &:&E_{7(7)}\supset so(6,6)\oplus sl(2,\mathbb{R})\supset _{ii}sl(6,\mathbb{%
R})\oplus sl(2,\mathbb{R})\oplus so(1,1);  \label{37} \\
\mathbf{133} &=&\left( \mathbf{66,1}\right) \oplus \left( \mathbf{1},\mathbf{%
3}\right) \oplus \left( \mathbf{32}^{\prime }\mathbf{,2}\right) =\left\{
\begin{array}{l}
\left( \mathbf{35,1}\right) _{0}\oplus \left( \mathbf{1,3}\right) _{0}\oplus
\left( \mathbf{1,1}\right) _{0}\oplus \left( \mathbf{1},\mathbf{2}\right)
_{3}\oplus (\mathbf{1},\mathbf{2})_{-3} \\
\oplus \left( \mathbf{15},\mathbf{2}\right) _{-1}\oplus \left( \mathbf{15}%
^{\prime },\mathbf{2}\right) _{1}\oplus \left( \mathbf{15,1}\right)
_{2}\oplus \left( \mathbf{15}^{\prime }\mathbf{,1}\right) _{-2}.%
\end{array}%
\right.  \label{38}
\end{eqnarray}%
\begin{eqnarray}
II &:&E_{7(7)}\supset ^{ns}sl(6,\mathbb{R})\oplus sl(3,\mathbb{R})\supset
sl(6,\mathbb{R})\oplus sl(2,\mathbb{R})\oplus so(1,1);  \label{39} \\
\mathbf{133} &=&\left( \mathbf{35,1}\right) \oplus \left( \mathbf{1},\mathbf{%
8}\right) \oplus \left( \mathbf{15,3}^{\prime }\right) \oplus \left( \mathbf{%
15}^{\prime }\mathbf{,3}\right)  \notag \\
&=&\left\{
\begin{array}{l}
\left( \mathbf{35,1}\right) _{0}\oplus \left( \mathbf{1,3}\right) _{0}\oplus
\left( \mathbf{1,1}\right) _{0}\oplus \left( \mathbf{1},\mathbf{2}\right)
_{3}\oplus (\mathbf{1},\mathbf{2})_{-3} \\
\oplus \left( \mathbf{15},\mathbf{2}\right) _{-1}\oplus \left( \mathbf{15}%
^{\prime },\mathbf{2}\right) _{1}\oplus \left( \mathbf{15,1}\right)
_{2}\oplus \left( \mathbf{15}^{\prime }\mathbf{,1}\right) _{-2}.%
\end{array}%
\right.  \label{40}
\end{eqnarray}%
The chains of embeddings $I$ and $II$ give rise to a $7$-grading with $%
\mathcal{N}_{2}^{7(7)+}=\left( \mathbf{1},\mathbf{2}\right) _{3}+\left(
\mathbf{15,1}\right) _{2}+\left( \mathbf{15}^{\prime },\mathbf{2}\right)
_{1}\,$.
%\begin{comment}
%Note that the 2-dimensional subspaces  ~$(\mathbf{1,2})_{\pm 3}$~
% of ~$\cn^{7(7)\pm}_2$~ are the zero-value eigenspaces of  ~$sl(6,\bbr)$~ factor of
%~$\cm^{7(7)}_2\,$. ~The 15-dimensional subspaces ~$(\mathbf{15^{(')},1})_{\pm 2}$~
%are the zero-value eigenspaces of the ~$sl(2,\bbr)$~ factor of
%~$\cm^{7(7)}_4\,$.
%The remaining 30-dimensional subspaces ~$(\mathbf{15^{(')},2})_{\pm 1}$\\  may be
%split in two subrepresentations which are two copies of
%a representation of ~$sl(6,\bbr)$. They are distinguished by the eigenvalue of ~$\a_1\,$~:
%it is ~$+1$~ on one subrepresentation and  ~$-1$~ on the other.
%\end{comment}
It is here worth remarking that the second steps of the chain $%
I $ do pertain to two \textit{different}, \textit{inequivalent} (maximal,
symmetric) embeddings of $gl(6,\mathbb{R})$ into $so(6,6)$, respectively
denoted by $i$ and $ii$; such two embeddings can \textit{e.g.} be
discriminated by the branching of the chiral spinor irreps. $\mathbf{32}$
and $\mathbf{32}^{\prime }$ of $so(6,6)$, namely\footnote{%
It is amusing to note that the branching $ii$ is overlooked in the otherwise
fairly comprehensive treatment \textit{e.g.} of \cite{Slansky} and \cite%
{Patera}; it is however considered \textit{e.g.} in \cite{Minchenko}.}:%
\begin{eqnarray}
i &:&so(6,6)\supset sl(6,\mathbb{R})\oplus so(1,1):\left\{
\begin{array}{l}
\mathbf{32}=\mathbf{20}_{0}\oplus \mathbf{6}_{-2}\oplus \mathbf{6}%
_{2}^{\prime }, \\
\mathbf{32}^{\prime }=\mathbf{15}_{-1}\oplus \mathbf{15}_{1}^{\prime }\oplus
\mathbf{1}_{3}\oplus \mathbf{1}_{-3};%
\end{array}%
\right.  \label{i} \\
ii &:&so(6,6)\supset sl(6,\mathbb{R})\oplus so(1,1):\left\{
\begin{array}{l}
\mathbf{32}=\mathbf{15}_{-1}\oplus \mathbf{15}_{1}^{\prime }\oplus \mathbf{1}%
_{3}\oplus \mathbf{1}_{-3}, \\
\mathbf{32}^{\prime }=\mathbf{20}_{0}\oplus \mathbf{6}_{-2}\oplus \mathbf{6}%
_{2}^{\prime }.%
\end{array}%
\right.  \label{ii}
\end{eqnarray}

Let us now consider the Jordan algebraic interpretation of the various
chains. We start by observing that%
\begin{equation}
gl(6,\mathbb{R})\cong str\left( J_{3}^{\mathbb{H}_{s}}\right) \cong
conf\left( J_{3}^{\mathbb{C}_{s}}\right) \oplus so(1,1);  \label{gg}
\end{equation}%
By also recalling (\ref{uno-1}), (\ref{uno-2}) and (\ref{uno-4}), possible
interpretations read as follows:%
\begin{equation}
I:conf\left( J_{3}^{\mathbb{O}_{s}}\right) \overset{J_{3}^{\mathbb{O}%
_{s}}\supset \mathbb{R}\oplus J_{2}^{\mathbb{O}_{s}}}{\supset }conf\left(
\mathbb{R}\oplus J_{2}^{\mathbb{O}_{s}}\right) \oplus \widetilde{\mathcal{A}}%
_{8}\supset _{i,ii}conf\left( J_{3}^{\mathbb{C}_{s}}\right) \oplus
\widetilde{\mathcal{A}}_{2}\oplus sl\left( 2,\mathbb{R}\right) ;
\end{equation}%
\begin{equation}
II:qconf\left( J_{3}^{\mathbb{H}_{s}}\right) \supset ^{ns}str_{0}\left(
J_{3}^{\mathbb{H}_{s}}\right) \oplus sl_{Ehlers}(3,\mathbb{R})\supset
str\left( J_{3}^{\mathbb{H}_{s}}\right) \oplus sl(2,\mathbb{R})\oplus
so(1,1).
\end{equation}

\md

\subsubsection{$\cp_{3}^{7(7)}$\label{713}}

The maximal parabolics ~$\cp_{3}^{7(7)}$ ~ from \eqref{maxseva} corresponds
to the Bruhat decomposition:
\begin{equation}
E_{7(7)}~=~\cn_{3}^{7(7)-}\oplus \left( sl(4,\bbr)\oplus sl(3,\bbr)\oplus
sl(2,\bbr)\right) \ \oplus so(1,1)\oplus \cn_{3}^{7(7)+}\ ,\newline
\end{equation}%
which can be obtained by \textit{at least} two chains of embeddings\footnote{%
We use the conventions of \cite{Kephart}.}, respectively denoted by $1$ and%
\footnote{%
In the second step of chain $2$ there is complete symmetry between the two $%
so(3,3)$ factors, thus the choice of which one is branched in the third step
is immaterial.} $2$ (recall that $so(3,3)\cong sl(4,\mathbb{R})$) :%
\begin{eqnarray}
1 &:&E_{7(7)}\supset ^{ns}sl(6,\mathbb{R})\oplus sl(3,\mathbb{R})\supset
sl(4,\mathbb{R})\oplus sl(3,\mathbb{R})\oplus sl(2,\mathbb{R})\oplus so(1,1);
\\
\mathbf{133} &=&\left( \mathbf{35,1}\right) \oplus \left( \mathbf{1},\mathbf{%
8}\right) \oplus \left( \mathbf{15,3}^{\prime }\right) \oplus \left( \mathbf{%
15}^{\prime }\mathbf{,3}\right)  \notag \\
&=&\left( \mathbf{15,1,1}\right) _{0}\oplus \left( \mathbf{1,1,3}\right)
_{0}\oplus \left( \mathbf{1,1,1}\right) _{0}\oplus \left( \mathbf{4,1,2}%
\right) _{3}\oplus \left( \mathbf{4}^{\prime }\mathbf{,1,2}\right) _{-3}
\notag \\
&&\oplus \left( \mathbf{1},\mathbf{8,1}\right) _{0}\oplus \left( \mathbf{1,3}%
^{\prime },\mathbf{1}\right) _{-4}\oplus \left( \mathbf{4,3}^{\prime },%
\mathbf{2}\right) _{-1}\oplus \left( \mathbf{6,3}^{\prime },\mathbf{1}%
\right) _{2}  \notag \\
&&\oplus \left( \mathbf{1,3},\mathbf{1}\right) _{4}\oplus \left( \mathbf{4}%
^{\prime }\mathbf{,3},\mathbf{2}\right) _{1}\oplus \left( \mathbf{6,3},%
\mathbf{1}\right) _{-2};
\end{eqnarray}%
\begin{eqnarray}
2 &:&E_{7(7)}\supset so(6,6)\oplus sl(2,\mathbb{R})\supset so(3,3)_{\mathbf{I%
}}\oplus so(3,3)_{\mathbf{II}}\oplus sl(2,\mathbb{R})  \notag \\
&\supset &sl(4,\mathbb{R})\oplus sl(3,\mathbb{R})\oplus sl(2,\mathbb{R}%
)\oplus so(1,1); \\
\mathbf{133} &=&\left( \mathbf{66,1}\right) \oplus \left( \mathbf{1},\mathbf{%
3}\right) \oplus \left( \mathbf{32}^{\prime }\mathbf{,2}\right) =\left(
\mathbf{15,1,1}\right) \oplus \left( \mathbf{1,15,1}\right) \oplus \left(
\mathbf{6,6,1}\right) \oplus \left( \mathbf{1},\mathbf{1},\mathbf{3}\right)
\oplus \left( \mathbf{4,4}^{\prime }\mathbf{,2}\right) \oplus \left( \mathbf{%
4}^{\prime }\mathbf{,4,2}\right)  \notag \\
&=&\left( \mathbf{15,1,1}\right) _{0}\oplus \left( \mathbf{1,8,1}\right)
_{0}\oplus \left( \mathbf{1,1,1}\right) _{0}\oplus \left( \mathbf{1,3,1}%
\right) _{4}\oplus \left( \mathbf{1,3}^{\prime }\mathbf{,1}\right) _{-4}
\notag \\
&&\oplus \left( \mathbf{6,3,1}\right) _{-2}\oplus \left( \mathbf{6,3}%
^{\prime }\mathbf{,1}\right) _{2}\oplus \left( \mathbf{1},\mathbf{1},\mathbf{%
3}\right) _{0}  \notag \\
&&\oplus \left( \mathbf{4,3}^{\prime }\mathbf{,2}\right) _{-1}\oplus \left(
\mathbf{4,1,2}\right) _{3}\oplus \left( \mathbf{4}^{\prime }\mathbf{,3,2}%
\right) _{1}\oplus \left( \mathbf{4}^{\prime }\mathbf{,1,2}\right) _{-3}.
\end{eqnarray}%
These chains of embeddings give rise to a $9$-grading, with $\mathcal{N}%
_{3\,;\,a)}^{7(7)+}\cong \mathcal{N}_{3\,;\,b)}^{7(7)+}=\left( \mathbf{1,3},%
\mathbf{1}\right) _{4}\oplus \left( \mathbf{4,1,2}\right) _{3}\oplus \left(
\mathbf{6,3}^{\prime },\mathbf{1}\right) _{2}\oplus \left( \mathbf{4}%
^{\prime }\mathbf{,3},\mathbf{2}\right) _{1}\,$.

A Jordan algebraic interpretation of the chains $1$ and $2$ is based on the
following identification :%
\begin{equation}
gl(4,\mathbb{R})\cong so(3,3)\oplus so(1,1)\cong str_{0}\left( \mathbb{%
R\oplus }J_{2}^{\mathbb{H}_{s}}\right) ,
\end{equation}%
and it reads%
\begin{equation}
1:\left.
\begin{array}{l}
qconf\left( J_{3}^{\mathbb{H}_{s}}\right) \supset ^{ns}str_{0}\left( J_{3}^{%
\mathbb{H}_{s}}\right) \oplus sl(3,\mathbb{R})_{Ehlers} \\
\text{\textit{or}} \\
conf\left( J_{3}^{\mathbb{O}_{s}}\right) \supset ^{ns}conf\left( J_{3}^{%
\mathbb{C}_{s}}\right) \oplus sl(3,\mathbb{R})%
\end{array}%
\right\} \supset str_{0}(\mathbb{R\oplus }J_{2}^{\mathbb{H}_{s}})\oplus sl(3,%
\mathbb{R})_{Ehlers}\oplus \widetilde{\mathcal{A}}_{4};  \label{j713a}
\end{equation}%
\begin{eqnarray}
2 &:&qconf\left( J_{3}^{\mathbb{H}_{s}}\right) \supset qconf\left( \mathbb{%
R\oplus }J_{2}^{\mathbb{H}_{s}}\right) \oplus \widetilde{\mathcal{A}}_{q=4}
\notag \\
&\supset &so(3,3)_{\mathbf{I}}\oplus so(3,3)_{\mathbf{II}}\oplus sl(2,%
\mathbb{R})\supset str_{0}\left( \mathbb{R\oplus }J_{2}^{\mathbb{H}%
_{s}}\right) \oplus sl(3,\mathbb{R})_{Ehlers}\oplus \widetilde{\mathcal{A}}%
_{4}.  \label{j713b}
\end{eqnarray}%
The second step of the Jordan algebraic interpretation (\ref{j713b}) of the
chain $2$ is intentionally left explicit : intriguingly, it provides an
\textit{enhancement} of the symmetry obtained by the JP embedding for the $%
\left( \mathbb{R\oplus }J_{2}^{\mathbb{H}_{s}}\right) $-based theory :
indeed, this latter reads%
\begin{equation}
qconf\left( J_{3}^{\mathbb{H}_{s}}\right) \supset str_{0}\left( \mathbb{%
R\oplus }J_{2}^{\mathbb{H}_{s}}\right) \oplus sl(3,\mathbb{R})_{Ehlers}\cong
str_{0}\left( J_{2}^{\mathbb{H}_{s}}\right) \oplus sl(3,\mathbb{R}%
)_{Ehlers}\times so(1,1),
\end{equation}%
and then the following enhancement takes place:%
\begin{equation}
sl(3,\mathbb{R})_{Ehlers}\oplus so(1,1)\rightarrow sl(4,\mathbb{R})\cong
so(3,3),  \label{enh-1}
\end{equation}%
where the $so(3,3)$ on the right-hand side is one of the two summands in the
second step of the chain $2$. The enhancement (\ref{enh-1}) actually hints
for another interpretation of the \textit{non-symmetric} embedding%
\begin{equation}
qconf\left( J_{3}^{\mathbb{H}_{s}}\right) \supset _{nm}^{ns}so(3,3)_{\mathbf{%
I}}\oplus so(3,3)_{\mathbf{II}}\oplus sl(2,\mathbb{R}),  \label{enh-2}
\end{equation}%
where the subscript \textquotedblleft $nm$" denotes\ its \textit{%
next-to-maximal} nature (\textit{i.e.}, the fact that it is realized by two
subsequent maximal embeddings, namely the first two steps of chain $2$). In
fact, (\ref{enh-2}) can be interpreted as the $D=6$ case of the Ehlers
embedding for the non-supersymmetric $\mathbb{H}_{s}$-based gravity theory
(see \cite{Magic-Non-Susy} and Refs. therein), and its
next-to-maximal nature is consistent with the treatment of \cite%
{super-Ehlers}; in other words, (\ref{enh-2}) can be interpreted as follows:%
\begin{equation}
qconf\left( J_{3}^{\mathbb{H}_{s}}\right) \supset _{nm}^{ns}str_{0}\left(
J_{2}^{\mathbb{H}_{s}}\right) \oplus \widetilde{\mathcal{A}}_{4}\oplus sl(4,%
\mathbb{R})_{Ehlers},  \label{enh-3}
\end{equation}%
where therefore one of the two $so(3,3)$ factors (say, $so(3,3)_{\mathbf{I}}$%
) is conceived as $str_{0}\left( J_{2}^{\mathbb{H}_{s}}\right) $, while the
other one (say, $so(3,3)_{\mathbf{II}}$) is nothing but the $D=6$ Ehlers
symmetry, and the commuting $sl(2,\mathbb{R})$ factor is seen as $\widetilde{%
\mathcal{A}}_{4}$. As a consequence of this interpretation, (\ref{j713b})
can concisely be rewritten as :%
\begin{eqnarray}
2 &:&qconf\left( J_{3}^{\mathbb{H}_{s}}\right) \supset qconf\left( \mathbb{%
R\oplus }J_{2}^{\mathbb{H}_{s}}\right) \oplus \widetilde{\mathcal{A}}_{4}
\notag \\
&\supset &str_{0}\left( J_{2}^{\mathbb{H}_{s}}\right) \oplus \widetilde{%
\mathcal{A}}_{4}\oplus sl(4,\mathbb{R})_{Ehlers}  \notag \\
&\supset &str_{0}\left( \mathbb{R\oplus }J_{2}^{\mathbb{H}_{s}}\right)
\oplus sl(3,\mathbb{R})_{Ehlers}\oplus \widetilde{\mathcal{A}}_{4},
\end{eqnarray}%
where%
\begin{equation}
str_{0}\left( \mathbb{R\oplus }J_{2}^{\mathbb{H}_{s}}\right) \cong
so(1,1)\oplus str_{0}\left( J_{2}^{\mathbb{H}_{s}}\right) .
\end{equation}

\md

\subsubsection{$\cp_{4}^{7(7)}$\label{714}}

The maximal parabolics ~$\cp_{4}^{7(7)}$ ~ from \eqref{maxseva} corresponds
to the Bruhat decomposition:
\begin{equation}
E_{7(7)}~=\left( \mathcal{N}_{4}^{-}\right) \oplus ~(sl(5,\bbr)\oplus sl(3,%
\bbr))\ \oplus so(1,1)\oplus \left( \mathcal{N}_{4}^{+}\right) \newline
\mathbf{,}
\end{equation}%
which can be obtained through \textit{at least} two chains of embeddings,
respectively denoted by $1$ and $2$ :%
\begin{eqnarray}
1 &:&E_{7(7)}\supset ^{ns}sl(6,\mathbb{R})\oplus sl(3,\mathbb{R})\supset
sl(5,\mathbb{R})\oplus sl(3,\mathbb{R})\oplus so(1,1); \\
\mathbf{133} &=&\left( \mathbf{35,1}\right) \oplus \left( \mathbf{1},\mathbf{%
8}\right) \oplus \left( \mathbf{15,3}^{\prime }\right) \oplus \left( \mathbf{%
15}^{\prime }\mathbf{,3}\right)  \notag \\
&=&\left( \mathbf{24,1}\right) _{0}\oplus \left( \mathbf{1,1}\right) \oplus
\left( \mathbf{5,1}\right) _{6}\oplus \left( \mathbf{5}^{\prime }\mathbf{,1}%
\right) _{-6}\oplus \left( \mathbf{1},\mathbf{8}\right) _{0}  \notag \\
&&\oplus \left( \mathbf{5,3}^{\prime }\right) _{-4}\oplus \left( \mathbf{10,3%
}^{\prime }\right) _{2}\oplus \left( \mathbf{5}^{\prime }\mathbf{,3}\right)
_{4}\oplus \left( \mathbf{10}^{\prime }\mathbf{,3}\right) _{-2};
\end{eqnarray}%
\begin{eqnarray}
2 &:&E_{7(7)}\supset sl(8,\mathbb{R})\supset sl(5,\mathbb{R})\oplus sl(3,%
\mathbb{R})\oplus so(1,1); \\
\mathbf{133} &=&\mathbf{63}\oplus \mathbf{70}=(\mathbf{24},\mathbf{1}%
)_{0}\oplus (\mathbf{1},\mathbf{8})_{0}\oplus (\mathbf{1},\mathbf{1}%
)_{0}\oplus (\mathbf{5},\mathbf{3}^{\prime })_{8}\oplus (\mathbf{5}^{\prime
},\mathbf{3})_{-8}  \notag \\
&&\oplus (\mathbf{5},\mathbf{1})_{-12}\oplus (\mathbf{5}^{\prime },\mathbf{1}%
)_{12}\oplus (\mathbf{10},\mathbf{3}^{\prime })_{-4}\oplus (\mathbf{10}%
^{\prime },\mathbf{3})_{4}.
\end{eqnarray}%
Both chains $1$ and $2$ give rise to a $7$-grading, with $\mathcal{N}%
_{4\,;\,a)}^{7(7)+}=\left( \mathbf{5,1}\right) _{6}\oplus \left( \mathbf{5}%
^{\prime }\mathbf{,3}\right) _{4}\oplus \left( \mathbf{10,3}^{\prime
}\right) _{2}$ and $\mathcal{N}_{4\,;\,b)}^{7(7)+}=(\mathbf{5}^{\prime },%
\mathbf{1})_{12}\oplus (\mathbf{5},\mathbf{3}^{\prime })_{8}\oplus (\mathbf{%
10}^{\prime },\mathbf{3})_{4}$ (both of real dimension $50$). Note that $%
1\leftrightarrow 2)$ \textit{iff} the weights of the parabolic $so(1,1)_{1}$
gets doubled, and \textit{iff} $\mathbf{5}\leftrightarrow \mathbf{5}^{\prime
}$ (yielding $\mathbf{10}\leftrightarrow \mathbf{10}^{\prime }$) and $%
\mathbf{3}\leftrightarrow \mathbf{3}^{\prime }$ in $sl(5,\mathbb{R})$ and $%
sl(3,\mathbb{R})$, respectively. A Jordan algebraic interpretation of the
first step of chains $1$ and $2$ coincides with the one of chain $1$ of $\cp%
_{3}^{7(7)}$, and of chain $1$ of case $\cp_{2}^{7(7)}$, respectively.

\md

\subsubsection{$\cp_{5}^{7(7)}$\label{715}}

The maximal parabolics ~$\cp_{5}^{7(7)}$ ~ from \eqref{maxseva} corresponds
to the Bruhat decomposition:
\begin{equation}
E_{7(7)}~=~\mathcal{N}_{5}^{7(7)-}\oplus so(5,5)\oplus sl(2,\bbr)\ \oplus
so(1,1)\oplus \mathcal{N}_{5}^{7(7)}\newline
\mathbf{,}
\end{equation}%
which can be obtained through the embedding chain :
\begin{eqnarray}
E_{7(7)} &\supset &so(6,6)\oplus sl(2,\mathbb{R})\supset so(5,5)\oplus sl(2,%
\mathbb{R})\oplus so(1,1);  \label{5-1} \\
\mathbf{133} &=&\left( \mathbf{66,1}\right) \oplus \left( \mathbf{1},\mathbf{%
3}\right) \oplus \left( \mathbf{32}^{\prime }\mathbf{,2}\right) = \\
&=&\left( \mathbf{45,1}\right) _{0}\oplus (\mathbf{1},\mathbf{1})_{0}\oplus (%
\mathbf{10},\mathbf{1})_{2}\oplus (\mathbf{10},\mathbf{1})_{-2}\oplus \left(
\mathbf{1},\mathbf{3}\right) _{0}\oplus \left( \mathbf{16,2}\right)
_{-1}\oplus \left( \mathbf{16}^{\prime }\mathbf{,2}\right) _{1},  \notag
\label{5-2}
\end{eqnarray}%
giving rise to a $5$-grading, with $\mathcal{N}_{5}^{7(7)}=(\mathbf{10},%
\mathbf{1})_{2}\oplus \left( \mathbf{16}^{\prime }\mathbf{,2}\right) _{1}\,$%
. Since%
\begin{equation}
so(5,5)\oplus so(1,1)\cong str_{0}\left( \mathbb{R}\oplus J_{2}^{\mathbb{O}%
_{s}}\right) \cong str\left( J_{2}^{\mathbb{O}_{s}}\right) ,
\end{equation}%
\textit{at least} two Jordan algebraic interpretations (denoted by $I$ and $%
II$) of the chain above can be given, namely :%
\begin{equation}
I:J_{3}^{\mathbb{O}_{s}}\supset \mathbb{R}\oplus J_{2}^{\mathbb{O}%
_{s}}\Rightarrow conf\left( J_{3}^{\mathbb{O}_{s}}\right) \supset conf\left(
\mathbb{R}\oplus J_{2}^{\mathbb{O}_{s}}\right) \oplus \widetilde{\mathcal{A}}%
_{8}\supset str\left( J_{2}^{\mathbb{O}_{s}}\right) \oplus sl(2,\mathbb{R}%
)\oplus \widetilde{\mathcal{A}}_{8};  \label{j715I}
\end{equation}%
\begin{equation}
II:qconf\left( J_{3}^{\mathbb{H}_{s}}\right) \supset conf\left( J_{3}^{%
\mathbb{H}_{s}}\right) \oplus sl(2,\mathbb{R})_{Ehlers}\supset str\left(
J_{2}^{\mathbb{O}_{s}}\right) \oplus sl(2,\mathbb{R})_{Ehlers}.
\label{j715II}
\end{equation}

\md

\subsubsection{$\cp_{6}^{7(7)}$\label{716}}

The maximal parabolics ~$\cp_{6}^{7(7)}$ ~ from \eqref{maxseva} corresponds
to the Bruhat decomposition:
\begin{equation}
E_{7(7)}~=~\mathcal{N}_{6}^{7(7)-}\oplus E_{6(6)}\oplus so(1,1)\oplus
\mathcal{N}_{6}^{7(7)+}\newline
\mathbf{,}
\end{equation}%
giving rise to a $3$-grading. \textit{At least} two Jordan algebraic
interpretations (denoted by $I$ and $II$) of the chain above can be given :%
\begin{equation}
I:\mathbb{H}_{s}\supset \mathbb{C}_{s}\Rightarrow qconf\left( J_{3}^{\mathbb{%
H}_{s}}\right) \supset qconf\left( J_{3}^{\mathbb{C}_{s}}\right) \oplus
\widetilde{\mathcal{A}}_{2};  \label{j716I}
\end{equation}%
\begin{equation}
II:conf\left( J_{3}^{\mathbb{O}_{s}}\right) \supset str\left( J_{3}^{\mathbb{%
O}_{s}}\right) \cong str_{0}\left( J_{3}^{\mathbb{O}_{s}}\right) \oplus
so(1,1)_{KK}.  \label{j716II}
\end{equation}%
It is worth remarking the different interpretation of the parabolic $so(1,1)$
(giving rise to the $3$-grading) in $I$ and $II$ : in $I$, it is identified
with $\widetilde{\mathcal{A}}_{2}$, whereas in $II$ it is the $so(1,1)_{KK}$
of the $S^{1}$-reduction $D=5\rightarrow 4$.

\md

\subsubsection{$\cp_{7}^{7(7)}$\label{717}}

The maximal parabolics ~$\cp_{7}^{7(7)}$ ~ from \eqref{maxseva} corresponds
to the Bruhat decomposition:
\begin{equation}
E_{7(7)}~=\mathcal{N}_{7}^{7(7)-}\oplus sl(7,\bbr)\oplus so(1,1)\oplus
\mathcal{N}_{7}^{7(7)+}\newline
\mathbf{,}
\end{equation}%
which can be obtained through the embedding chain
\begin{eqnarray}
E_{7(7)} &\supset &sl\left( 8,\mathbb{R}\right) \supset sl\left( 7,\mathbb{R}%
\right) \oplus so(1,1); \\
\mathbf{133} &=&\mathbf{63}\oplus \mathbf{70}=\mathbf{48}_{0}\oplus \mathbf{1%
}_{0}\oplus \mathbf{7}_{8}\oplus \mathbf{7}_{-8}^{\prime }\oplus \mathbf{35}%
_{-4}\oplus \mathbf{35}_{4}^{\prime }.
\end{eqnarray}%
Thus, it gives rise to a $5$-grading, with $\mathcal{N}_{7}^{7(7)+}=\mathbf{%
35}_{4}^{\prime }\oplus \mathbf{7}_{8}\,$.\newline

\textit{At least} two Jordan algebraic interpretations (denoted by $I$ and $%
II$) of the first step of the chain above can be given, respectively
pertaining to $\mathbb{H}_{s}$ and $\mathbb{O}_{s}$:%
\begin{equation}
I:qconf\left( J_{3}^{\mathbb{H}_{s}}\right) \supset \left. sl\left( q+4,%
\mathbb{R}\right) \right\vert _{q=4};  \label{I-one}
\end{equation}%
\begin{equation}
II:conf\left( J_{3}^{\mathbb{O}_{s}}\right) \supset \left. sl\left( q,%
\mathbb{R}\right) \right\vert _{q=8}.  \label{II-two}
\end{equation}

\md

\subsection{$E_{7(-5)}$}

The Jordan interpretation of\ $E_{7(-5)}$ reads:%
\begin{equation}
E_{7(-5)}\cong qconf\left( J_{3}^{\mathbb{H}}\right) .
\end{equation}

\md

\subsubsection{$\cp_{1}^{7(-5)}$\label{721}}

The maximal parabolics ~$\cp_{1}^{7(-5)}$ ~ from \eqref{maxsevb} corresponds
to the Bruhat decomposition:
\begin{equation}
E_{7(-5)}~=\mathcal{N}_{1}^{7(-5)-}\oplus ~so^{\ast }(12)\oplus
so(1,1)\oplus \mathcal{N}_{1}^{7(-5)}\newline
\mathbf{,}
\end{equation}%
which can be obtained through the embedding chain
\begin{eqnarray}
E_{7(-5)} &\supset &so^{\ast }(12)\oplus sl(2,\mathbb{R})\supset so^{\ast
}(12)\oplus so(1,1);  \label{1-1-1} \\
\mathbf{133} &=&\left( \mathbf{66,1}\right) \oplus \left( \mathbf{1},\mathbf{%
3}\right) \oplus \left( \mathbf{32}^{\prime }\mathbf{,2}\right) =\mathbf{66}%
_{0}\oplus \mathbf{1}_{0}\oplus \mathbf{1}_{2}\oplus \mathbf{1}_{-2}\oplus
\mathbf{32}_{1}^{\prime }\oplus \mathbf{32}_{-1}^{\prime },  \label{1-2-1}
\end{eqnarray}%
thus giving rise to a $5$-grading of contact type, with $\mathcal{N}%
_{1}^{7(-5)+}=\mathbf{1}_{2}\oplus \mathbf{32}_{1}^{\prime }$. By observing
that%
\begin{equation}
so^{\ast }(12)\cong conf\left( J_{3}^{\mathbb{H}}\right) \cong der\left(
\mathbf{F}\left( J_{3}^{\mathbb{H}}\right) \right) ,
\end{equation}%
a Jordan algebraic interpretation (of the first step) of the chain above is
given by the Ehlers embedding for the $\mathbb{H}$-based magic supergravity
theory :%
\begin{equation}
qconf\left( J_{3}^{\mathbb{H}}\right) \supset conf\left( J_{3}^{\mathbb{H}%
}\right) \oplus sl(2,\mathbb{R})_{Ehlers};  \label{j721}
\end{equation}%
thus, the parabolic $so(1,1)$ (determining the $5$-grading) is the
non-compact Cartan of the $D=4$ Ehlers $sl(2,\mathbb{R})_{Ehlers}$.

\md

\subsubsection{$\cp_{2}^{7(-5)}$\label{722}}

The maximal parabolics ~$\cp_{2}^{7(-5)}$ ~ from \eqref{maxsevb} corresponds
to the Bruhat decomposition:
\begin{equation}
E_{7(-5)}~=\mathcal{N}_{2}^{7(-5)-}\oplus ~so(7,3)\oplus su(2)\oplus
so(1,1)\oplus \mathcal{N}_{2}^{7(-5)+}\newline
\mathbf{,}
\end{equation}%
which can be obtained through the embedding chain
\begin{eqnarray}
E_{7(-5)} &\supset &so(8,4)\oplus su(2)\supset so(7,3)\oplus su(2)\oplus
so(1,1);  \label{5-1-1} \\
\mathbf{133} &=&\left( \mathbf{66,1}\right) \oplus \left( \mathbf{1},\mathbf{%
3}\right) \oplus \left( \mathbf{32}^{\prime }\mathbf{,2}\right)
\label{5-2-1} \\
&=&\left( \mathbf{45,1}\right) _{0}\oplus (\mathbf{1},\mathbf{1})_{0}\oplus (%
\mathbf{10},\mathbf{1})_{2}\oplus (\mathbf{10},\mathbf{1})_{-2}\oplus \left(
\mathbf{1},\mathbf{3}\right) _{0}\oplus \left( \mathbf{16,2}\right)
_{-1}\oplus \left( \mathbf{16}^{\prime }\mathbf{,2}\right) _{1}\ .  \notag
\end{eqnarray}%
Thus, it gives rise to a $5$-grading (\textit{cfr.} (\ref{5-1})-(\ref{5-2}%
)), with $\mathcal{N}_{5}^{7(-5)+}=(\mathbf{10},\mathbf{1})_{2}\oplus \left(
\mathbf{16}^{\prime }\mathbf{,2}\right) _{1}\,$. Since
\begin{equation}
so(8,4)\cong qconf\left( \mathbb{R}\oplus J_{2}^{\mathbb{H}}\right) ,
\end{equation}%
a Jordan algebraic interpretation (of the first step) of the chain above is
given by the embedding $J_{3}^{\mathbb{H}}\supset \mathbb{R}\oplus J_{2}^{%
\mathbb{H}}$ considered at $qconf$ level :%
\begin{equation}
qconf\left( J_{3}^{\mathbb{H}}\right) \supset conf\left( \mathbb{R}\oplus
J_{2}^{\mathbb{H}}\right) \oplus \mathcal{A}_{4}.  \label{j722}
\end{equation}

\md

\subsubsection{$\cp_{3}^{7(-5)}$\label{723}}

The maximal parabolics ~$\cp_{3}^{7(-5)}$ ~ from \eqref{maxsevb} corresponds
to the Bruhat decomposition:
\begin{equation}
E_{7(-5)}~=\mathcal{N}_{3}^{7(-5)-}\oplus ~su^{\ast }(6)\oplus sl(2,\bbr)\
\oplus so(1,1)\oplus \mathcal{N}_{3}^{7(-5)+}\newline
\mathbf{,}
\end{equation}%
which can be obtained \textit{at least} through two chains of embeddings,
respectively denoted by $1$ and $2$ :%
\begin{eqnarray}
1 &:&E_{7(-5)}\supset ^{ns}su^{\ast }(6)\oplus sl(3,\mathbb{R})\supset
su^{\ast }(6)\oplus sl(2,\mathbb{R})\oplus so(1,1); \\
\mathbf{133} &=&\left( \mathbf{35,1}\right) \oplus \left( \mathbf{1},\mathbf{%
8}\right) \oplus \left( \mathbf{15,3}^{\prime }\right) \oplus \left( \mathbf{%
15}^{\prime }\mathbf{,3}\right) \\
&=&\left\{
\begin{array}{l}
\left( \mathbf{35,1}\right) _{0}\oplus \left( \mathbf{1,3}\right) _{0}\oplus
\left( \mathbf{1,1}\right) _{0}\oplus \left( \mathbf{1},\mathbf{2}\right)
_{3}\oplus (\mathbf{1},\mathbf{2})_{-3} \\
\oplus \left( \mathbf{15},\mathbf{2}\right) _{-1}\oplus \left( \mathbf{15}%
^{\prime },\mathbf{2}\right) _{1}\oplus \left( \mathbf{15,1}\right)
_{2}\oplus \left( \mathbf{15}^{\prime }\mathbf{,1}\right) _{-2};%
\end{array}%
\right.  \notag
\end{eqnarray}%
\begin{eqnarray}
2 &:&E_{7(-5)}\supset so^{\ast }(12)\oplus sl(2,\mathbb{R})\supset
_{ii}su^{\ast }(6)\oplus sl(2,\mathbb{R})\oplus so(1,1); \\
\mathbf{133} &=&\left( \mathbf{66,1}\right) \oplus \left( \mathbf{1},\mathbf{%
3}\right) \oplus \left( \mathbf{32}^{\prime }\mathbf{,2}\right)  \notag \\
&=&\left\{
\begin{array}{l}
\left( \mathbf{35,1}\right) _{0}\oplus \left( \mathbf{1,3}\right) _{0}\oplus
\left( \mathbf{1,1}\right) _{0}\oplus \left( \mathbf{1},\mathbf{2}\right)
_{3}\oplus (\mathbf{1},\mathbf{2})_{-3} \\
\oplus \left( \mathbf{15},\mathbf{2}\right) _{-1}\oplus \left( \mathbf{15}%
^{\prime },\mathbf{2}\right) _{1}\oplus \left( \mathbf{15,1}\right)
_{2}\oplus \left( \mathbf{15}^{\prime }\mathbf{,1}\right) _{-2}.%
\end{array}%
\right.
\end{eqnarray}%
Both chains $1$ and $2$ give rise to a $7$-grading, with and $\ \mathcal{N}%
_{3}^{7(-5)+}=\left( \mathbf{1},\mathbf{2}\right) _{3}\oplus \left( \mathbf{%
15,1}\right) _{2}\oplus \left( \mathbf{15}^{\prime },\mathbf{2}\right)
_{1}\, $. Once again, it is here worth remarking that the second step of the
chain $2$ do pertain to two \textit{different}, \textit{inequivalent}
(maximal, symmetric) embeddings of $su^{\ast }(6,\mathbb{R})\times so(1,1)$
into $so^{\ast }(12)$, respectively denoted by $i$ and $ii$; such two
embeddings can \textit{e.g.} be discriminated by (a different non-compact,
real form of) the branching of the chiral spinor irreps. $\mathbf{32}$ and $%
\mathbf{32}^{\prime }$ of $so^{\ast }(12)$, \textit{cfr.} (\ref{i})-(\ref{ii}%
).

Let us now consider the Jordan algebraic interpretation of the various
chains. We observe that%
\begin{equation}
su^{\ast }(6)\cong str_{0}\left( J_{3}^{\mathbb{H}}\right) \cong der\left(
J_{3}^{\mathbb{H}},J_{3}^{\mathbb{H}\prime }\right) \ominus so(1,1).
\end{equation}%
Thus, Jordan interpretations can be given as follows:%
\begin{equation}
1:qconf\left( J_{3}^{\mathbb{H}}\right) \supset ^{ns}str\left( J_{3}^{%
\mathbb{H}}\right) \oplus sl_{Ehlers}(3,\mathbb{R})\supset str\left( J_{3}^{%
\mathbb{H}}\right) \oplus sl(2,\mathbb{R})\oplus so(1,1).
\end{equation}%
\begin{equation}
2:qconf\left( J_{3}^{\mathbb{H}}\right) \supset conf\left( J_{3}^{\mathbb{H}%
}\right) \oplus sl_{Ehlers}(2,\mathbb{R})\supset _{i,ii}str_{0}\left( J_{3}^{%
\mathbb{H}}\right) \oplus sl(2,\mathbb{R})_{Ehlers}\oplus so(1,1)_{KK}.
\end{equation}%
It is here worth commenting that the first step of the interpretation of $1$
is the JP\textit{\ }embedding for the $\mathbb{H}$-based theory (determining
the $D=5$ Ehlers $sl(3,\mathbb{R})_{Ehlers}$), and the resulting parabolic $%
so(1,1)$ (generating the $7$-grading) is the $so(1,1)$ commuting factor in
the right-hand side of the maximal, symmetric embedding $sl_{Ehlers}(3,%
\mathbb{R})\supset sl(2,\mathbb{R})\oplus so(1,1)$. On the other hand, the
first step of the interpretation of $2$ is the Ehlers embedding for the $%
\mathbb{H}$-based theory (determining the $D=4$ Ehlers $sl(2,\mathbb{R}%
)_{Ehlers}$ through the inverse $c^{\ast }$-map \cite{BGM,c-map}), and the
second step consists in an inverse $R^{\ast }$-map \cite{R-map}, which thus
introduces the $so(1,1)_{KK}$ of the $S^{1}$-reduction $D=5\rightarrow 4$;
in this case, this latter is the parabolic $so(1,1)$ (which generates the $7$%
-grading in $c)$).

\md

\subsubsection{$\cp_{4}^{7(-5)}$\label{724}}

The maximal parabolics ~$\cp_{4}^{7(-5)}$ ~ from \eqref{maxsevb} corresponds
to the Bruhat decomposition:
\begin{equation}
E_{7(-5)}~=~\mathcal{N}_{4}^{7(-5)-}\oplus so(5,1)\oplus sl(3,\bbr)\oplus
su(2)\oplus so(1,1)\oplus \mathcal{N}_{4}^{7(-5)+}\newline
\mathbf{,}
\end{equation}%
which can be obtained through \textit{at least} two chains of embeddings,
respectively denoted by $1$ and $2$ (recall that $su^{\ast }(4)\cong so(5,1)$%
, $su^{\ast }(2)\cong su(2)$) :%
\begin{eqnarray}
1 &:&E_{7(-5)}\supset ^{ns}su^{\ast }(6)\oplus sl(3,\mathbb{R})\supset
so(5,1)\oplus su(2)\oplus sl(3,\mathbb{R})\oplus so(1,1) \\
\mathbf{133} &=&\left( \mathbf{35,1}\right) \oplus \left( \mathbf{1},\mathbf{%
8}\right) \oplus \left( \mathbf{15,3}^{\prime }\right) \oplus \left( \mathbf{%
15}^{\prime }\mathbf{,3}\right)  \notag \\
&=&\left( \mathbf{15,1,1}\right) _{0}\oplus \left( \mathbf{1,3,1}\right)
_{0}\oplus \left( \mathbf{1,1,1}\right) _{0}\oplus \left( \mathbf{4,2,1}%
\right) _{3}\oplus \left( \mathbf{4}^{\prime }\mathbf{,2,1}\right) _{-3}
\notag \\
&&\oplus \left( \mathbf{1},\mathbf{1,8}\right) _{0}\oplus \left( \mathbf{1,1}%
,\mathbf{3}^{\prime }\right) _{-4}\oplus \left( \mathbf{4,2},\mathbf{3}%
^{\prime }\right) _{-1}\oplus \left( \mathbf{6,1},\mathbf{3}^{\prime
}\right) _{2}  \notag \\
&&\oplus \left( \mathbf{1,1},\mathbf{3}\right) _{4}\oplus \left( \mathbf{4}%
^{\prime }\mathbf{,2},\mathbf{3}\right) _{1}\oplus \left( \mathbf{6,1},%
\mathbf{3}\right) _{-2};
\end{eqnarray}%
\begin{eqnarray}
2 &:&E_{7(-5)}\supset so(8,4)\oplus su(2)  \notag \\
&\supset &so(5,1)\oplus su(2)\oplus sl(4,\mathbb{R})\supset so(5,1)\oplus
su(2)\oplus sl(3,\mathbb{R})\oplus so(1,1); \\
\mathbf{133} &=&\left( \mathbf{66,1}\right) \oplus \left( \mathbf{1},\mathbf{%
3}\right) \oplus \left( \mathbf{32}^{\prime }\mathbf{,2}\right)  \notag \\
&=&\left( \mathbf{15,1,1}\right) \oplus \left( \mathbf{1,15,1}\right) \oplus
\left( \mathbf{6,6,1}\right) \oplus \left( \mathbf{1},\mathbf{1},\mathbf{3}%
\right) \oplus \left( \mathbf{4,4}^{\prime }\mathbf{,2}\right) \oplus \left(
\mathbf{4}^{\prime }\mathbf{,4,2}\right)  \notag \\
&=&\left( \mathbf{15,1,1}\right) _{0}\oplus \left( \mathbf{1,3,1}\right)
_{0}\oplus \left( \mathbf{1,1,1}\right) _{0}\oplus \left( \mathbf{4,2,1}%
\right) _{3}\oplus \left( \mathbf{4}^{\prime }\mathbf{,2,1}\right) _{-3} \\
&&\oplus \left( \mathbf{1},\mathbf{1,8}\right) _{0}\oplus \left( \mathbf{1,1}%
,\mathbf{3}^{\prime }\right) _{-4}\oplus \left( \mathbf{4,2},\mathbf{3}%
^{\prime }\right) _{-1}\oplus \left( \mathbf{6,1},\mathbf{3}^{\prime
}\right) _{2}\oplus \left( \mathbf{1,1},\mathbf{3}\right) _{4}\oplus \left(
\mathbf{4}^{\prime }\mathbf{,2},\mathbf{3}\right) _{1}\oplus \left( \mathbf{%
6,1},\mathbf{3}\right) _{-2}.  \notag
\end{eqnarray}%
Note that $1$ and $2$ are different non-compact real forms of the chains $1$
and $2$ pertaining to $\cp_{3}^{7(7)}$\textbf{\ }(and the same holds for the
maximal parabolics under consideration). The chains of embeddings $1$ and $2$
give rise to a $7$-grading, with $\mathcal{N}_{4}^{7(-5)+}=\left( \mathbf{1,1%
},\mathbf{3}\right) _{4}\oplus \left( \mathbf{4,2,1}\right) _{3}\oplus
\left( \mathbf{6,1},\mathbf{3}^{\prime }\right) _{2}\oplus \left( \mathbf{4}%
^{\prime }\mathbf{,2},\mathbf{3}\right) _{1}\,$.

A Jordan algebraic interpretation of the chains $1$ and $2$ is based on the
following identification :%
\begin{equation}
so(5,1)\oplus so(1,1)\cong str_{0}\left( \mathbb{R\oplus }J_{2}^{\mathbb{H}%
}\right) ,
\end{equation}%
and it reads as follows:%
\begin{equation}
1:qconf\left( J_{3}^{\mathbb{H}}\right) \supset ^{ns}str_{0}\left( J_{3}^{%
\mathbb{H}}\right) \oplus sl(3,\mathbb{R})_{Ehlers}\supset str_{0}(\mathbb{%
R\oplus }J_{2}^{\mathbb{H}})\oplus \mathcal{A}_{4}\oplus sl(3,\mathbb{R}%
)_{Ehlers};  \label{j724a}
\end{equation}%
\begin{eqnarray}
2 &:&qconf\left( J_{3}^{\mathbb{H}}\right) \supset qconf\left( \mathbb{%
R\oplus }J_{2}^{\mathbb{H}}\right) \oplus \mathcal{A}_{4}  \notag
\label{j724b} \\
&\supset &so(5,1)\oplus so(3,3)\oplus sl(2,\mathbb{R})\supset str_{0}\left(
\mathbb{R\oplus }J_{2}^{\mathbb{H}}\right) \oplus sl(3,\mathbb{R}%
)_{Ehlers}\oplus \mathcal{A}_{q=4}.
\end{eqnarray}%
The second step of the Jordan-algebraic interpretation (\ref{j724b}) of the
chain $2$ is intentionally left explicit : intriguingly, it constitutes an
\textit{enhancement} of the symmetry obtained by the $D=5$ case of the Ehlers%
\textit{\ }embedding for the $\left( \mathbb{R\oplus }J_{2}^{\mathbb{H}%
}\right) $-based theory. In fact, this latter reads%
\begin{equation}
qconf\left( J_{3}^{\mathbb{H}}\right) \supset str_{0}\left( \mathbb{R\oplus }%
J_{2}^{\mathbb{H}}\right) \oplus sl(3,\mathbb{R})_{Ehlers}\cong
str_{0}\left( J_{2}^{\mathbb{H}}\right) \oplus sl(3,\mathbb{R}%
)_{Ehlers}\oplus so(1,1),
\end{equation}%
and then the following enhancement takes place:%
\begin{equation}
sl(3,\mathbb{R})_{Ehlers}\oplus so(1,1)\rightarrow sl(4,\mathbb{R})\cong
so(3,3).  \label{enh-1-2}
\end{equation}%
The enhancement (\ref{enh-1-2}) actually hints for another interpretation of
the \textit{non-symmetric} and \textit{next-to-maximal} embedding%
\begin{equation}
qconf\left( J_{3}^{\mathbb{H}_{s}}\right) \supset _{nm}^{ns}so(5,1)\oplus
so(3,3)\oplus sl(2,\mathbb{R}),  \label{enh-2-2}
\end{equation}%
In fact, (\ref{enh-2-2}) can be interpreted as the $D=6$ case of the Ehlers
embedding for the $\mathbb{H}$-based magic supergravity theory (which also
enjoys a twin $\mathcal{N}=6$ fermionic completion), and its next-to-maximal
nature is consistent with the treatment of \cite{super-Ehlers}; in other
words, (\ref{enh-2-2}) can be interpreted as follows:%
\begin{equation}
qconf\left( J_{3}^{\mathbb{H}}\right) \supset _{nm}^{ns}str_{0}\left( J_{2}^{%
\mathbb{H}}\right) \oplus \mathcal{A}_{4}\oplus sl(4,\mathbb{R})_{Ehlers},
\label{enh-3-2}
\end{equation}%
where therefore one of the $so(5,1)$ factor is conceived as $str_{0}\left(
J_{2}^{\mathbb{H}}\right) $, while $so(3,3)$ is nothing but the $D=4$ Ehlers
symmetry, and the commuting $su(2)$ factor is seen as $\mathcal{A}_{4}$. As
a consequence of this interpretation, (\ref{j724b}) can concisely be
rewritten as :%
\begin{eqnarray}
2 &:&qconf\left( J_{3}^{\mathbb{H}}\right) \supset qconf\left( \mathbb{%
R\oplus }J_{2}^{\mathbb{H}}\right) \oplus \mathcal{A}_{q=4}  \notag \\
&\supset &str_{0}\left( J_{2}^{\mathbb{H}}\right) \oplus \mathcal{A}%
_{4}\oplus sl(4,\mathbb{R})_{Ehlers}  \notag \\
&\supset &str_{0}\left( \mathbb{R\oplus }J_{2}^{\mathbb{H}}\right) \oplus
sl(3,\mathbb{R})_{Ehlers}\oplus \mathcal{A}_{4},
\end{eqnarray}%
where%
\begin{equation}
str_{0}\left( \mathbb{R\oplus }J_{2}^{\mathbb{H}}\right) \cong so(1,1)\oplus
str_{0}\left( J_{2}^{\mathbb{H}}\right) .
\end{equation}

\md

\subsection{$E_{7(-25)}$}

This is the minimally non-compact real form of $E_{7}$. Its Jordan
interpretation reads:%
\begin{equation}
E_{7(-25)}\cong conf\left( J_{3}^{\mathbb{O}}\right) \cong der\left( \mathbf{%
F}\left( J_{3}^{\mathbb{O}}\right) \right) .
\end{equation}

\md

\subsubsection{$\cp_{1}^{7(-25)}$\label{731}}

The maximal parabolics ~$\cp_{1}^{7(-25)}$ ~ from \eqref{maxsevc}
corresponds to the Bruhat decomposition:
\begin{equation}
E_{7(-25)}~=~\mathcal{N}_{1}^{7(-25)-}\oplus so(10,2)\oplus so(1,1)\oplus
\mathcal{N}_{1}^{7(-25)+}\newline
\mathbf{,}
\end{equation}%
which can be obtained through the embedding chain
\begin{eqnarray}
E_{7(-25)} &\supset &so(2,10)\oplus sl(2,\mathbb{R})\supset so(2,10)\oplus
so(1,1);  \label{pom-1} \\
\mathbf{133} &=&\left( \mathbf{66,1}\right) \oplus \left( \mathbf{1},\mathbf{%
3}\right) \oplus \left( \mathbf{32}^{\prime }\mathbf{,2}\right) =\mathbf{66}%
_{0}\oplus \mathbf{1}_{0}\oplus \mathbf{1}_{2}\oplus \mathbf{1}_{-2}\oplus
\mathbf{32}_{1}^{\prime }\oplus \mathbf{32}_{-1}^{\prime },
\end{eqnarray}%
which gives rise to a $5$-grading with $\mathcal{N}_{1}^{7(-25)+}=\mathbf{1}%
_{2}\oplus \mathbf{32}_{1}^{\prime }\,$.

A Jordan algebraic interpretation (of the first step) of the chain above is
given by:%
\begin{equation}
conf\left( J_{3}^{\mathbb{O}}\right) \overset{J_{3}^{\mathbb{O}}\supset
\mathbb{R}\oplus J_{2}^{\mathbb{O}}}{\supset }conf\left( \mathbb{R}\oplus
J_{2}^{\mathbb{O}}\right) \oplus \mathcal{A}_{8}.
\end{equation}%
In this case,, the parabolic $so(1,1)$ (giving rise to the $5$-grading) is
the non-compact Cartan of the $sl(2,\mathbb{R})$ factor in the first step of
(\ref{pom-1}), namely of the axio-dilatonic ($S$-duality) factor of $%
conf\left( \mathbb{R}\oplus \mathfrak{J}_{2}^{\mathbb{O}}\right) $.

\md

\subsubsection{$\cp_{2}^{7(-25)}$\label{732}}

The maximal parabolics ~$\cp_{2}^{7(-25)}$ ~ from \eqref{maxsevc}
corresponds to the Bruhat decomposition:
\begin{equation}
E_{7(-25)}~=\mathcal{N}_{2}^{7(-25)-}\oplus ~(so(9,1)\oplus sl(2,\bbr))\
\oplus so(1,1)\oplus \mathcal{N}_{2}^{7(-25)+}\newline
\mathbf{,}
\end{equation}%
which can be obtained through the embedding chain
\begin{eqnarray}
E_{7(-25)} &\supset &so(2,10)\oplus sl(2,\mathbb{R})\supset so(1,9)\oplus
sl(2,\mathbb{R})\oplus so(1,1); \\
\mathbf{133} &=&\left( \mathbf{66,1}\right) \oplus \left( \mathbf{1},\mathbf{%
3}\right) \oplus \left( \mathbf{32}^{\prime }\mathbf{,2}\right) \\
&=&\left( \mathbf{45,1}\right) _{0}\oplus (\mathbf{1},\mathbf{1})_{0}\oplus (%
\mathbf{10},\mathbf{1})_{2}\oplus (\mathbf{10},\mathbf{1})_{-2}\oplus \left(
\mathbf{1},\mathbf{3}\right) _{0}\oplus \left( \mathbf{16,2}\right)
_{-1}\oplus \left( \mathbf{16}^{\prime }\mathbf{,2}\right) _{1}\ .  \notag
\end{eqnarray}%
Thus, it gives rise to a $5$-grading (\textit{cfr.} (\ref{5-1})-(\ref{5-2}),
as well as (\ref{5-1-1})-(\ref{5-2-1})), with $\mathcal{N}_{2}^{7(-25)+}=(%
\mathbf{10},\mathbf{1})_{2}\oplus \left( \mathbf{16}^{\prime }\mathbf{,2}%
\right) _{1}\,$. Since
\begin{eqnarray}
so(2,10)\oplus sl(2,\mathbb{R}) &\cong &conf\left( \mathbb{R}\oplus J_{2}^{%
\mathbb{O}}\right) \cong der\left( \mathbf{F}\left( \mathbb{R}\oplus J_{2}^{%
\mathbb{O}}\right) \right) ; \\
so(1,9) &\cong &str_{0}\left( J_{2}^{\mathbb{O}}\right)  \notag \\
&\cong &str\left( J_{2}^{\mathbb{O}}\right) \ominus so(1,1)\cong
str_{0}\left( \mathbb{R}\oplus \mathbb{J}_{2}^{\mathbb{O}}\right) \ominus
so(1,1),
\end{eqnarray}%
a Jordan algebraic interpretation of the chain above is given by%
\begin{equation}
conf\left( \mathfrak{J}_{3}^{\mathbb{O}}\right) \overset{J_{3}^{\mathbb{O}%
}\supset \mathbb{R}\oplus J_{2}^{\mathbb{O}}}{\supset }conf\left( \mathbb{R}%
\oplus J_{2}^{\mathbb{O}}\right) \oplus \mathcal{A}_{8}\supset str_{0}\left(
\mathbb{R}\oplus J_{2}^{\mathbb{O}}\right) \oplus sl(2,\mathbb{R})\oplus
\mathcal{A}_{8}.  \label{Interpr}
\end{equation}

\md

\subsubsection{$\cp_{3}^{7(-25)}$\label{733}}

The maximal parabolics ~$\cp_{3}^{7(-25)}$ ~ from \eqref{maxsevc}
corresponds to the Bruhat decomposition:
\begin{equation}
E_{7(-25)}~=\mathcal{N}_{3}^{7(-25)-}\oplus ~E_{6(-26)}\oplus so(1,1)\oplus
\mathcal{N}_{3}^{7(-25)}\newline
\mathbf{,}
\end{equation}%
giving rise to a $3$-grading with $\mathcal{N}_{3}^{7(-25)+}=\mathbf{27}%
_{2}^{\prime }\,$. A Jordan algebraic interpretation of the chain above reads%
\begin{equation}
conf\left( J_{3}^{\mathbb{O}}\right) \supset str\left( J_{3}^{\mathbb{O}%
}\right) \cong str_{0}\left( J_{3}^{\mathbb{O}}\right) \oplus so(1,1)_{KK};
\label{j733}
\end{equation}%
therefore, the parabolic $so(1,1)$ (generating the $3$-grading) is the $%
so(1,1)_{KK}$ of the $S^{1}$-reduction $D=5\rightarrow 4$.

\md

\subsection{$E_{8(8)}$}

This is the \textit{split} real form of $E_{8}$. Its Jordan interpretation
reads:%
\begin{equation}
E_{8(8)}\cong qconf\left( J_{3}^{\mathbb{O}_{s}}\right) .
\end{equation}

\md

\subsubsection{$\cp_{1}^{8(8)}$\label{811}}

The maximal parabolics ~$\cp_{1}^{8(8)}$ ~ from \eqref{maxeighta}
corresponds to the Bruhat decomposition:
\begin{equation}
E_{8(8)}~=~\cn_{1}^{\,8(8)-}\oplus so(7,7)\oplus so(1,1)\oplus \cn%
_{1}^{\,8(8)+}\newline
\mathbf{,}
\end{equation}%
which can be obtained through embedding chain:
\begin{eqnarray}
E_{8(8)} &\supset &so(8,8)\supset so(7,7)\oplus so(1,1); \\
\mathbf{248} &=&\mathbf{120}\oplus \mathbf{128}=\mathbf{91}_{0}\oplus
\mathbf{1}_{0}+\mathbf{14}_{2}\oplus \mathbf{14}_{-2}\oplus \mathbf{64}%
_{-1}\oplus \mathbf{64}_{1}^{\prime }.
\end{eqnarray}%
It gives rise to a $5$-grading, with $\cn_{1}^{\,8(8)+}=\mathbf{14}%
_{2}\oplus \mathbf{64}_{1}^{\prime }\,$. ~ By observing that
\begin{equation}
so(8,8)\cong qconf\left( \mathbb{R}\oplus J_{2}^{\mathbb{O}_{s}}\right) ,
\label{so(8,8)}
\end{equation}%
a Jordan algebraic interpretation (of the first step) of the chain above is
given by the embedding $J_{3}^{\mathbb{O}_{s}}\supset \mathbb{R}\oplus
J_{2}^{\mathbb{O}_{s}}$ considered at the $qconf$ level :%
\begin{equation}
qconf\left( J_{3}^{\mathbb{O}_{s}}\right) \supset qconf\left( \mathbb{R}%
\oplus J_{2}^{\mathbb{O}_{s}}\right) \oplus \widetilde{\mathcal{A}}_{8}.
\label{j811b}
\end{equation}

\md

\subsubsection{$\cp_{2}^{8(8)}$\label{812}}

The maximal parabolics ~$\cp_{2}^{8(8)}$ ~ from \eqref{maxeighta}
corresponds to the Bruhat decomposition:
\begin{equation}
E_{8(8)}~=\cn_{2}^{\,8(8)-}\oplus ~sl(7,\bbr)\oplus sl(2,\bbr)\oplus
so(1,1)\oplus \cn_{2}^{\,8(8)+}\newline
\mathbf{,}
\end{equation}%
which can be obtained through the embedding chain
\begin{eqnarray}
E_{8(8)} &\supset &^{ns}sl(9,\mathbb{R})\supset sl(7,\mathbb{R})\oplus sl(2,%
\mathbb{R})\oplus so(1,1); \\
\mathbf{248} &=&\mathbf{80}\oplus \mathbf{84\oplus 84}^{\prime }=\left\{
\begin{array}{l}
\left( \mathbf{48,1}\right) _{0}\oplus \left( \mathbf{1,3}\right) _{0}\oplus
\left( \mathbf{1,1}\right) _{0}\oplus \left( \mathbf{7,2}\right) _{9}\oplus
\left( \mathbf{7}^{\prime }\mathbf{,2}\right) _{-9} \\
\oplus \left( \mathbf{7,1}\right) _{12}\oplus \left( \mathbf{7}^{\prime }%
\mathbf{,1}\right) _{-12} \\
\oplus \left( \mathbf{21,2}\right) _{-3}\oplus \left( \mathbf{21}^{\prime }%
\mathbf{,2}\right) _{3}\oplus \left( \mathbf{35,1}\right) _{6}\oplus \left(
\mathbf{35}^{\prime }\mathbf{,1}\right) _{-6}.%
\end{array}%
\right.
\end{eqnarray}%
Thus, it gives rise to a $9$-grading, with ~$\cn_{2}^{\,8(8)+}=\left(
\mathbf{7,1}\right) _{12}\oplus \left( \mathbf{7,2}\right) _{9}\oplus \left(
\mathbf{35,1}\right) _{6}\oplus \left( \mathbf{21}^{\prime }\mathbf{,2}%
\right) _{3}\,$. ~The interpretation of the first step of the chain above is
provided by the $D=11$ case of the Ehlers embedding, which is actually
relevant for $M$-theory (whose $g_{D=11}=\varnothing $) (\textit{cfr.} \cite%
{super-Ehlers}, and Refs. therein) :
\begin{equation}
qconf\left( \mathfrak{J}_{3}^{\mathbb{O}_{s}}\right) \supset
^{ns}g_{D=11}\oplus \left. sl(D-2,\mathbb{R})\right\vert _{D=11}.  \label{jj}
\end{equation}

\md

\subsubsection{$\cp_{3}^{8(8)}$\label{813}}

The maximal parabolics ~$\cp_{3}^{8(8)}$ ~ from \eqref{maxeighta}
corresponds to the Bruhat decomposition:
\begin{equation}
E_{8(8)}~=~\cn_{3}^{\,8(8)-}\oplus sl(5,\bbr)\oplus sl(3,\bbr)\oplus sl(2,%
\bbr)\oplus so(1,1)\oplus \cn_{3}^{\,8(8)+}\newline
\mathbf{,}
\end{equation}%
which can be obtained by \textit{at least} four chains of embeddings%
\footnote{%
Once again, we use the conventions of \cite{Kephart}.}, respectively denoted
by $1$, $2$, $3$ and $4$ :%
\begin{eqnarray}
1 &:&E_{8(8)}\supset ^{ns}sl(5,\mathbb{R})_{\mathbf{I}}\oplus sl(5,\mathbb{R}%
)_{\mathbf{II}}  \notag \\
&\supset &sl(5,\mathbb{R})\oplus sl(3,\mathbb{R})\oplus sl(2,\mathbb{R}%
)\oplus so(1,1); \\
\mathbf{248} &=&(\mathbf{24},\mathbf{1})\oplus (\mathbf{1},\mathbf{24}%
)\oplus (\mathbf{10},\mathbf{5})\oplus (\mathbf{10}^{\prime },\mathbf{5}%
^{\prime })\oplus \left( \mathbf{5},\mathbf{10}^{\prime }\right) \oplus (%
\mathbf{5}^{\prime },\mathbf{10}) \\
&=&(\mathbf{24},\mathbf{1,1})_{0}\oplus (\mathbf{1},\mathbf{8,1})_{0}\oplus (%
\mathbf{1},\mathbf{1,3})_{0}\oplus (\mathbf{1},\mathbf{1,1})_{0}\oplus (%
\mathbf{1},\mathbf{3,2})_{5}\oplus (\mathbf{1},\mathbf{3}^{\prime }\mathbf{,2%
})_{-5}  \notag \\
&&\oplus (\mathbf{10},\mathbf{1,2})_{-3}\oplus (\mathbf{10},\mathbf{3,1}%
)_{2}\oplus (\mathbf{10}^{\prime },\mathbf{1,2})_{3}\oplus (\mathbf{10}%
^{\prime },\mathbf{3}^{\prime }\mathbf{,1})_{-2}  \notag \\
&&\oplus (\mathbf{5}^{\prime },\mathbf{1,1})_{-6}\oplus (\mathbf{5}^{\prime
},\mathbf{3}^{\prime }\mathbf{,1})_{4}\oplus (\mathbf{5}^{\prime },\mathbf{%
3,2})_{-1}\oplus (\mathbf{5},\mathbf{1,1})_{6}\oplus (\mathbf{5},\mathbf{3,1}%
)_{-4}\oplus (\mathbf{5},\mathbf{3}^{\prime }\mathbf{,2})_{1};  \notag
\end{eqnarray}%
\begin{eqnarray}
2 &:&E_{8(8)}\supset E_{7(7)}\oplus sl(2,\mathbb{R})  \notag \\
&\supset &sl(8,\mathbb{R})\oplus sl(2,\mathbb{R})\supset sl(5,\mathbb{R}%
)\oplus sl(3,\mathbb{R})\oplus sl(2,\mathbb{R})\oplus so(1,1); \\
\mathbf{248} &=&(\mathbf{133},\mathbf{1})\oplus (\mathbf{1},\mathbf{3}%
)\oplus (\mathbf{56},\mathbf{2})=(\mathbf{63},\mathbf{1})\oplus (\mathbf{70},%
\mathbf{1})\oplus (\mathbf{1},\mathbf{3})\oplus (\mathbf{28},\mathbf{2}%
)\oplus (\mathbf{28}^{\prime },\mathbf{2}) \\
&=&(\mathbf{24},\mathbf{1,1})_{0}\oplus (\mathbf{1},\mathbf{8,1})_{0}\oplus (%
\mathbf{1},\mathbf{1},\mathbf{1})_{0}\oplus (\mathbf{5},\mathbf{3}^{\prime },%
\mathbf{1})_{8}\oplus (\mathbf{5}^{\prime },\mathbf{3},\mathbf{1})_{-8}
\notag \\
&&\oplus (\mathbf{5,1},\mathbf{1})_{-12}\oplus (\mathbf{5}^{\prime }\mathbf{%
,1},\mathbf{1})_{12}\oplus (\mathbf{10,3}^{\prime },\mathbf{1})_{-4}\oplus (%
\mathbf{10}^{\prime }\mathbf{,3},\mathbf{1})_{4}\oplus (\mathbf{1},\mathbf{1}%
,\mathbf{3})_{0}  \notag \\
&&\oplus (\mathbf{1,3}^{\prime },\mathbf{2})_{-10}\oplus (\mathbf{5,3},%
\mathbf{2})_{-2}\oplus (\mathbf{10,1},\mathbf{2})_{6}\oplus (\mathbf{1,3},%
\mathbf{2})_{10}\oplus (\mathbf{5}^{\prime }\mathbf{,3}^{\prime },\mathbf{2}%
)_{2}\oplus (\mathbf{10}^{\prime }\mathbf{,1},\mathbf{2})_{-6};  \notag
\end{eqnarray}%
\begin{eqnarray}
3 &:&E_{8(8)}\supset E_{7(7)}\oplus sl(2,\mathbb{R})  \notag \\
&\supset &^{ns}sl(6,\mathbb{R})\oplus sl(3,\mathbb{R})\oplus sl(2,\mathbb{R}%
)\supset sl(5,\mathbb{R})\oplus sl(3,\mathbb{R})\oplus sl(2,\mathbb{R}%
)\oplus so(1,1); \\
\mathbf{248} &=&(\mathbf{133},\mathbf{1})\oplus (\mathbf{1},\mathbf{3}%
)\oplus (\mathbf{56},\mathbf{2}) \\
&=&\left( \mathbf{35,1,1}\right) \oplus \left( \mathbf{1},\mathbf{8,1}%
\right) \oplus \left( \mathbf{15,3}^{\prime },\mathbf{1}\right) \oplus
\left( \mathbf{15}^{\prime }\mathbf{,3,1}\right) \oplus (\mathbf{1},\mathbf{1%
},\mathbf{3})\oplus (\mathbf{6},\mathbf{3},\mathbf{2})\oplus (\mathbf{6}%
^{\prime },\mathbf{3}^{\prime },\mathbf{2})\oplus (\mathbf{20},\mathbf{1},%
\mathbf{2})  \notag \\
&=&\left( \mathbf{24,1,1}\right) _{0}\oplus \left( \mathbf{1,1,1}\right)
_{0}\oplus \left( \mathbf{5,1,1}\right) _{6}\oplus \left( \mathbf{5}^{\prime
}\mathbf{,1,1}\right) _{-6}\oplus \left( \mathbf{1},\mathbf{8,1}\right) _{0}
\notag \\
&&\oplus \left( \mathbf{5,3}^{\prime },\mathbf{1}\right) _{-4}\oplus \left(
\mathbf{10,3}^{\prime },\mathbf{1}\right) _{2}\oplus \left( \mathbf{5}%
^{\prime }\mathbf{,3},\mathbf{1}\right) _{4}\oplus \left( \mathbf{10}%
^{\prime }\mathbf{,3},\mathbf{1}\right) _{-2}  \notag \\
&&\oplus (\mathbf{1},\mathbf{1},\mathbf{3})_{0}\oplus (\mathbf{1},\mathbf{3},%
\mathbf{2})_{-5}\oplus (\mathbf{5},\mathbf{3},\mathbf{2})_{1}\oplus (\mathbf{%
1},\mathbf{3}^{\prime },\mathbf{2})_{5}\oplus (\mathbf{5}^{\prime },\mathbf{3%
}^{\prime },\mathbf{2})_{-1}\oplus (\mathbf{10},\mathbf{1},\mathbf{2}%
)_{-3}\oplus (\mathbf{10}^{\prime },\mathbf{1},\mathbf{2})_{3}\ ;  \notag
\end{eqnarray}%
\begin{eqnarray}
4 &:&E_{8(8)}\supset ^{ns}E_{6(6)}\oplus sl(3,\mathbb{R})  \notag \\
&\supset &sl(6,\mathbb{R})\oplus sl(3,\mathbb{R})\oplus sl(2,\mathbb{R}%
)\supset sl(5,\mathbb{R})\oplus sl(3,\mathbb{R})\oplus sl(2,\mathbb{R}%
)\oplus so(1,1); \\
\mathbf{248} &=&(\mathbf{78},\mathbf{1})\oplus (\mathbf{1},\mathbf{8})\oplus
(\mathbf{27},\mathbf{3})\oplus (\mathbf{27}^{\prime },\mathbf{3}^{\prime })
\\
&=&(\mathbf{35},\mathbf{1,1})\oplus (\mathbf{1},\mathbf{1,3})\oplus (\mathbf{%
20},\mathbf{1,2})\oplus (\mathbf{1},\mathbf{8},\mathbf{1})  \notag \\
&&\oplus (\mathbf{6},\mathbf{3,2})\oplus (\mathbf{15}^{\prime },\mathbf{3,1}%
)\oplus (\mathbf{6}^{\prime },\mathbf{3}^{\prime }\mathbf{,2})\oplus (%
\mathbf{15},\mathbf{3}^{\prime }\mathbf{,1})  \notag \\
&=&(\mathbf{24},\mathbf{1,1})_{0}\oplus (\mathbf{1},\mathbf{1,1})_{0}\oplus (%
\mathbf{5},\mathbf{1,1})_{6}\oplus (\mathbf{5}^{\prime },\mathbf{1,1}%
)_{-6}\oplus (\mathbf{1},\mathbf{1,3})_{0}\oplus (\mathbf{10},\mathbf{1,2}%
)_{-3}\oplus (\mathbf{10}^{\prime },\mathbf{1,2})_{3}  \notag \\
&&\oplus (\mathbf{1},\mathbf{8},\mathbf{1})_{0}\oplus (\mathbf{1},\mathbf{3,2%
})_{-5}\oplus (\mathbf{5},\mathbf{3,2})_{1}\oplus (\mathbf{5},\mathbf{3}%
^{\prime }\mathbf{,1})_{-4}\oplus (\mathbf{10},\mathbf{3}^{\prime }\mathbf{,1%
})_{2}  \notag \\
&&\oplus (\mathbf{5}^{\prime },\mathbf{3,1})_{4}\oplus (\mathbf{10}^{\prime
},\mathbf{3,1})_{-2}\oplus (\mathbf{1},\mathbf{3}^{\prime }\mathbf{,2}%
)_{5}\oplus (\mathbf{5}^{\prime },\mathbf{3}^{\prime }\mathbf{,2})_{-1}\ .
\notag
\end{eqnarray}%
The chains of embeddings $1$, $2$, $3$ and $4$ give rise to a $13$-grading,
with
\begin{eqnarray}
\cn_{3\,;\,1}^{\,8(8)+} &=&(\mathbf{5},\mathbf{1,1})_{6}\oplus (\mathbf{1},%
\mathbf{3,2})_{5}\oplus (\mathbf{5}^{\prime },\mathbf{3}^{\prime }\mathbf{,1}%
)_{4}\oplus (\mathbf{10}^{\prime },\mathbf{1,2})_{3}\oplus (\mathbf{10},%
\mathbf{3,1})_{2}\oplus (\mathbf{5},\mathbf{3}^{\prime }\mathbf{,2})_{1}; \\
\cn_{3\,;\,2}^{\,8(8)+} &=&(\mathbf{5}^{\prime }\mathbf{,1},\mathbf{1}%
)_{12}\oplus (\mathbf{1,3},\mathbf{2})_{10}\oplus (\mathbf{5},\mathbf{3}%
^{\prime },\mathbf{1})_{8}\oplus (\mathbf{10,1},\mathbf{2})_{6}\oplus (%
\mathbf{10}^{\prime }\mathbf{,3},\mathbf{1})_{4}\oplus (\mathbf{5}^{\prime }%
\mathbf{,3}^{\prime },\mathbf{2})_{2}; \\
\cn_{3\,;\,3}^{\,8(8)+} &\cong &\cn_{3\,;\,4}^{\,8(8)\oplus }=(\mathbf{5},%
\mathbf{1,1})_{6}\oplus (\mathbf{1},\mathbf{3}^{\prime }\mathbf{,2}%
)_{5}\oplus (\mathbf{5}^{\prime },\mathbf{3,1})_{4}\oplus (\mathbf{10}%
^{\prime },\mathbf{1,2})_{3}\oplus (\mathbf{10},\mathbf{3}^{\prime }\mathbf{%
,1})_{2}\oplus (\mathbf{5},\mathbf{3,2})_{1},  \notag \\
&&
\end{eqnarray}%
all of real dimension $106$. Note that $1\leftrightarrow 3\equiv d)$ \textit{%
iff} $\mathbf{3}\leftrightarrow \mathbf{3}^{\prime }$ in $sl(3,\mathbb{R})$,
and $2\leftrightarrow 3\equiv 4$ \textit{iff} the weights of the parabolic $%
so(1,1)$ of $2$ are multiplied by $1/2$, and \textit{iff} $\mathbf{5}%
\leftrightarrow \mathbf{5}^{\prime }$ and $\mathbf{3}\leftrightarrow \mathbf{%
3}^{\prime }$ in $sl(5,\mathbb{R})$ and in $sl(3,\mathbb{R})$, respectively.

The Jordan algebraic interpretation of the chains $1$, $2$, $3$ and $4$ goes
as follows. The first (\textit{non-symmetric}) embedding in chain $1$ can be
interpreted as the $D=7$ case of the Ehlers embedding for the $\mathbb{O}_{s}
$-based theory (maximal supergravity) \cite{super-Ehlers}, thus determining
the $D=7$ Ehlers $sl(5,\mathbb{R})_{Ehlers}$; then, the subsequent embedding
enjoys a(n \textit{at least}) twofold interpretation : it can be conceived
as the uplift to $D=8$, where the $U$-duality Lie algebra is $g_{D=8}\left(
\mathbb{O}_{s}\right) =sl(3,\mathbb{R})\oplus sl(2,\mathbb{R})$ (as such,
the parabolic $so(1,1)$ in this chain is nothing but the KK $so(1,1)_{KK}$
in the $S^{1}$-reduction $D=8\rightarrow 7$ of maximal supergravity), or it
can be seen as a further decomposition of $sl(5,\mathbb{R})_{Ehlers}$ into $%
sl(3,\mathbb{R})\oplus sl(2,\mathbb{R})\oplus so(1,1)$ :%
\begin{equation}
1:qconf\left( J_{3}^{\mathbb{O}_{s}}\right) \supset ^{ns}g_{D=7}\left(
\mathbb{O}_{s}\right) \oplus sl(5,\mathbb{R})_{Ehlers}\supset \left\{
\begin{array}{l}
g_{D=8}\left( \mathbb{O}_{s}\right) \oplus sl(5,\mathbb{R})_{Ehlers}\oplus
so(1,1)_{KK}; \\
\text{\textit{or}} \\
g_{D=7}\left( \mathbb{O}_{s}\right) \oplus sl(3,\mathbb{R})\oplus sl(2,%
\mathbb{R})\oplus so(1,1).%
\end{array}%
\right.   \label{interpr-a}
\end{equation}%
On the other hand, the first embedding in chain $2$ enjoys a(n \textit{at
least}) twofold interpretation, by virtue of the twofold characterization of
$E_{7(7)}$ given by (\ref{E7(7)-1})-(\ref{E7(7)-2}) (ultimately due to the
symmetry of the double-split Magic Square $\mathcal{L}_{3}\left( \mathbb{A}%
_{s},\mathbb{B}_{s}\right) $ \cite{BS}) : it can be interpreted as the
Ehlers embedding for the $\mathbb{O}_{s}$-based theory (corresponding to the
uplift $D=3\rightarrow 4$ of maximal supergravity, and giving rise to the $%
D=4$ Ehlers symmetry $sl(2,\mathbb{R})_{Ehlers}$), or it can be conceived as
a consequence of the embedding $\mathbb{O}_{s}\supset \mathbb{H}_{s}$,
evaluated at the $qconf$ level. The second step also has an \textit{at least}
twofold interpretation, given by (\ref{I-one})-(\ref{II-two}) :%
\begin{eqnarray}
2 &:&qconf\left( J_{3}^{\mathbb{O}_{s}}\right) \supset \left\{
\begin{array}{l}
conf\left( J_{3}^{\mathbb{O}_{s}}\right) \oplus sl(2,\mathbb{R})_{Ehlers} \\
\text{\textit{or}} \\
qconf\left( J_{3}^{\mathbb{H}_{s}}\right) \oplus \widetilde{\mathcal{A}}_{4}%
\end{array}%
\right. \supset \left\{
\begin{array}{l}
\left. sl(q,\mathbb{R})\right\vert _{q=8}\oplus sl(2,\mathbb{R})_{Ehlers} \\
\text{\textit{or}} \\
\left. sl(q+4,\mathbb{R})\right\vert _{q=4}\oplus \widetilde{\mathcal{A}}_{4}%
\end{array}%
\right.   \notag \\
&\supset &\left\{
\begin{array}{l}
g_{D=8}\left( \mathbb{O}_{s}\right) \oplus sl(5,\mathbb{R})_{Ehlers}\oplus
so(1,1)_{KK}; \\
\text{\textit{or}} \\
g_{D=7}\left( \mathbb{O}_{s}\right) \oplus sl(3,\mathbb{R})\oplus sl(2,%
\mathbb{R})\oplus so(1,1).%
\end{array}%
\right.   \label{interpr-b}
\end{eqnarray}%
The first embedding in chain $3$ is the same as the first embedding of chain
$2$. The second step consists of a maximal, \textit{non-symmetric}
embedding, and it has a(n \textit{at least}) twofold interpretation, by
virtue of the twofold characterization of $sl(6,\mathbb{R})$ given by (\ref%
{gg}) (once again due to the symmetry of $\mathcal{L}_{3}\left( \mathbb{A}%
_{s},\mathbb{B}_{s}\right) $ \cite{BS}) : it can be interpreted as the JP%
\textit{\ }embedding for the $\mathbb{H}_{s}$-based theory (thus introducing
the $D=5$ Ehlers $sl(3,\mathbb{R})_{Ehlers}$), or as the consequence of the
(non-maximal) embedding $\mathbb{O}_{s}\supset \mathbb{C}_{s}$, evaluated at
the level of conformal symmetries :%
\begin{eqnarray}
3 &:&qconf\left( J_{3}^{\mathbb{O}_{s}}\right) \supset \left\{
\begin{array}{l}
conf\left( J_{3}^{\mathbb{O}_{s}}\right) \oplus sl(2,\mathbb{R})_{Ehlers} \\
\text{\textit{or}} \\
QConf\left( \mathfrak{J}_{3}^{\mathbb{H}_{s}}\right) \times \widetilde{%
\mathcal{A}}_{q=4}%
\end{array}%
\right.   \notag \\
&\supset &^{ns}\left\{
\begin{array}{l}
conf\left( J_{3}^{\mathbb{C}_{s}}\right) \oplus sl(3,\mathbb{R})\oplus sl(2,%
\mathbb{R})_{Ehlers} \\
\text{\textit{or}} \\
str_{0}\left( J_{3}^{\mathbb{H}_{s}}\right) \oplus sl(3,\mathbb{R}%
)_{Ehlers}\oplus \widetilde{\mathcal{A}}_{4}%
\end{array}%
\right.   \notag \\
&\supset &\left\{
\begin{array}{l}
g_{D=8}\left( \mathbb{O}_{s}\right) \oplus sl(5,\mathbb{R})_{Ehlers}\oplus
so(1,1)_{KK}; \\
\text{\textit{or}} \\
g_{D=7}\left( \mathbb{O}_{s}\right) \oplus sl(3,\mathbb{R})\oplus sl(2,%
\mathbb{R})\oplus so(1,1).%
\end{array}%
\right.   \label{interpr-c}
\end{eqnarray}%
Finally, each of the steps of the chain $4$ has an \textit{at least} twofold
interpretation. As a consequence of the twofold characterization of $E_{6(6)}
$ given by (\ref{uno-3})-(\ref{uno-4}) (once again, due to the symmetry of $%
\mathcal{L}_{3}\left( \mathbb{A}_{s},\mathbb{B}_{s}\right) $ \cite{BS}), a
first line of interpretation (the upper one in (\ref{interpr-c})) conceives
the first embedding as the JP embedding for the $\mathbb{O}_{s}$-based
theory (corresponding to the uplift $D=3\rightarrow 5$ of maximal
supergravity, and thus giving rise to the $D=5$ Ehlers symmetry $sl(3,%
\mathbb{R})_{Ehlers}$), followed by the consequence of the embedding $%
\mathbb{O}_{s}\supset \mathbb{H}_{s}$ at the level of reduced structure
symmetry. A second line of interpretation (the lower one in (\ref{interpr-c}%
)) sees the first embedding as a consequence of the non-maximal embedding $%
\mathbb{O}_{s}\supset \mathbb{C}_{s}$, evaluated at the $qconf$ level, and
then followed an Ehlers embedding, corresponding to an uplift $%
D=3\rightarrow 4$, and therefore introducing the $D=4$ $sl(2,\mathbb{R}%
)_{Ehlers}$ group. Note that the first two steps of chain $4$ realize the $%
D=8$ case of the Ehlers embedding for the $\mathbb{O}_{s}$-based theory \cite%
{super-Ehlers}; thus, in this view the $sl(6,\mathbb{R})$ occurring in the
second step is nothing but the $D=8$ Ehlers group :%
\begin{eqnarray}
4 &:&\left\{
\begin{array}{l}
qconf\left( J_{3}^{\mathbb{O}_{s}}\right) \supset ^{ns}str_{0}\left( J_{3}^{%
\mathbb{O}_{s}}\right) \oplus sl(3,\mathbb{R})_{Ehlers}\supset str_{0}\left(
J_{3}^{\mathbb{H}_{s}}\right) \oplus sl(3,\mathbb{R})_{Ehlers}\oplus
\widetilde{\mathcal{A}}_{4} \\
\text{\textit{or}} \\
qconf\left( J_{3}^{\mathbb{O}_{s}}\right) \supset ^{ns}qconf\left( J_{3}^{%
\mathbb{C}_{s}}\right) \oplus sl(3,\mathbb{R})\supset conf\left( J_{3}^{%
\mathbb{C}_{s}}\right) \oplus sl(3,\mathbb{R})\oplus sl(2,\mathbb{R}%
)_{Ehlers}%
\end{array}%
\right.   \notag \\
&\supset &\left\{
\begin{array}{l}
g_{D=8}\left( \mathbb{O}_{s}\right) \oplus sl(5,\mathbb{R})_{Ehlers}\oplus
so(1,1)_{KK}; \\
\text{\textit{or}} \\
g_{D=7}\left( \mathbb{O}_{s}\right) \oplus sl(3,\mathbb{R})\oplus sl(2,%
\mathbb{R})\oplus so(1,1).%
\end{array}%
\right.   \label{interpr-d}
\end{eqnarray}

\md

\subsubsection{$\cp_{4}^{8(8)}$\label{814}}

The maximal parabolics ~$\cp_{4}^{8(8)}$ ~ from \eqref{maxeighta}
corresponds to the Bruhat decomposition:
\begin{equation}
E_{8(8)}~=~\cn_{4}^{\,8(8)-}\oplus sl(5,\bbr)\oplus sl(4,\bbr)\oplus
so(1,1)\oplus \cn_{4}^{\,8(8)+}\newline
\mathbf{,}
\end{equation}%
which can be obtained by \textit{at least} two chains of embeddings,
respectively denoted by $1$ and $2$ :%
\begin{eqnarray}
1 &:&E_{8(8)}\supset ^{ns}sl(5,\mathbb{R})_{\mathbf{I}}\oplus sl(5,\mathbb{R}%
)_{\mathbf{II}}\supset sl(5,\mathbb{R})\oplus sl(4,\mathbb{R})\oplus so(1,1);
\\
\mathbf{248} &=&(\mathbf{24},\mathbf{1})\oplus (\mathbf{1},\mathbf{24}%
)\oplus (\mathbf{10},\mathbf{5})\oplus (\mathbf{10}^{\prime },\mathbf{5}%
^{\prime })\oplus \left( \mathbf{5},\mathbf{10}^{\prime }\right) \oplus (%
\mathbf{5}^{\prime },\mathbf{10})  \notag \\
&=&(\mathbf{24},\mathbf{1})_{0}\oplus (\mathbf{1},\mathbf{15})_{0}\oplus (%
\mathbf{1},\mathbf{1})_{0}\oplus (\mathbf{1},\mathbf{4})_{5}\oplus (\mathbf{1%
},\mathbf{4}^{\prime })_{-5}  \notag \\
&&\oplus (\mathbf{10},\mathbf{1})_{-4}\oplus (\mathbf{10},\mathbf{4}%
)_{1}\oplus (\mathbf{10}^{\prime },\mathbf{1})_{4}\oplus (\mathbf{10}%
^{\prime },\mathbf{4}^{\prime })_{-1}  \notag \\
&&\oplus (\mathbf{5}^{\prime },\mathbf{4})_{-3}\oplus (\mathbf{5}^{\prime },%
\mathbf{6})_{2}\oplus (\mathbf{5},\mathbf{4}^{\prime })_{3}\oplus (\mathbf{5}%
,\mathbf{6})_{-2};
\end{eqnarray}%
\begin{eqnarray}
2 &:&E_{8(8)}\supset ^{ns}sl(9,\mathbb{R})\supset sl(5,\mathbb{R})\oplus
sl(4,\mathbb{R})\oplus so(1,1); \\
\mathbf{248} &=&\mathbf{80}\oplus \mathbf{84\oplus 84}^{\prime }=(\mathbf{24}%
,\mathbf{1})_{0}\oplus (\mathbf{1},\mathbf{15})_{0}\oplus (\mathbf{1},%
\mathbf{1})_{0}\oplus (\mathbf{5},\mathbf{4}^{\prime })_{9}\oplus (\mathbf{5}%
^{\prime },\mathbf{4})_{-9}  \notag \\
&&\oplus (\mathbf{1},\mathbf{4}^{\prime })_{-15}\oplus (\mathbf{5},\mathbf{6}%
)_{-6}\oplus (\mathbf{10}^{\prime },\mathbf{1})_{12}\oplus (\mathbf{10},%
\mathbf{4})_{3}  \notag \\
&&\oplus (\mathbf{1},\mathbf{4})_{15}\oplus (\mathbf{5}^{\prime },\mathbf{6}%
)_{6}\oplus (\mathbf{10},\mathbf{1})_{-12}\oplus (\mathbf{10}^{\prime },%
\mathbf{4}^{\prime })_{-3}.
\end{eqnarray}%
The chains of embeddings $1$ and $2$ give rise to an $11$-grading, with $\cn%
_{4\,;\,a)}^{\,8(8)+}=(\mathbf{1},\mathbf{4})_{5}\oplus (\mathbf{10}^{\prime
},\mathbf{1})_{4}\oplus (\mathbf{5},\mathbf{4}^{\prime })_{3}\oplus (\mathbf{%
5}^{\prime },\mathbf{6})_{2}\oplus (\mathbf{10},\mathbf{4})_{1}\,$, ~$\cn%
_{4\,;\,b)}^{\,8(8)+}=(\mathbf{1},\mathbf{4})_{15}\oplus (\mathbf{10}%
^{\prime },\mathbf{1})_{12}\oplus (\mathbf{5},\mathbf{4}^{\prime
})_{9}\oplus (\mathbf{5}^{\prime },\mathbf{6})_{6}\oplus (\mathbf{10},%
\mathbf{4})_{3}$, both of real dimension $104$. Note that $1\leftrightarrow 2
$ \textit{iff} the weights of the parabolic $so(1,1)$ of $1$ are multiplied
by $3$. The first steps of chains $1$ and $2$ are the same as the first
steps of chains $1$ of subsubsection \ref{813} and of subsubsection \ref{812}
above, respectively, and thus they correspondingly enjoy the same Jordan
algebraic interpretation (in particular, (\ref{jj}) for the latter).

\md

\subsubsection{$\cp_{5}^{8(8)}$\label{815}}

The maximal parabolics ~$\cp_{5}^{8(8)}$ ~ from \eqref{maxeighta}
corresponds to the Bruhat decomposition:
\begin{equation}
E_{8(8)}~=~\cn_{5}^{\,8(8)-}\oplus so(5,5)\oplus sl(3,\bbr)\oplus
so(1,1)\oplus \cn_{5}^{\,8(8)+}\newline
\mathbf{,}
\end{equation}%
which can be obtained through \textit{at least} two embedding chains,
respectively denoted by $1$ and $2$ (recall that $so(3,3)\cong sl(4,\mathbb{R%
})$):%
\begin{eqnarray}
1 &:&E_{8(8)}\supset so(8,8)\supset so(5,5)\oplus so(3,3)\supset
so(5,5)\oplus sl(3,\mathbb{R})\oplus so(1,1); \\
\mathbf{248} &=&\mathbf{120}\oplus \mathbf{128}=(\mathbf{45},\mathbf{1}%
)\oplus (\mathbf{1},\mathbf{15})\oplus (\mathbf{10},\mathbf{6})\oplus (%
\mathbf{16},\mathbf{4})\oplus (\mathbf{16}^{\prime },\mathbf{4}^{\prime })
\notag \\
&=&(\mathbf{45},\mathbf{1})_{0}\oplus (\mathbf{1},\mathbf{1})_{0}\oplus (%
\mathbf{1},\mathbf{3})_{4}\oplus (\mathbf{1},\mathbf{3}^{\prime
})_{-4}\oplus (\mathbf{1},\mathbf{8})_{0}  \notag \\
&&\oplus (\mathbf{10},\mathbf{3})_{-2}\oplus (\mathbf{10},\mathbf{3}^{\prime
})_{2}\oplus (\mathbf{16},\mathbf{1})_{-3}\oplus \left( \mathbf{16},\mathbf{3%
}\right) _{1}\oplus (\mathbf{16}^{\prime },\mathbf{1})_{3}\oplus (\mathbf{16}%
^{\prime },\mathbf{3}^{\prime })_{-1};
\end{eqnarray}%
\begin{eqnarray}
2 &:&E_{8(8)}\supset ^{ns}E_{6(6)}\oplus sl(3,\mathbb{R})\supset
so(5,5)\oplus sl(3,\mathbb{R})\oplus so(1,1); \\
\mathbf{248} &=&(\mathbf{78},\mathbf{1})\oplus (\mathbf{1},\mathbf{8})\oplus
(\mathbf{27},\mathbf{3})\oplus (\mathbf{27}^{\prime },\mathbf{3}^{\prime })
\notag \\
&=&(\mathbf{45},\mathbf{1})_{0}\oplus (\mathbf{1},\mathbf{1})_{0}\oplus (%
\mathbf{16},\mathbf{1})_{3}\oplus (\mathbf{16}^{\prime },\mathbf{1}%
)_{-3}\oplus (\mathbf{1},\mathbf{1},\mathbf{8})_{0}  \notag \\
&&\oplus (\mathbf{1},\mathbf{3})_{-4}\oplus (\mathbf{10},\mathbf{3}%
)_{2}\oplus (\mathbf{16},\mathbf{3})_{-1}\oplus (\mathbf{1},\mathbf{3}%
^{\prime })_{4}\oplus (\mathbf{10},\mathbf{3}^{\prime })_{-2}\oplus (\mathbf{%
16}^{\prime },\mathbf{3}^{\prime })_{1}.
\end{eqnarray}%
The chains of embeddings $1$ and $2$ give rise to a $9$-grading, with $\cn%
_{3\,;\,a)}^{\,8(8)+}=(\mathbf{1},\mathbf{3})_{4}\oplus (\mathbf{16}^{\prime
},\mathbf{1})_{3}\oplus (\mathbf{10},\mathbf{3}^{\prime })_{2}\oplus (%
\mathbf{16},\mathbf{3})_{1}$, ~ $\cn_{3\,;\,b)}^{\,8(8)+}=(\mathbf{1},%
\mathbf{3}^{\prime })_{4}\oplus (\mathbf{16},\mathbf{1})_{3}\oplus (\mathbf{%
10},\mathbf{3})_{2}\oplus (\mathbf{16}^{\prime },\mathbf{3}^{\prime })_{1}$,
both of real dimension $97$. Note that $1\leftrightarrow 2$ \textit{iff} the
weights of the parabolic $so(1,1)$'s gets flipped (or, equivalently, \textit{%
iff} all $\ \cm_{5}^{8(8)}$-irreps. get conjugated). The Jordan-algebraic
interpretation of chains $1$ and $2$ is based on the characterization (\ref%
{so(8,8)}) and on%
\begin{eqnarray}
so(5,5) &\cong &str_{0}(J_{2}^{\mathbb{O}_{s}})\cong str_{0}(\mathbb{R}%
\oplus J_{2}^{\mathbb{O}})\ominus so(1,1) \\
&\cong &qconf\left( \mathbb{R}\oplus J_{2}^{\mathbb{C}_{s}}\right) .
\label{qc}
\end{eqnarray}%
The first step of chain $1$ is a consequence of the embedding $J_{3}^{%
\mathbb{O}_{s}}\supset \mathbb{R}\oplus J_{2}^{\mathbb{O}_{s}}$, considered
at the $qconf$ level; then, by virtue of the quasi-conformal interpretation (%
\ref{qc}) of $so(5,5)$, the second step is a consequence (for cubic
semi-simple Jordan algebras) of the non-maximal embedding $\mathbb{O}%
_{s}\supset \mathbb{C}_{s}$ (note the commuting factor $sl(4,\mathbb{R})$).
Note also that the \textit{next-to-maximal} (\textit{non-symmetric})
embedding%
\begin{equation}
E_{8(8)}\supset so(5,5)\oplus so(3,3)\cong so(5,5)\oplus sl(4,\mathbb{R})
\end{equation}%
is the $D=6$ case of the Ehlers embedding \cite{super-Ehlers} for the $%
\mathbb{O}_{s}$-based theory (maximal supergravity), thus characterizing $%
sl(4,\mathbb{R})$ as the $D=6$ Ehlers symmetry (\textit{cfr.} (\ref{aa})
below). Concerning the chain $2$, it has ( \textit{at least}) twofold
interpretation : in the first interpretation, its first step is the JP
embedding for the $\mathbb{O}_{s}$-based theory (maximal supergravity), thus
giving rise to the $D=5$ Ehlers $sl(3,\mathbb{R})_{Ehlers}$; then, in the
second step a further uplift to $D=6$ is performed, introducing the KK $%
so(1,1)_{KK}$ of the $S^{1}$-reduction $D=6\rightarrow 5$, which is the
parabolic $so(1,1)$ in this chain. In the second interpretation, the first
step is a consequence (for cubic simple Jordan algebras) of the non-maximal
embedding $\mathbb{O}_{s}\supset \mathbb{C}_{s}$ (note the commuting factor $%
sl(3,\mathbb{R})$), while the second step is a consequence of the embedding $%
J_{3}^{\mathbb{C}_{s}}\supset \mathbb{R}\oplus J_{2}^{\mathbb{C}_{s}}$,
considered for quasi-conformal symmetries.
\begin{eqnarray}
1 &:&qconf\left( J_{3}^{\mathbb{O}_{s}}\right) \supset qconf\left( \mathbb{R}%
\oplus J_{2}^{\mathbb{O}_{s}}\right) \oplus \widetilde{\mathcal{A}}%
_{8}\supset qconf\left( \mathbb{R}\oplus J_{2}^{\mathbb{C}_{s}}\right)
\oplus \widetilde{\mathcal{A}}_{8}\oplus sl(4,\mathbb{R})_{Ehlers}  \notag \\
&\supset &qconf\left( \mathbb{R}\oplus J_{2}^{\mathbb{C}_{s}}\right) \oplus
\widetilde{\mathcal{A}}_{8}\oplus sl(3,\mathbb{R})\oplus so(1,1);  \label{aa}
\end{eqnarray}%
\begin{equation}
2:\left\{
\begin{array}{l}
qconf(J_{3}^{\mathbb{O}_{s}})\supset ^{ns}str_{0}(J_{3}^{\mathbb{O}%
_{s}})\oplus sl(3,\mathbb{R})_{Ehlers}\supset str_{0}(J_{2}^{\mathbb{O}%
_{s}})\oplus \widetilde{\mathcal{A}}_{8}\oplus sl(3,\mathbb{R}%
)_{Ehlers}\oplus so(1,1)_{KK}\ . \\
\text{\textit{or}} \\
qconf\left( J_{3}^{\mathbb{O}_{s}}\right) \supset ^{ns}qconf\left( J_{3}^{%
\mathbb{C}_{s}}\right) \oplus sl(3,\mathbb{R})\supset qconf\left( \mathbb{R}%
\oplus J_{2}^{\mathbb{C}_{s}}\right) \oplus sl(3,\mathbb{R})\oplus
\widetilde{\mathcal{A}}_{2}\ .%
\end{array}%
\right.   \label{bb}
\end{equation}

\md

\subsubsection{$\cp_{6}^{8(8)}$\label{816}}

The maximal parabolics ~$\cp_{6}^{8(8)}$ ~ from \eqref{maxeighta}
corresponds to the Bruhat decomposition:
\begin{equation}
E_{8(8)}~=\cn_{6}^{\,8(8)-}\oplus E_{6(6)}\oplus sl(2,\bbr)\oplus
so(1,1)\oplus \cn_{6}^{\,8(8)+}\newline
\mathbf{,}
\end{equation}%
which can be obtained through \textit{at least} two chains of embeddings,
respectively denoted by $1$ and $2$ :%
\begin{eqnarray}
1 &:&E_{8(8)}\supset E_{7(7)}\oplus sl(2,\mathbb{R})\supset E_{6(6)}\oplus
sl(2,\mathbb{R})\oplus so(1,1); \\
\mathbf{248} &=&(\mathbf{133},\mathbf{1})\oplus (\mathbf{1},\mathbf{3}%
)\oplus (\mathbf{56},\mathbf{2})  \notag \\
&=&\left( \mathbf{78,1}\right) _{0}\oplus \left( \mathbf{1,1}\right)
_{0}\oplus \left( \mathbf{27,1}\right) _{-2}\oplus \left( \mathbf{27}%
^{\prime },\mathbf{1}\right) _{2}\oplus (\mathbf{1},\mathbf{3})_{0}  \notag
\\
&&\oplus (\mathbf{27},\mathbf{2})_{1}\oplus (\mathbf{27}^{\prime },\mathbf{2}%
)_{-1}\oplus (\mathbf{1},\mathbf{2})_{3}\oplus (\mathbf{1},\mathbf{2})_{-3};
\end{eqnarray}%
\begin{eqnarray}
2 &:&E_{8(8)}\supset ^{ns}E_{6(6)}\oplus sl(3,\mathbb{R})\supset
E_{6(6)}\oplus sl(2,\mathbb{R})\oplus so(1,1); \\
\mathbf{248} &=&(\mathbf{78},\mathbf{1})\oplus (\mathbf{1},\mathbf{8})\oplus
(\mathbf{27},\mathbf{3})\oplus (\mathbf{27}^{\prime },\mathbf{3}^{\prime })
\notag \\
&=&(\mathbf{78},\mathbf{1})_{0}\oplus (\mathbf{1},\mathbf{3})_{0}\oplus (%
\mathbf{1},\mathbf{1})_{0}\oplus (\mathbf{1},\mathbf{2})_{3}\oplus (\mathbf{1%
},\mathbf{2})_{-3}  \notag \\
&&\oplus (\mathbf{27},\mathbf{1})_{-2}\oplus (\mathbf{27},\mathbf{2}%
)_{1}\oplus (\mathbf{27}^{\prime },\mathbf{1})_{2}\oplus (\mathbf{27}%
^{\prime },\mathbf{2})_{1}.
\end{eqnarray}%
Both chains of embeddings $1$ and $2$ give rise to a $7$-grading, with ~$\cn%
_{3\,;\,a)}^{\,8(8)+}\cong \cn_{3\,;\,b)}^{\,8(8)+}=(\mathbf{1},\mathbf{2}%
)_{3}\oplus \left( \mathbf{27}^{\prime },\mathbf{1}\right) _{2}\oplus (%
\mathbf{27},\mathbf{2})_{1}$, with real dimension $83$. By recalling \textit{%
e.g. }(\ref{interpr-b}), the Jordan-algebraic interpretation of the chains $1
$ and $2$ goes as follows :%
\begin{equation}
1:\left\{
\begin{array}{l}
qconf\left( J_{3}^{\mathbb{O}_{s}}\right) \supset conf\left( J_{3}^{\mathbb{O%
}_{s}}\right) \oplus sl(2,\mathbb{R})_{Ehlers}\supset str_{0}\left( J_{3}^{%
\mathbb{O}_{s}}\right) \oplus sl(2,\mathbb{R})_{Ehlers}\oplus so(1,1)_{KK};
\\
\text{\textit{or}} \\
qconf\left( J_{3}^{\mathbb{O}_{s}}\right) \supset qconf\left( J_{3}^{\mathbb{%
H}_{s}}\right) \oplus \widetilde{\mathcal{A}}_{q=4}\supset qconf\left(
J_{3}^{\mathbb{C}_{s}}\right) \oplus \widetilde{\mathcal{A}}_{4}\oplus
\widetilde{\mathcal{A}}_{2};%
\end{array}%
\right.   \label{j816a}
\end{equation}%
\begin{equation}
2:\left\{
\begin{array}{l}
qconf\left( J_{3}^{\mathbb{O}_{s}}\right) \supset ^{ns}str_{0}\left( J_{3}^{%
\mathbb{O}_{s}}\right) \oplus sl(3,\mathbb{R})_{Ehlers}\supset str_{0}\left(
J_{3}^{\mathbb{O}_{s}}\right) \oplus sl(2,\mathbb{R})\oplus so(1,1); \\
\text{\textit{or}} \\
qconf\left( J_{3}^{\mathbb{O}_{s}}\right) \supset ^{ns}qconf\left( J_{3}^{%
\mathbb{C}_{s}}\right) \oplus sl(3,\mathbb{R})\supset qconf\left( J_{3}^{%
\mathbb{C}_{s}}\right) \oplus sl(2,\mathbb{R})\oplus so(1,1).%
\end{array}%
\right.   \label{j816b}
\end{equation}

\md

\subsubsection{$\cp_{7}^{8(8)}$\label{817}}

The maximal parabolics ~$\cp_{7}^{8(8)}$ ~ from \eqref{maxeighta}
corresponds to the Bruhat decomposition:
\begin{equation}
E_{8(8)}~=~\cn_{7}^{\,8(8)-}\oplus E_{7(7)}\oplus so(1,1)\oplus \cn%
_{7}^{\,8(8)+}\newline
\mathbf{,}
\end{equation}%
which can be obtained through the embedding chain
\begin{eqnarray}
E_{8(8)} &\supset &E_{7(7)}\oplus sl(2,\mathbb{R})\supset E_{7(7)}\oplus
so(1,1); \\
\mathbf{248} &=&(\mathbf{133},\mathbf{1})\oplus (\mathbf{1},\mathbf{3}%
)\oplus (\mathbf{56},\mathbf{2})=\mathbf{133}_{0}\oplus \mathbf{1}_{0}\oplus
\mathbf{1}_{2}\oplus \mathbf{1}_{-2}\oplus \mathbf{56}_{1}\oplus \mathbf{56}%
_{-1},
\end{eqnarray}%
giving rise to a $3$-grading, with $\cn_{7}^{\,8(8)+}=\mathbf{1}_{2}\oplus
\mathbf{56}_{1}$, with real dimension $57$. The Jordan algebraic
interpretation of the first step is the same as the chain $1$ of case \ref%
{816} above.

\md

\subsubsection{$\cp_{8}^{8(8)}$\label{818}}

The maximal parabolics ~$\cp_{8}^{8(8)}$ ~ from \eqref{maxeighta}
corresponds to the Bruhat decomposition:
\begin{equation}
E_{8(8)}~=~\cn_{8}^{\,8(8)-}sl(8,\bbr)\oplus so(1,1)\oplus \oplus \cn%
_{8}^{\,8(8)+}\newline
\mathbf{,}
\end{equation}%
which can be obtained through the embedding chain
\begin{eqnarray}
&&E_{8(8)}\supset ^{ns}sl(9,\mathbb{R})\supset sl(8,\mathbb{R})\oplus
so(1,1); \\
\mathbf{248} &=&\mathbf{80}\oplus \mathbf{84\oplus 84}^{\prime }=\mathbf{63}%
_{0}\oplus \mathbf{1}_{0}\oplus \mathbf{8}_{9}\oplus \mathbf{8}_{-9}^{\prime
}\oplus \mathbf{28}_{-6}\oplus \mathbf{56}_{3}\oplus \mathbf{28}_{6}^{\prime
}\oplus \mathbf{56}_{-3}^{\prime },  \notag
\end{eqnarray}%
giving rise to a $7$-grading, with $\cn_{8}^{\,8(8)+}=\mathbf{8}_{9}\oplus
\mathbf{28}_{6}^{\prime }\oplus \mathbf{56}_{3}$, with real dimension $92$.

The Jordan algebraic interpretation of the chain goes as follows (recall (%
\ref{jj}))
\begin{eqnarray}
&&qconf\left( J_{3}^{\mathbb{O}_{s}}\right) \supset ^{ns}g_{D=11}\oplus
\left. sl(D-2,\mathbb{R})\right\vert _{D=11}  \notag \\
&\supset &g_{D=11}\oplus \left. sl(D-2,\mathbb{R})\right\vert _{D=10}\oplus
so(1,1)  \notag \\
&&\cong g_{D=10,IIA}\oplus \left. sl(D-2,\mathbb{R})\right\vert _{D=10}.
\end{eqnarray}

\md

\subsection{$E_{8(-24)}$}

This is the minimally non-compact real form of $E_{8}$. Its Jordan
interpretation reads:%
\begin{equation}
E_{8(-24)}\cong qconf\left( J_{3}^{\mathbb{O}}\right) .
\end{equation}

\md

\subsubsection{$\cp_{1}^{8(-24)}$\label{821}}

The maximal parabolics ~$\cp_{1}^{8(-24)}$ ~ from \eqref{maxeightb}
corresponds to the Bruhat decomposition:
\begin{equation}
E_{8(-24)}~=~\cn_{1}^{\,8(-24)-}\oplus so(11,3)\oplus so(1,1)\oplus \cn%
_{1}^{\,8(-24)+}\newline
\mathbf{,}
\end{equation}

which can be obtained through the embedding chain
\begin{eqnarray}
E_{8(-24)} &\supset &so(12,4)\supset so(11,3)\oplus so(1,1); \\
\mathbf{248} &=&\mathbf{120}\oplus \mathbf{128}=\mathbf{91}_{0}\oplus
\mathbf{1}_{0}\oplus \mathbf{14}_{2}\oplus \mathbf{14}_{-2}\oplus \mathbf{64}%
_{-1}\oplus \mathbf{64}_{1}^{\prime },
\end{eqnarray}%
giving rise to a $5$-grading, with $\cn_{8}^{\,8(-24)+}=\mathbf{14}%
_{2}\oplus \mathbf{64}_{1}^{\prime }$. Since
\begin{equation}
so(12,4)\cong qconf\left( \mathbb{R}\oplus J_{2}^{\mathbb{O}}\right) ,
\label{so(12,4)}
\end{equation}%
a Jordan algebraic interpretation (of the first step) of the chain above is
given by the embedding $J_{3}^{\mathbb{O}}\supset \mathbb{R}\oplus J_{2}^{%
\mathbb{O}}$ considered at the $qconf$ level :%
\begin{equation}
qconf\left( J_{3}^{\mathbb{O}}\right) \supset qconf\left( \mathbb{R}\oplus
J_{2}^{\mathbb{O}}\right) \oplus \mathcal{A}_{8}.  \label{j821}
\end{equation}

\md

\subsubsection{$\cp_{2}^{8(-24)}$\label{822}}

The maximal parabolics ~$\cp_{2}^{8(-24)}$ ~ from \eqref{maxeightb}
corresponds to the Bruhat decomposition:
\begin{equation}
E_{8(8)}~=\cn_{2}^{\,8(-24)-}\oplus so(9,1)\oplus sl(3,\bbr)\oplus
so(1,1)\oplus \cn_{2}^{\,8(-24)+}\newline
\mathbf{,}
\end{equation}%
which can be obtained through \textit{at least} two chains of embeddings,
respectively denoted by $1$ and $2$ (recall that $so(3,3)\cong sl(4,\mathbb{R%
})$) :%
\begin{eqnarray}
1 &:&E_{8(-24)}\supset so(12,4)  \notag \\
&\supset &so(9,1)\oplus so(3,3)\supset so(9,1)\oplus sl(3,\mathbb{R})\oplus
so(1,1); \\
\mathbf{248} &=&\mathbf{120}\oplus \mathbf{128}=(\mathbf{45},\mathbf{1}%
)\oplus (\mathbf{1},\mathbf{15})\oplus (\mathbf{10},\mathbf{6})\oplus (%
\mathbf{16},\mathbf{4})\oplus (\mathbf{16}^{\prime },\mathbf{4}^{\prime })
\notag \\
&=&(\mathbf{45},\mathbf{1})_{0}\oplus (\mathbf{1},\mathbf{1})_{0}\oplus (%
\mathbf{1},\mathbf{3})_{4}\oplus (\mathbf{1},\mathbf{3}^{\prime
})_{-4}\oplus (\mathbf{1},\mathbf{8})_{0}  \notag \\
&&\oplus (\mathbf{10},\mathbf{3})_{-2}\oplus (\mathbf{10},\mathbf{3}^{\prime
})_{2}\oplus (\mathbf{16},\mathbf{1})_{-3}\oplus \left( \mathbf{16},\mathbf{3%
}\right) _{1}\oplus (\mathbf{16}^{\prime },\mathbf{1})_{3}\oplus (\mathbf{16}%
^{\prime },\mathbf{3}^{\prime })_{-1};
\end{eqnarray}%
\begin{eqnarray}
2 &:&E_{8(-24)}\supset ^{ns}E_{6(-26)}\oplus sl(3,\mathbb{R})\supset
so(9,1)\oplus sl(3,\mathbb{R})\oplus so(1,1); \\
\mathbf{248} &=&(\mathbf{78},\mathbf{1})\oplus (\mathbf{1},\mathbf{8})\oplus
(\mathbf{27},\mathbf{3})\oplus (\mathbf{27}^{\prime },\mathbf{3}^{\prime })
\notag \\
&=&(\mathbf{45},\mathbf{1})_{0}\oplus (\mathbf{1},\mathbf{1})_{0}\oplus (%
\mathbf{16},\mathbf{1})_{3}\oplus (\mathbf{16}^{\prime },\mathbf{1}%
)_{-3}\oplus (\mathbf{1},\mathbf{1},\mathbf{8})_{0}  \notag \\
&&\oplus (\mathbf{1},\mathbf{3})_{-4}\oplus (\mathbf{10},\mathbf{3}%
)_{2}\oplus (\mathbf{16},\mathbf{3})_{-1}\oplus (\mathbf{1},\mathbf{3}%
^{\prime })_{4}\oplus (\mathbf{10},\mathbf{3}^{\prime })_{-2}\oplus (\mathbf{%
16}^{\prime },\mathbf{3}^{\prime })_{1}.
\end{eqnarray}%
Both chains $1$ and $2$ give rise to a $9$-grading, with ~$\cn%
_{2\,;\,a)}^{\,8(-24)+}=(\mathbf{1},\mathbf{3})_{4}\oplus (\mathbf{16}%
^{\prime },\mathbf{1})_{3}\oplus (\mathbf{10},\mathbf{3}^{\prime
})_{2}\oplus (\mathbf{16},\mathbf{3})_{1}$, ~$\cn_{2\,;\,b)}^{\,8(-24)+}=(%
\mathbf{1},\mathbf{3}^{\prime })_{4}\oplus (\mathbf{16},\mathbf{1}%
)_{3}\oplus (\mathbf{10},\mathbf{3})_{2}\oplus (\mathbf{16}^{\prime },%
\mathbf{3}^{\prime })_{1}$, both of real dimension $97$. Note that $%
1\leftrightarrow 2$ \textit{iff} the weights of the parabolic $so(1,1)$'s
gets flipped (or, equivalently, \textit{iff} all $\mathcal{M}_{5}^{8(8)}$%
-irreps. get conjugated). The Jordan-algebraic interpretation of chains $1$
and $2$ is based on the identification (\ref{so(12,4)}) and on%
\begin{equation}
so(9,1)\cong str_{0}(J_{2}^{\mathbb{O}})\cong str_{0}(\mathbb{R}\oplus
J_{2}^{\mathbb{O}})\ominus so(1,1).
\end{equation}%
Note that $so(9,1)$, differently from the split form $so(5,5)$ (\textit{cfr.}
(\ref{qc})), does not admit a quasi-conformal interpretation. The first step
of chain $1$ is then a consequence of the embedding $J_{3}^{\mathbb{O}%
}\supset \mathbb{R}\oplus J_{2}^{\mathbb{O}}$, considered at the $qconf$
level. Note that the \textit{next-to-maximal} (\textit{non-symmetric})
embedding%
\begin{equation}
E_{8(-24)}\supset so(9,1)\oplus so(3,3)\cong so(9,1)\oplus sl(4,\mathbb{R})
\end{equation}%
is the $D=6$ case of the Ehlers embedding \cite{super-Ehlers} for the $%
\mathbb{O}$-based theory, thus characterizing $sl(4,\mathbb{R})$ as the $D=6$
Ehlers symmetry (\textit{cfr.} (\ref{aa-2}) below). Concerning the chain $2$%
, its first step is the JP embedding for the $\mathbb{O}$-based theory, thus
giving rise to the $D=5$ Ehlers $sl(3,\mathbb{R})_{Ehlers}$; then, in the
second step a further uplift to $D=6$ is performed, introducing the KK $%
so(1,1)_{KK}$ of the $S^{1}$-reduction $D=6\rightarrow 5$, which is the
parabolic $so(1,1)$ in this chain:
\begin{eqnarray}
1 &:&qconf\left( J_{3}^{\mathbb{O}}\right) \supset qconf\left( \mathbb{R}%
\oplus J_{2}^{\mathbb{O}}\right) \oplus \mathcal{A}_{8}\supset
str_{0}(J_{2}^{\mathbb{O}})\oplus \mathcal{A}_{8}\oplus sl(4,\mathbb{R}%
)_{Ehlers}  \notag \\
&&  \notag \\
&\supset &str_{0}(J_{2}^{\mathbb{O}})\oplus \mathcal{A}_{8}\oplus sl(3,%
\mathbb{R})\oplus so(1,1);  \label{aa-2}
\end{eqnarray}%
\begin{equation}
2:qconf(J_{3}^{\mathbb{O}})\supset ^{ns}str_{0}(J_{3}^{\mathbb{O}})\oplus
sl(3,\mathbb{R})_{Ehlers}\supset str_{0}(J_{2}^{\mathbb{O}})\oplus \mathcal{A%
}_{8}\oplus sl(3,\mathbb{R})_{Ehlers}\oplus so(1,1)_{KK}.  \label{bb-2}
\end{equation}

\md

\subsubsection{$\cp_{3}^{8(-24)}$\label{823}}

The maximal parabolics ~$\cp_{3}^{8(-24)}$ ~ from \eqref{maxeightb}
corresponds to the Bruhat decomposition:
\begin{equation}
E_{8(8)}~=\cn_{3}^{\,8(-24)-}\oplus E_{6(-26)}\oplus sl(2,\bbr)\oplus
so(1,1)\oplus \cn_{3}^{\,8(-24)+}\newline
\mathbf{,}
\end{equation}%
which can be obtained through \textit{at least} two chains of embeddings,
respectively denoted by $1$ and $2$ :%
\begin{eqnarray}
1 &:&E_{8(-24)}\supset E_{7(-25)}\oplus sl(2,\mathbb{R})\supset
E_{6(-26)}\oplus sl(2,\mathbb{R})\oplus so(1,1); \\
\mathbf{248} &=&(\mathbf{133},\mathbf{1})\oplus (\mathbf{1},\mathbf{3}%
)\oplus (\mathbf{56},\mathbf{2})  \notag \\
&=&\left( \mathbf{78,1}\right) _{0}\oplus \left( \mathbf{1,1}\right)
_{0}\oplus \left( \mathbf{27,1}\right) _{-2}\oplus \left( \mathbf{27}%
^{\prime },\mathbf{1}\right) _{2}\oplus (\mathbf{1},\mathbf{3})_{0}  \notag
\\
&&\oplus (\mathbf{27},\mathbf{2})_{1}\oplus (\mathbf{27}^{\prime },\mathbf{2}%
)_{-1}\oplus (\mathbf{1},\mathbf{2})_{3}\oplus (\mathbf{1},\mathbf{2})_{-3};
\end{eqnarray}%
\begin{eqnarray}
2 &:&E_{8(-24)}\supset ^{ns}E_{6(-26)}\oplus sl(3,\mathbb{R})\supset
E_{6(-26)}\oplus sl(2,\mathbb{R})\oplus so(1,1); \\
\mathbf{248} &=&(\mathbf{78},\mathbf{1})\oplus (\mathbf{1},\mathbf{8})\oplus
(\mathbf{27},\mathbf{3})\oplus (\mathbf{27}^{\prime },\mathbf{3}^{\prime })
\notag \\
&=&(\mathbf{78},\mathbf{1})_{0}\oplus (\mathbf{1},\mathbf{3})_{0}\oplus (%
\mathbf{1},\mathbf{1})_{0}\oplus (\mathbf{1},\mathbf{2})_{3}\oplus (\mathbf{1%
},\mathbf{2})_{-3}  \notag \\
&&\oplus (\mathbf{27},\mathbf{1})_{-2}\oplus (\mathbf{27},\mathbf{2}%
)_{1}\oplus (\mathbf{27}^{\prime },\mathbf{1})_{2}\oplus (\mathbf{27}%
^{\prime },\mathbf{2})_{1}.
\end{eqnarray}%
Both chains $1$ and $2$ give rise to a $7$-grading, with $\cn%
_{2\,;\,1}^{\,8(-24)+}\cong \cn_{2\,;\,1}^{\,8(-24)+}=(\mathbf{1},\mathbf{2}%
)_{3}\oplus (\mathbf{27}^{\prime },\mathbf{1})_{2}\oplus (\mathbf{27}%
^{\prime },\mathbf{2})_{1}$, with real dimension $83$.

The Jordan algebraic interpretation of the chains $1$ and $2$ goes as
follows :%
\begin{equation}
1:qconf\left( J_{3}^{\mathbb{O}}\right) \supset conf\left( J_{3}^{\mathbb{O}%
}\right) \oplus sl(2,\mathbb{R})_{Ehlers}\supset str_{0}\left( J_{3}^{%
\mathbb{O}}\right) \oplus sl(2,\mathbb{R})_{Ehlers}\oplus so(1,1)_{KK};
\label{j823a}
\end{equation}%
\begin{equation}
2:qconf\left( J_{3}^{\mathbb{O}}\right) \supset ^{ns}str_{0}\left( J_{3}^{%
\mathbb{O}}\right) \oplus sl(3,\mathbb{R})_{Ehlers}\supset str_{0}\left(
J_{3}^{\mathbb{O}}\right) \oplus sl(2,\mathbb{R})\oplus so(1,1).
\label{j823b}
\end{equation}

\md

\subsubsection{$\cp_{4}^{8(-24)}$\label{824}}

The maximal parabolics ~$\cp_{4}^{8(-24)}$ ~ from \eqref{maxeightb}
corresponds to the Bruhat decomposition:
\begin{equation}
E_{8(8)}~=~\cn_{4}^{\,8(-24)-}\oplus E_{7(-25)}\oplus so(1,1)\oplus \cn%
_{4}^{\,8(-24)+}\newline
\mathbf{,}
\end{equation}%
which can be obtained through the embedding chain
\begin{eqnarray}
E_{8(-24)} &\supset &E_{7(-25)}\oplus sl(2,\mathbb{R})\supset
E_{7(-25)}\oplus so(1,1); \\
\mathbf{248} &=&(\mathbf{133},\mathbf{1})\oplus (\mathbf{1},\mathbf{3}%
)\oplus (\mathbf{56},\mathbf{2})=\mathbf{133}_{0}\oplus \mathbf{1}_{0}\oplus
\mathbf{1}_{2}\oplus \mathbf{1}_{-2}\oplus \mathbf{56}_{1}\oplus \mathbf{56}%
_{-1},
\end{eqnarray}%
giving rise to a $5$-grading, with $\cn_{4}^{\,8(-24)+}=\mathbf{1}_{2}\oplus
\mathbf{56}_{1}$, with real dimension $57$.

The Jordan algebraic interpretation of the first step is the same as the
chain $1$ of subsubsection \ref{823}, and it is analogous to the chain from
the parabolically related case ~$\cp_{7}^{8(8)}\,$, Subsubsection \ref{817}.

\md

\subsection{$F_{4(4)}$}

This is the split real form of $F_{4}$. Its Jordan interpretation is twofold
(due to the symmetry of the double-split Magic Square ${\mathcal{L}}_{3}(%
\mathbb{A}_{s},\mathbb{B}_{s})$ \cite{BS}, reported in Table 3) :%
\begin{eqnarray}
F_{4(4)} &\cong &qconf\left( J_{3}^{\mathbb{R}}\right)   \label{unos} \\
&\cong &der\left( J_{3}^{\mathbb{O}_{s}}\right) .  \label{dues}
\end{eqnarray}

\md

\subsubsection{$\cp_{1}^{4(4)}$\label{411}}

The maximal parabolics ~$\cp_{1}^{4(4)}$ ~ from \eqref{maxfoura} corresponds
to the Bruhat decomposition:
\begin{equation}
F_{4(4)}~=~\cn_{1}^{\,4(4)-}\oplus sl(3,\bbr)_{S}\oplus sl(2,\bbr)_{L}\oplus
so(1,1)\oplus \cn_{1}^{\,4(4)+}\newline
\mathbf{,}
\end{equation}%
which can be obtained through \textit{at least} two chains of embeddings
respectively denoted by $1$ and $2$ :
\begin{eqnarray}
1 &:&F_{4(4)}\supset so(5,4)  \notag \\
&\supset &so(3,3)\oplus so(2,1)\cong sl(4,\mathbb{R})\oplus sl(2,\mathbb{R}%
)\supset sl(3,\mathbb{R})_{S}\oplus sl(2,\mathbb{R})_{L}\oplus so(1,1);
\notag \\
&& \\
\mathbf{52} &=&\mathbf{36}\oplus \mathbf{16}=(\mathbf{15},\mathbf{1})\oplus (%
\mathbf{1},\mathbf{3})\oplus (\mathbf{6},\mathbf{3})\oplus (\mathbf{4},%
\mathbf{2})\oplus (\mathbf{4}^{\prime },\mathbf{2})  \notag \\
&=&(\mathbf{8},\mathbf{1})_{0}\oplus (\mathbf{1},\mathbf{3})_{0}\oplus (%
\mathbf{1},\mathbf{1})_{0}\oplus (\mathbf{3}^{\prime },\mathbf{1}%
)_{-4}\oplus (\mathbf{1},\mathbf{2})_{-3}\oplus   \notag \\
&&\oplus (\mathbf{3},\mathbf{3})_{-2}\oplus (\mathbf{3}^{\prime },\mathbf{2}%
)_{-1}\oplus (\mathbf{3},\mathbf{1})_{4}\oplus (\mathbf{1},\mathbf{2}%
)_{3}\oplus (\mathbf{3}^{\prime },\mathbf{3})_{2}\oplus (\mathbf{3},\mathbf{2%
})_{1};
\end{eqnarray}%
\begin{equation*}
2:F_{4(4)}\supset ^{ns}sl(3,\mathbb{R})_{L}\oplus sl(3,\mathbb{R}%
)_{S}\supset sl(3,\mathbb{R})_{S}\oplus sl(2,\mathbb{R})_{R}\oplus so(1,1);
\end{equation*}%
\begin{equation*}
\mathbf{52}=(\mathbf{8},\mathbf{1})\oplus (\mathbf{1},\mathbf{8})\oplus (%
\mathbf{6},\mathbf{3}^{\prime })\oplus (\mathbf{6}^{\prime },\mathbf{3})=%
\begin{array}{l}
(\mathbf{8},\mathbf{1})_{0}\oplus (\mathbf{1},\mathbf{3})_{0}\oplus (\mathbf{%
1},\mathbf{1})_{0}\oplus (\mathbf{3}^{\prime },\mathbf{1})_{-4}\ \oplus  \\
\oplus (\mathbf{1},\mathbf{2})_{-3}\oplus (\mathbf{3},\mathbf{3})_{-2}\oplus
(\mathbf{3}^{\prime },\mathbf{2})_{-1}\oplus (\mathbf{3},\mathbf{1})_{4}\
\oplus  \\
\oplus (\mathbf{1},\mathbf{2})_{3}\oplus (\mathbf{3}^{\prime },\mathbf{3}%
)_{2}\oplus (\mathbf{3},\mathbf{2})_{1};%
\end{array}%
\end{equation*}%
Both chains $1$ and $2$ give rise to a $9$-grading, with ~$\cm%
_{1}^{4(4)}~\cong ~sl(3,\bbr)_{S}\oplus sl(2,\bbr)_{L}\,$ and ~$\cn%
_{1\,;\,1}^{\,4(4)+}\cong \cn_{1\,;\,2}^{\,4(4)+}=(\mathbf{3},\mathbf{1}%
)_{4}\oplus (\mathbf{1},\mathbf{2})_{3}\oplus (\mathbf{3}^{\prime },\mathbf{3%
})_{2}\oplus (\mathbf{3},\mathbf{2})_{1}$, of real dimension $20$.

Since
\begin{equation}
so(5,4)\cong qconf\left( \mathbb{R}\oplus J_{2}^{\mathbb{R}}\right) ,
\end{equation}%
a Jordan algebraic interpretation (of the first step) of the chain $1$ is
given by the embedding $J_{3}^{\mathbb{R}}\supset \mathbb{R}\oplus J_{2}^{%
\mathbb{R}}$ considered at the $qconf$ level :%
\begin{equation}
1:qconf\left( J_{3}^{\mathbb{R}}\right) \supset qconf\left( \mathbb{R}\oplus
J_{2}^{\mathbb{R}}\right) \oplus \mathcal{A}_{1}.
\end{equation}%
The second step realizes (with respect to $F_{4(4)}$) the $D=6$ case of the
Ehlers embedding for the $\mathbb{R}$-based theory, thus giving rise to the $%
D=6$ Ehlers symmetry $sl(4,\mathbb{R})_{Ehlers}$ (which then gets further
branched in order to generate the parabolic $so(1,1)$).

Finally, the first step of the chain $2$ can be interpreted as the JP
embedding for the $\mathbb{R}$-based theory, where $sl(3,\mathbb{R}%
)_{L}\cong str_{0}\left( J_{3}^{\mathbb{R}}\right) $, and $sl(3,\mathbb{R}%
)_{S}\cong sl(3,\mathbb{R})_{Ehlers}$ is the $D=5$ Ehlers symmetry.
Furthermore, $sl(3,\mathbb{R})_{L}\cong str_{0}\left( J_{3}^{\mathbb{R}%
}\right) $ branches into $sl(2,\mathbb{R})_{L}\cong so(2,1)\cong
str_{0}\left( J_{2}^{\mathbb{R}}\right) \cong str_{0}\left( \mathbb{R}\oplus
J_{2}^{R}\right) \ominus so(1,1)$, thus allowing for the identification of
the parabolic $so(1,1)$ with the $so(1,1)_{KK}$ of the $D=6\rightarrow 5$ $%
S^{1}$-reduction.

%###

\md

\subsubsection{$\cp_{2}^{4(4)}$\label{412}}

The maximal parabolics ~$\cp_{2}^{4(4)}$ ~ from \eqref{maxfoura} corresponds
to the Bruhat decomposition:
\begin{equation}
F_{4(4)}~=\cn_{2}^{\,4(4)-}\oplus sl(3,\bbr)_{L}\oplus sl(2,\bbr)_{S}\oplus
so(1,1)\oplus \cn_{2}^{\,4(4)+}\newline
\mathbf{,}
\end{equation}%
which can be obtained through \textit{at least} two chains of embeddings,
respectively denoted by $1$ and $2$ :%
\begin{eqnarray}
1 &:&F_{4(4)}\supset sp(3,\mathbb{R})\oplus sl(2,\mathbb{R})\supset sl(3,%
\mathbb{R})_{L}\oplus sl(2,\mathbb{R})_{S}\oplus so(1,1); \\
\mathbf{52} &=&(\mathbf{21},\mathbf{1})\oplus (\mathbf{1},\mathbf{3})\oplus (%
\mathbf{14}^{\prime },\mathbf{2})  \notag \\
&=&(\mathbf{8},\mathbf{1})_{0}\oplus (\mathbf{1},\mathbf{3})_{0}\oplus (%
\mathbf{1},\mathbf{1})_{0}\oplus (\mathbf{1},\mathbf{2})_{-3}\oplus (\mathbf{%
6}^{\prime },\mathbf{1})_{-2}\oplus   \notag \\
&&\oplus (\mathbf{6},\mathbf{2})_{-1}\oplus (\mathbf{1},\mathbf{2}%
)_{3}\oplus (\mathbf{6},\mathbf{1})_{2}\oplus (\mathbf{6}^{\prime },\mathbf{2%
})_{1};
\end{eqnarray}%
\begin{equation*}
2:F_{4(4)}\supset ^{ns}sl(3,\mathbb{R})_{L}\oplus sl(3,\mathbb{R}%
)_{S}\supset sl(3,\mathbb{R})_{L}\oplus sl(2,\mathbb{R})_{S}\oplus so(1,1);
\end{equation*}%
\begin{equation*}
\mathbf{52}=(\mathbf{8},\mathbf{1})\oplus (\mathbf{1},\mathbf{8})\oplus (%
\mathbf{6},\mathbf{3}^{\prime })\oplus (\mathbf{6}^{\prime },\mathbf{3})=%
\begin{array}{l}
(\mathbf{8},\mathbf{1})_{0}\oplus (\mathbf{1},\mathbf{3})_{0}\oplus (\mathbf{%
1},\mathbf{1})_{0}\oplus (\mathbf{1},\mathbf{2})_{-3}\ \oplus  \\
\oplus (\mathbf{6}^{\prime },\mathbf{1})_{-2}\oplus (\mathbf{6},\mathbf{2}%
)_{-1}\oplus (\mathbf{1},\mathbf{2})_{3}\ \oplus  \\
\oplus (\mathbf{6},\mathbf{1})_{2}\oplus (\mathbf{6}^{\prime },\mathbf{2}%
)_{1};%
\end{array}%
\end{equation*}%
Both chains of embeddings $1$ and $2$ give rise to a $7$-grading, with%
\newline
~$\cm_{2}^{4(4)}~\cong ~sl(3,\bbr)_{L}\oplus sl(2,\bbr)_{S}\,$ and ~$\cn%
_{2\,;\,b)}^{\,4(4)+}\cong \cn_{2\,;\,c.1)}^{\,4(4)+}=(\mathbf{1},\mathbf{2}%
)_{3}\oplus (\mathbf{6},\mathbf{1})_{2}\oplus (\mathbf{6}^{\prime },\mathbf{2%
})_{1}$, of real dimension $20$.

Since%
\begin{eqnarray}
sp(3,\mathbb{R}) &\cong &conf\left( J_{3}^{\mathbb{R}}\right) \cong
der\left( \mathbf{F}\left( J_{3}^{\mathbb{R}}\right) \right) ; \\
sl(3,\mathbb{R}) &\cong &str_{0}\left( J_{3}^{\mathbb{R}}\right) ,
\end{eqnarray}%
the Jordan algebraic interpretation of the chain $1$ reads:%
\begin{equation}
1:qconf\left( J_{3}^{\mathbb{R}}\right) \supset conf\left( J_{3}^{\mathbb{R}%
}\right) \oplus sl(2,\mathbb{R})_{Ehlers}\supset str_{0}\left( J_{3}^{%
\mathbb{R}}\right) \oplus sl(2,\mathbb{R})_{Ehlers}\oplus so(1,1),
\end{equation}%
where $sl(2,\mathbb{R})_{Ehlers}$ is the $D=4$ Ehlers symmetry.

Finally, the first step of the chain $2$ can be interpreted as the JP
embedding for the $\mathbb{R}$-based theory, where $sl(3,\mathbb{R}%
)_{L}\cong str_{0}\left( J_{3}^{\mathbb{R}}\right) $, and $sl(3,\mathbb{R}%
)_{S}\cong sl(3,\mathbb{R})_{Ehlers}$ is the $D=5$ Ehlers symmetry. (This is
the same as in Subsubsection \ref{411}.) Then, $sl(3,\mathbb{R})_{S}=sl(3,%
\mathbb{R})_{Ehlers}$ branches to ~$sl(2,\mathbb{R})_{S}$~ to generate the
parabolic $so(1,1)$.

\md

\subsubsection{$\cp_{3}^{4(4)}$\label{413}}

The maximal parabolics ~$\cp_{3}^{4(4)}$ ~ from \eqref{maxfoura} corresponds
to the Bruhat decomposition:
\begin{equation}
F_{4(4)}~=~\cn_{3}^{\,4(4)-}\oplus sp(3,\bbr)\oplus so(1,1)\oplus \cn%
_{3}^{\,4(4)+}\newline
\mathbf{,}
\end{equation}%
which can be obtained through the embedding chain
\begin{eqnarray}
F_{4(4)} &\supset &sp(3,\mathbb{R})\oplus sl(2,\mathbb{R})\supset sp(3,%
\mathbb{R})\oplus so(1,1); \\
\mathbf{52} &=&(\mathbf{21},\mathbf{1})\oplus (\mathbf{1},\mathbf{3})\oplus (%
\mathbf{14}^{\prime },\mathbf{2})  \notag \\
&=&\mathbf{21}_{0}\oplus \mathbf{1}_{0}\oplus \mathbf{1}_{2}\oplus \mathbf{1}%
_{-2}\oplus \mathbf{14}_{1}^{\prime }\oplus \mathbf{14}_{-1}^{\prime },
\end{eqnarray}%
giving rise to a $5$-grading, with ~ $\cn_{3}^{\,4(4)+}=\mathbf{1}_{2}\oplus
\mathbf{14}_{1}^{\prime }$, with real dimension $15$.

The Jordan algebraic interpretation of the first step is the same as the
chain $1$ of subsubsection \ref{412} above. Thus, the parabolic $so(1,1)$
can here be interpreted as the non-compact Cartan generator of the $D=4$
Ehlers symmetry $sl(2,\mathbb{R})_{Ehlers}$.

\md

\subsubsection{$\cp_{4}^{4(4)}$\label{414}}

The maximal parabolics ~$\cp_{4}^{4(4)}$ ~ from \eqref{maxfoura} corresponds
to the Bruhat decomposition:
\begin{equation}
F_{4(4)}~=~\cn_{4}^{\,4(4)-}\oplus so(4,3)\oplus so(1,1)\oplus \cn%
_{4}^{\,4(4)+}\newline
\mathbf{,}
\end{equation}%
which can be obtained through the embedding chain
\begin{eqnarray}
F_{4(4)} &\supset &so(5,4)\supset so(4,3)\oplus so(1,1); \\
\mathbf{52} &=&\mathbf{36}\oplus \mathbf{16}=\mathbf{21}_{0}\oplus \mathbf{1}%
_{0}\oplus \mathbf{7}_{2}\oplus \mathbf{7}_{-2}\oplus \mathbf{8}_{1}\oplus
\mathbf{8}_{-1},
\end{eqnarray}%
giving rise to a $5$-grading, with ~$\cn_{4}^{\,4(4)+}=\mathbf{7}_{2}\oplus
\mathbf{8}_{1}$, with real dimension $15$.

The Jordan algebraic interpretation of the first step is the same as the
chain $1$ of subsubsection \ref{411} above. Moreover, by observing that%
\begin{equation}
so(4,3)\cong qconf\left( \mathbb{R}\oplus \mathbb{R}\right) ,
\end{equation}%
the second step can be interpreted as the consequence, at the $qconf$ level,
of the embedding $\mathbb{R}\oplus J_{2}^{\mathbb{R}}\supset $ $\mathbb{R}%
\oplus \mathbb{R}$, corresponding to the embedding of the $c$-map \cite%
{c-map} of the so-called $ST^{2}$ model of $\mathcal{N}=2$, $D=4$
supergravity into the $c$-map of the $\left( \mathbb{R}\oplus J_{2}^{\mathbb{%
R}}\cong \mathbb{R}\oplus \Gamma _{1,2}\right) $-based model of $\mathcal{N}%
=2$, $D=4$ supergravity.

\md

\subsection{$F_{4(-20)}$}

The Jordan interpretation of\ $F_{4(-20)}$ reads:%
\begin{equation}
F_{4(-20)}\cong der\left( J_{2,1}^{\mathbb{O}}\right) \cong \mathcal{S}%
\left( J_{3}^{\mathbb{O}}\right) ,
\end{equation}%
where $J_{2,1}^{\mathbb{O}}\cong J_{1,2}^{\mathbb{O}}$ denotes the rank-$3$
Lorentzian Jordan algebra over $\mathbb{O}$ (\textit{cfr. e.g.} \cite%
{GZ,Squaring-Magic}, and Refs. therein); indeed, the non-Euclidean nature of
$J_{2,1}^{\mathbb{O}}$ generally implies the non-compactness of its
automorphism group (differently from the automorphism symmetry of Euclidean
Jordan algebras over division algebras, which is always compact). Moreover, $%
\mathcal{S}\left( J_{3}^{\mathbb{O}}\right) $ denotes the stabilizer group
of the rank-$3$ orbit of the action of $Str_{0}\left( J_{3}^{\mathbb{O}%
}\right) $ on its (fundamental) irrep. $\mathbf{27}$ with non-vanishing
cubic invariant $I_{3}\neq 0$ and representative \textquotedblleft $++-$"
(for further detail, \textit{cfr. e.g.} \cite{FG-D=5,Small-Orbits}; the
relation between $\mathcal{S}$ and $\mathcal{K}$ is investigated in \cite%
{CFM1}).

The maximal parabolics ~$\cp_{2}^{4(-20)}$ ~ from \eqref{maxfourb}
corresponds to the Bruhat decomposition (which is both maximal and minimal):
\begin{equation}
F_{4(-20)}~=~\cn^{\,4(-20)-}\oplus so(7)\oplus so(1,1)\oplus \cn^{\,4(-20)+}%
\newline
\mathbf{,}
\end{equation}%
which can be obtained through the embedding chain
\begin{eqnarray}
F_{4(-20)} &\supset &so(8,1)\supset so(7)\oplus so(1,1); \\
\mathbf{52} &=&\mathbf{36}\oplus \mathbf{16}=\mathbf{21}_{0}\oplus \mathbf{1}%
_{0}\oplus \mathbf{7}_{2}\oplus \mathbf{7}_{-2}\oplus \mathbf{8}_{1}\oplus
\mathbf{8}_{-1},
\end{eqnarray}%
giving rise to a $5$-grading, with $\cn^{\,4(-20)+}=\mathbf{7}_{2}\oplus
\mathbf{8}_{1}$, with real dimension $15$. Note that $so(8,1)$ does not have
a quasi-conformal interpretation. However, we observe that%
\begin{equation}
so(8,1)\cong \mathcal{S}\left( \mathbb{R}\oplus J_{2}^{\mathbb{O}}\right) ,
\end{equation}%
which is the $q=8$ case of the general result%
\begin{equation}
so(q,1)\cong \mathcal{S}\left( \mathbb{R}\oplus J_{2}^{\mathbb{A}}\right) .
\end{equation}%
Therefore, the first step of the chain can be interpreted as the
consequence, at the level of the $\mathcal{S}$-symmetry, of the embedding $%
J_{3}^{\mathbb{R}}\supset $ $\mathbb{R}\oplus J_{2}^{\mathbb{R}}$ :%
\begin{equation}
\mathcal{S}\left( J_{3}^{\mathbb{O}}\right) \supset \mathcal{S}\left(
\mathbb{R}\oplus J_{2}^{\mathbb{O}}\right) \cong \left. so(1,q)\right\vert
_{q=8}\supset \left. so(q-1)\right\vert _{q=8}\oplus so(1,1).
\end{equation}

\subsection{$G_{2(2)}$}

The Jordan interpretation of\ $G_{2(2)}$, split real form of $G_{2}$, reads:%
\begin{equation}
G_{2(2)}\cong qconf\left( \mathbb{R}\right) ,
\end{equation}%
where $\mathbb{R}$ here denotes the real numbers, conceived as the simplest
example of cubic simple Jordan algebra with cubic norm (the parameter $q$
for this case has the \textit{effective} value $q=-2/3$; see \textit{e.g.}
\cite{Vogel,LM-1}, and Refs. therein).

The maximal parabolics ~$\cp_{LS}^{2(2)}$ ~ from \eqref{maxge}
corresponds to the Bruhat decomposition:
\begin{equation}
G_{2(2)}~=~\cn_{LS}^{\,2(2)-}\oplus sl(2,\bbr)_{LS}\oplus so(1,1)\oplus \cn_{LS}^{\,2(2)+}\newline
\mathbf{,}
\end{equation}%
which can be obtained through \textit{at least} two chains of embeddings,
respectively denoted by $1$ and $2$ :%
\begin{eqnarray}
1 &:&G_{2(2)}\supset ^{ns}sl(3,\mathbb{R})_{S}\supset sl(2,\mathbb{R}%
)_{S}\oplus so(1,1); \\
\mathbf{14} &=&\mathbf{8}\oplus \mathbf{3}\oplus \mathbf{3}^{\prime }=%
\mathbf{3}_{0}\oplus \mathbf{1}_{0}\oplus \mathbf{2}_{-3}\oplus \mathbf{1}%
_{-2}\oplus \mathbf{2}_{-1}\oplus \mathbf{2}_{3}\oplus \mathbf{1}_{2}\oplus
\mathbf{2}_{1};
\end{eqnarray}%
\begin{eqnarray}
2 &:&G_{2(2)}\supset sl(2,\mathbb{R})_{L}\oplus sl(2,\mathbb{R})_{S}\supset
\left\{
\begin{array}{l}
2.1:sl(2,\mathbb{R})_{L}\oplus so(1,1); \\
\text{\textit{or}} \\
2.2:sl(2,\mathbb{R})_{S}\oplus so(1,1).%
\end{array}%
\right.  \\
\mathbf{14} &=&(\mathbf{3},\mathbf{1})\oplus (\mathbf{1},\mathbf{3})\oplus (%
\mathbf{4},\mathbf{2})=\left\{
\begin{array}{l}
2.1:\mathbf{3}_{0}\oplus \mathbf{1}_{0}\oplus \mathbf{1}_{-2}\oplus \mathbf{4%
}_{-1}\oplus \mathbf{1}_{2}\oplus \mathbf{4}_{1}; \\
\text{\textit{or}} \\
2.2:\mathbf{3}_{0}\oplus \mathbf{1}_{0}\oplus \mathbf{2}_{-3}\oplus \mathbf{1%
}_{-2}\oplus \mathbf{2}_{-1}\oplus \mathbf{2}_{3}\oplus \mathbf{1}_{2}\oplus
\mathbf{2}_{1};%
\end{array}%
\right.
\end{eqnarray}%
The chains of embeddings $1$ and $2.2$ give rise to a $7$-grading, with ~$\cm%
_{S}^{2(2)}=sl(2,\bbr)_{S}$~ and ~$\cn_{1)}^{\,2(2)+}\cong \cn%
_{2.2)}^{\,2(2)+}\cong \cn_{S}^{\,2(2)+}=\mathbf{2}_{3}\oplus \mathbf{1}%
_{2}\oplus \mathbf{2}_{1}\,$, with real dimension $5$. On the other hand,
the chain of embeddings $2.1$ gives rise to a $5$-grading, with ~$\cm%
_{L}^{2(2)}=sl(2,\bbr)_{L}$~ and ~$\cn_{2.1)}^{\,2(2)+}\cong \cn%
_{L}^{\,2(2)+}=\mathbf{1}_{2}\oplus \mathbf{4}_{1}\,$, once again with real
dimension $5$.

The Jordan algebraic interpretation of the first step of chain $1$ is
provided by the JP embedding for the $\mathbb{R}$-based theory under
consideration, recalling that $str_{0}\left( \mathbb{R}\right) \cong
\varnothing $. On the other hand, by observing that
\begin{equation}
sl(2,\mathbb{R})\cong conf\left( \mathbb{R}\right) \cong der\left( \mathbf{F}%
\left( \mathbb{R}\right) \right) ,
\end{equation}%
the first step of chain $2$ can be conceived as the Ehlers embedding for the
$\mathbb{R}$-based theory under consideration, followed by two possible
second steps : in $2.1$ the $D=4$ Ehlers group $sl(2,\mathbb{R})_{L}\cong
sl(2,\mathbb{R})_{Ehlers}$ is branched, and thus the parabolic $so(1,1)$ can
be interpreted as its non-compact Cartan; in $2.2$ $sl(2,\mathbb{R}%
)_{S}\cong conf\left( \mathbb{R}\right) $ gets branched, and therefore the
parabolic $so(1,1)$ is the KK $so(1,1)$ of the $S^{1}$-reduction $%
D=5\rightarrow 4$ for the theory at stake (which is the so-called $T^{3}$
model of $\mathcal{N}=2$, $D=4$ supergravity \cite{T^3}).%
\begin{equation}
1:qconf\left( \mathbb{R}\right) \supset ^{ns}str_{0}\left( \mathbb{R}\right)
\oplus sl(3,\mathbb{R})_{Ehlers}\supset str_{0}\left( \mathbb{R}\right)
\oplus sl(2,\mathbb{R})\oplus so(1,1);
\end{equation}%
\begin{equation}
2:qconf\left( \mathbb{R}\right) \supset conf\left( \mathbb{R}\right) \oplus
sl(2,\mathbb{R})_{Ehlers}\supset \left\{
\begin{array}{l}
2.1:conf\left( \mathbb{R}\right) \oplus so(1,1); \\
\text{\textit{or}} \\
2.2:sl(2,\mathbb{R})_{Ehlers}\oplus str_{0}\left( \mathbb{R}\right) \oplus
so(1,1)_{KK}.%
\end{array}%
\right.
\end{equation}

\vskip 1cm

\section{Outlook}

In the present paper we initiated the investigation of the relations between
representation theory and Jordan algebras, focussing on non-compact real
forms of finite-dimensional exceptional Lie algebras. We provided a
derivation of the maximal parabolic subalgebras in terms of chains of
maximal (symmetric or non-symmetric) embeddings of Lie algebras, which were
then interpreted in terms of symmetries of Jordan algebras (in particular,
we focussed on the rank-2 and rank-3 classes). This also allowed to provide
a complete Jordan algebraic characterization (classified in Table B) of the
maximally parabolical relations between exceptional Lie algebras (classified
in Table A).

There is a number of possible venues for further future research. For
instance, in light of the Jordan algebraic interpretation provided in this
paper, the relevance of maximal parabolic subalgebras to the theory of
induced representations might have interesting consequences in
Maxwell-Einstein (super)gravity $\mathcal{N}$-extended theories in various
space-time dimensions. Moreover, the present analysis might be extended to
classical Lie algebras (and to higher-rank Jordan algbras), and slight
generalizations of the definition of maximal parabolic subalgebras are
possible within the Borel-de Siebenthal theory. We plan to deal with such
issues in future works \cite{DobMar2}.

%\np

\bigskip

\section*{Acknowledgments}

\nt
V.K.D. acknowledges partial support from Bulgarian NSF Grant DN-18/1.


\begin{thebibliography}{99}
\bibitem{Jordan:1933a} P.~Jordan, \textit{\"{U}ber die multiplikation
quanten-mechanischer grossen}, \emph{\ }Zschr. f. Phys\emph{.} \textbf{80}
(1933) 285.

\bibitem{Jordan:1933b} P.~Jordan, \textit{\"{U}ber verallgemeinerungsm\"{o}%
glichkeiten des formalismus der quantenmechanik}, Nachr. Ges. Wiss.
Gottingen (1933) 209--214.

\bibitem{Jordan:1933vh} P.~Jordan, J.~von Neumann, and E.~P. Wigner, \textit{%
On an algebraic generalization of the quantum mechanical formalism}, {Ann.
Math. \textbf{35} (1934) no. 1, 29--64}.

\bibitem{Jacobson:1961} N.~Jacobson, \textit{Some groups of transformations
defined by Jordan algebras}, J. Reine Angew. Math. \textbf{207} (1961)
61--85.

\bibitem{Jacobson:1968} N.~Jacobson : \textquotedblleft \textit{Structure
and Representations of Jordan Algebras"}, vol.~\textbf{39}, American
Mathematical Society Colloquium Publications, 1968.

\bibitem{Elkies:1996} N.~Elkies and B.~H. Gross, \textit{The exceptional
cone and the Leech lattice}, Internat. Math. Res. Notices \textbf{14} (1996)
665--698.

\bibitem{Gross:1996} B.~H. Gross, \textit{Groups over }$\mathbb{Z}$, Invent.
Math. \textbf{124} (1996) 263--279.

\bibitem{Krutelevich:2002} S.~Krutelevich, \textit{On a canonical form of a }%
$3\times 3$\textit{\ Herimitian matrix over the ring of integral split
octonions}, {J. Algebra \textbf{253} (2002) no.~2, 276--295}.

\bibitem{Krutelevich:2004} S.~Krutelevich, \textit{Jordan algebras,
exceptional groups, and Bhargava composition}, {J. Algebra \textbf{314}
(2007) no.~2, 924--977}, {\texttt{arXiv:math/0411104}}.

\bibitem{Springer:1962} T.~A. Springer, \textit{Characterization of a class
of cubic forms}, Nederl. Akad. Wetensch. Proc. Ser.\emph{\ }\textbf{A24}
(1962) 259--265.

\bibitem{McCrimmon:1969} K.~McCrimmon, \textit{The Freudenthal-Springer-Tits
construction of exceptional Jordan algebras}, {Trans. Amer. Math. Soc.
\textbf{139} (1969) 495--510}.

\bibitem{McCrimmon:2004} K.~McCrimmon : \textit{\textquotedblleft A Taste of
Jordan Algebras"}, Springer-Verlag New York Inc., New York, 2004.

\bibitem{Baez:2001dm} J.~C. Baez, \textit{The octonions},\
%\href{http://dx.doi.org/10.1090/S0273-0979-01-00934-X}%
{Bull. Amer. Math. Soc. \textbf{39} (2002) 145--205}, {\texttt{%
arXiv:math/0105155}}.

\bibitem{loos1} O. Loos : \textit{\textquotedblleft Jordan Pairs"}, Lect.
Notes Math. \textbf{460}, (Springer, 1975).

\bibitem{tits1} J. Tits,\textit{\ Une classe d'alg\`{e}bres de Lie en
relation avec les alg\`{e}bres de Jordan}, Nederl. Akad. Wetensch. Proc.
Ser. \textbf{A 65} = Indagationes Mathematicae \textbf{24}, 530 (1962).

\bibitem{kantor1} I. L. Kantor, \textit{Classification of irreducible
transitive differential groups}, Doklady Akademiii Nauk SSSR \textbf{158},
1271 (1964).

\bibitem{koecher1} M. Koecher, \textit{Imbedding of Jordan algebras into Lie
algebras I.}, Am. J. Math. \textbf{89}, 787 (1967).

\bibitem{faulk} J. R. Faulkner, \textit{Jordan pairs and Hopf algebras}, J.
of Algebra \textbf{232}, 152 (2000).

\bibitem{Schafer:1966} R.~Schafer : \textit{\textquotedblleft Introduction
to Nonassociative Algebras"}, Academic Press Inc., New York, 1966.

\bibitem{Schafer-2} R. D. Schafer, \textit{Inner derivations of non
associative algebra}s, Bull. Amer. Math. Soc. \textbf{55}, 769 (1949).

\bibitem{Brown:1969} R.~B. Brown, \textit{Groups of type }$E_{7}$,\ J. Reine
Angew. Math. \textbf{236} (1969) 79--102.

\bibitem{Koecher} M. Koecher : \textquotedblleft \textit{An elementary
approach to bounded symmetric domains"}, Rice Univ. lectures, Houston (1969).

\bibitem{Loos} O. Loos : \textit{\textquotedblleft Bounded symmetric domains
and Jordan pairs"}, lectures, Univ. California, Irvine (1977).

\bibitem{FRT} H. Freudenthal, \textit{Lie groups in the foundations of
geometry}, Adv. Math. \textbf{1}, 145 (1963). J. Tits, \textit{Alg\`{e}bres
alternatives, alg\`{e}bres de Jordan et alg\`{e}bres de Lie exceptionnelles.
I. Construction}, (French), Nederl. Akad. Wetensch. Proc. Ser. \textbf{A 69}%
, 223 (1966). B. A. Rozenfeld, \textit{Geometrical interpretation of the
compact simple Lie groups of the class }$\mathit{E}$ (in Russian), Dokl.
Akad. Nauk. SSSR \textbf{106}, 600 (1956).

\bibitem{Freudenthal:1954} H.~Freudenthal, \textit{Beziehungen der }$E_{7}$%
\textit{\ und }$E_{8}$\textit{\ zur oktavenebene I-II}, Nederl. Akad.
Wetensch. Proc. Ser. \textbf{57} (1954) 218--230.

\bibitem{Faulkner:1971} J.~R. Faulkner, \textit{A Construction of Lie
Algebras from a Class of Ternary Algebras}, {\textit{Trans. Amer. Math. Soc.}
\textbf{155} (1971) no. 2, 397--408}.

\bibitem{Ferrar:1972} C.~J. Ferrar, \textit{Strictly Regular Elements in
Freudenthal Triple Systems},\ {\textit{Trans. Amer. Math. Soc.} \textbf{174}
(1972) 313--331}.

\bibitem{F-Gimon-K} S. Ferrara, E. G. Gimon and R. Kallosh, \textit{Magic
supergravities, }$\mathit{N=8}$\textit{\ and black hole composites}, Phys.
Rev. \textbf{D74} (2006) 125018, \texttt{hep-th/0606211}.

\bibitem{Wissanji} K. Dasgupta, V. Hussin and A. Wissanji, \textit{%
Quaternionic K\"{a}hler Manifolds, Constrained Instantons and the Magic
Square. I}, Nucl. Phys. \textbf{B793} (2008) 34-82, \texttt{arXiv:0708.1023
[hep-th]}.

\bibitem{Truini-1} P. Truini, \textit{Exceptional Lie Algebras, }$\mathit{%
SU(3)}$\textit{\ and Jordan Pairs}, Pacific J. Math. \textbf{260}, 227
(2012), \texttt{arXiv:1112.1258 [math-ph]}.

\bibitem{FMZ} S. Ferrara, A. Marrani and B. Zumino, \textit{Jordan Pairs, }$%
\mathit{E}_{6}$\textit{\ and }$\mathit{U}$\textit{-Duality in Five Dimensions%
}, J. Phys. \textbf{A46} (2013) 065402, \texttt{arXiv:1208.0347 [math-ph]}.

\bibitem{Truini-2} A. Marrani and P. Truini, \textit{Exceptional Lie
Algebras, }$\mathit{SU(3)}$\textit{\ and Jordan Pairs Part 2: Zorn-type
Representations}, J. Phys. \textbf{A47} (2014) 265202, \texttt{%
arXiv:1403.5120 [math-ph]}.

\bibitem{Gunaydin:2000xr} M.~G{\"{u}}naydin, K.~Koepsell, and H.~Nicolai,
\textit{Conformal and quasiconformal realizations of exceptional Lie groups}%
,\ Commun. Math. Phys. \textbf{\ 221} (2001) 57--76, {\texttt{hep-th/0008063}%
}.

\bibitem{Gun-Bars} I. Bars and M. G{\"{u}}naydin, \textit{Dynamical Theory
Of Subconstituents Based On Ternary Algebras}, Phys. Rev. \textbf{D22}
(1980), 1403-1413.

\bibitem{GST} M. G{\"{u}}naydin, G. Sierra and P. K. Townsend, \textit{%
Exceptional Supergravity Theories and the Magic Square}, Phys. Lett. B133,
72 (1983). M. G{\"{u}}naydin, G. Sierra and P. K. Townsend, \textit{The
Geometry of }$\mathit{N=2}$\textit{\ Maxwell-Einstein Supergravity and
Jordan Algebras}, Nucl. Phys. \textbf{B242}, 244 (1984).

\bibitem{Squaring-Magic} S. L. Cacciatori, B. L. Cerchiai and A. Marrani,
\textit{Squaring the Magic}, Adv. Theor. Math. Phys. \textbf{19} (2015)
923-954, \texttt{arXiv:1208.6153 [math-ph]}.

\bibitem{BS} C. H. Barton and A. Sudbery, \textit{Magic Squares of Lie
Algebras}, \texttt{arXiv:math/0001083}.

\bibitem{BaSu} C.H. Barton and A. Sudbery, Magic squares and matrix models
of Lie algebras, Adv. Math. \textbf{180} (2003) 596-647.

\bibitem{Gunaydin:1975mp} M.~G{\"{u}}naydin, \textit{Exceptional
realizations of Lorentz group: Supersymmetries and leptons}, Nuovo Cim.
\textbf{A29} (1975) 467.

\bibitem{Gunaydin:1989dq} M.~G{\"{u}}naydin, \textit{The exceptional
superspace and the quadratic Jordan formulation of quantum mechanics}, in :
\textit{\textquotedblleft Elementary particles and the universe: Essays in
honor of Murray Gell-Mann\textquotedblright }, Pasadena 1989, pp. 99-119.,
ed. by J. Schwarz, Cambridge University Press.

\bibitem{Gunaydin:1992zh} M.~G{\"{u}}naydin, \textit{Generalized conformal
and superconformal group actions and Jordan algebras},\ Mod. Phys. Lett.
\textbf{A8} (1993) 1407--1416, {\texttt{hep-th/9301050}}.

\bibitem{Gunaydin:2005zz} M.~G{\"{u}}naydin and O.~Pavlyk, \textit{%
Generalized spacetimes defined by cubic forms and the minimal unitary
realizations of their quasiconformal groups},\ \emph{\ }JHEP \textbf{08}
(2005) 101, {\texttt{hep-th/0506010}}.

\bibitem{j} E. Cremmer and B. Julia, \textit{The }$\mathit{N=8}$\textit{\
Supergravity Theory. 1. The Lagrangian}, Phys. Lett. \textbf{B80}, 48
(1978). E. Cremmer and B. Julia, \textit{The }$\mathit{SO(8)}$\textit{\
Supergravity}, Nucl. Phys. \textbf{B159}, 141 (1979).

\bibitem{HT} C. Hull and P. K. Townsend, \textit{Unity of Superstring
Dualities}, Nucl. Phys. \textbf{B438}, 109 (1995), \texttt{hep-th/9410167}.

\bibitem{Gun-2} M. G{\"{u}}naydin and O. Pavlyk, \textit{Spectrum Generating
Conformal and Quasiconformal }$\mathit{U}$\textit{-Duality Groups,
Supergravity and Spherical Vectors}, JHEP \textbf{1004} (2010) 070, \texttt{%
arXiv:0901.1646 [hep-th]}.

\bibitem{R-map} B. de Wit, F. Vanderseypen and A Van Proeyen, \textit{%
Symmetry structure of special geometries}, Nucl. Phys. \textbf{B400} (1993)
463-524, \texttt{hep-th/9210068}. D. V. Alekseevsky and V. Cort\'{e}s,
\textit{Geometric construction of the }$r$\textit{-map: from affine special
real to special K\"{a}hler manifolds}, Comm. Math.Phys. \textbf{291} (2009)
579-590, \texttt{arXiv:0811.1658 [math.DG]}.

\bibitem{timelike-reduction} E. Cremmer, I. V. Lavrinenko, H. Lu , C. N.
Pope, K.S. Stelle and T. A. Tran, \textit{Euclidean signature
supergravities, dualities and instantons}, Nucl. Phys. \textbf{B534} (1998)
40-82, \texttt{hep-th/9803259}.

\bibitem{BGM} P. Breitenlohner, G. W. Gibbons and D. Maison, \textit{%
Four-Dimensional Black Holes from Kaluza-Klein Theories}, Commun. Math.
Phys. \textbf{120} (1988) 295.

\bibitem{c-map} S. Cecotti, S. Ferrara and L. Girardello, \textit{Geometry
of Type II Superstrings and the Moduli of Superconformal Field Theories},
Int. J. Mod. Phys. \textbf{A4} (1989) 2475.

\bibitem{super-Ehlers} S. Ferrara, A. Marrani and M. Trigiante, \textit{%
Super-Ehlers in Any Dimension}, JHEP \textbf{1211} (2012) 068, \texttt{%
arXiv:1206.1255 [hep-th]}.

\bibitem{Keurentjes} A. Keurentjes, \textit{The Group theory of oxidation},
Nucl. Phys. \textbf{B658}, 303 (2003), \texttt{hep-th/0210178}. A.
Keurentjes, \textit{The Group theory of oxidation 2: Cosets of nonsplit
groups}, Nucl. Phys. \textbf{B658}, 348 (2003), \texttt{hep-th/0212024}.

\bibitem{CFMZ1-D=5} B. L. Cerchiai, S. Ferrara, A. Marrani and B. Zumino,
\textit{Charge Orbits of Extremal Black Holes in Five Dimensional
Supergravity}, Phys. Rev. \textbf{D82}, 085010 (2010), \texttt{%
arXiv:1006.3101 [hep-th]}.

\bibitem{MCD} S. L. Cacciatori, B. L. Cerchiai and A. Marrani, \textit{Magic
Coset Decompositions}, Adv. Theor. Math. Phys. \textbf{17} (2013) 1077-1128,
\texttt{arXiv:1201.6314 [hep-th]}.

\bibitem{Kugo-Townsend} T. Kugo and P. K. Townsend, \textit{Supersymmetry
and the Division Algebras}, Nucl. Phys. \textbf{B221} (1983) 357.

\bibitem{magic-D=6} M. G\"{u}naydin, H. Samtleben and E. Sezgin, \textit{On
the Magical Supergravities in Six Dimensions}, Nucl. Phys. \textbf{B848}
(2011) 62-89, \texttt{arXiv:1012.1818 [hep-th]}.

\bibitem{Mkrtchyan-Nersessian} R. Mkrtchyan, A. Nersessian and V. Yeghikyan,
\textit{Hopf Maps and Wigner's Little Groups}, Mod. Phys. Lett. \textbf{A26}%
, 1393 (2011), \texttt{arXiv:1008.2589 [hep-th]}.

\bibitem{DV-1}
  M.~Dubois-Violette, \textit{Exceptional quantum geometry and particle physics},
  Nucl.\ Phys.\ B {\bf 912} (2016) 426, \texttt{arXiv:1604.01247 [math.QA]}.

\bibitem{DV-2}
  A.~Carotenuto, L.~Dabrowski and M.~Dubois-Violette, \textit{Differential calculus on Jordan algebra and Jordan modules},
  Lett.\ Math.\ Phys.\  {\bf 109} (2019) no.1,  113, \texttt{arXiv:1803.08373 [math.QA]}.

\bibitem{Todorov-1}
  I.~Todorov and S.~Drenska, \textit{Octonions, exceptional Jordan algebra and the role of the group $F_4$ in particle physics},
  Adv.\ Appl.\ Clifford Algebras {\bf 28} (2018) no.4,  82, \texttt{arXiv:1805.06739 [hep-th]}.

\bibitem{Todorov-2}
  I.~Todorov and M.~Dubois-Violette, \textit{Deducing the symmetry of the standard model from the automorphism and structure groups of the exceptional Jordan algebra}, Int.\ J.\ Mod.\ Phys.\ A {\bf 33} (2018) no.20,  1850118, \texttt{arXiv:1806.09450 [hep-th]}.

\bibitem{Todorov-3}
  M.~Dubois-Violette and I.~Todorov, \textit{Exceptional quantum geometry and particle physics II},
  Nucl.\ Phys.\ B {\bf 938} (2019) 751, \texttt{arXiv:1808.08110 [hep-th]}.


\bibitem{GeNa} I.M. Gelfand and M.A. Naimark, Unitary representations of the
Lorentz group, Acad. Sci. USSR. J. Phys. \textbf{10} (1946) 93-94.

\bibitem{GeV} I.M. Gelfand, M.I. Graev and N.Y. Vilenkin, \textit{%
Generalised Functions}, vol. 5 (Academic Press, New York, 1966).

\bibitem{HC} Harish-Chandra, "Representations of semisimple Lie groups:
IV,V", Am. J. Math. \textbf{77} (1955) 743-777, \textbf{78} (1956) 1-41.

\bibitem{Har} Harish-Chandra, ``Discrete series for semisimple Lie groups:
II'', Ann. Math. \textbf{116} (1966) 1-111.

\bibitem{Bru} F. Bruhat, ``Sur les representations induites des groups de
Lie'', Bull. Soc. Math. France, \textbf{84} (1956) 97-205.

\bibitem{Lan} R.P. Langlands, \textit{On the classification of irreducible
representations of real algebraic groups}, in: "Representation theory and
harmonic analysis on semi-simple Lie groups", eds. P. Sally and D. Vogan,
Math. Surveys and Monographs, Vol. 31 (AMS, 1989) pp. 101--170; (first as
IAS Princeton preprint, 1973).

\bibitem{KnZu} A.W. Knapp and G.J. Zuckerman, ``Classification theorems for
representations of semisimple groups'', in: Lecture Notes in Math., Vol. 587
(Springer, Berlin, 1977) pp. 138-159; ~``Classification of irreducible
tempered representations of semisimple groups'', Ann. Math. \textbf{116}
(1982) 389-501.

\bibitem{War} G. Warner, \textit{Harmonic Analysis on Semi-Simple Lie Groups
I}, (Springer, Berlin, 1972).

\bibitem{Vog} D. Vogan, \textit{Representations of Real Reductive Lie Groups}%
, Progr. Math., Vol. 15 (Boston-Basel-Stuttgart, Birkh\"auser, 1981).

\bibitem{Sat} I. Satake, ``On representations and compactifications of
symmetric Riemannian spaces'', Ann. Math. \textbf{71} (1960) 77-110.

\bibitem{Bor} A. Borel, "Linear algebraic groups", Benjamin (1969)

\bibitem{BoTi} A. Borel, J. Tits, "Groupes reductifs" Publ. Math. IHES, 27
(1965) pp. 55-150

\bibitem{BoTi2} A. Borel, J. Tits, "Elements unipotents et sous-groupes
paraboliques de groupes reductifs I" Invent. Math., 12 (1971) pp. 95-104

\bibitem{Hum} J.E. Humphreys, "Linear algebraic groups", Springer (1975)

\bibitem{Bou} N. Bourbaki, "Groupes et alg\~{A}%
%TCIMACRO{\U{a8}}%
%BeginExpansion
\"{}%
%EndExpansion
bres de Lie", Hermann (1975) pp. Chapts. VII-VIII

\bibitem{Dobpar} V.K. Dobrev, Invariant Differential Operators for
Non-Compact Lie Groups: Parabolic Subalgebras, Rev. Math. Phys. \textbf{20}
(2008) 407-449, \texttt{hep-th/0702152}.

\bibitem{Dobparab} V.K. Dobrev, Invariant Differential Operators for
Non-Compact Lie Algebras Parabolically Related to Conformal Lie Algebras, J.
High Energy Phys. 02 (2013) 015, \texttt{arXiv:1208.0409}.

\bibitem{Dobk} Vladimir K. Dobrev, \textit{Invariant Differential Operators,
Volume 1: Noncompact Semisimple Lie Algebras and Groups}, De Gruyter Studies
in Mathematical Physics vol. 35 (De Gruyter, Berlin, Boston, 2016).

\bibitem{jaco}
  N. Jacobson: \textit{"Exceptional Lie Algebras"}, Lecture Notes in Pure and Applied Mathematics
{\bf 1} (M. Dekker, 1971).


 \bibitem{myb}
  K. Meyberg, \textit{Jordan-Triplesysteme und die Koecher-Konstruktion von Lie Algebren}, Math.
Z. {\bf 115}, 58 (1970).

\bibitem{Gun-CQC}
  M.~Gunaydin, \textit{Generalized conformal and superconformal group actions and Jordan algebras},
  Mod.\ Phys.\ Lett.\ A {\bf 8} (1993) 1407, \texttt{hep-th/9301050}.

\bibitem{Magic-Non-Susy}
  A.~Marrani, G.~Pradisi, F.~Riccioni and L.~Romano, \textit{Nonsupersymmetric magic theories and Ehlers truncations},
  Int.\ J.\ Mod.\ Phys.\ A {\bf 32} (2017) no.19n20,  1750120, \texttt{arXiv:1701.03031 [hep-th]}.

\bibitem{DF-E6} M. J. Duff and S. Ferrara, $E_{6}$\textit{\ and the
bipartite entanglement of three qutrits}, Phys. Rev. \textbf{D76} (2007)
124023, \texttt{arXiv:0704.0507 [hep-th]}.

\bibitem{Slansky} R. Slansky, \textit{Group Theory for Unified Model Building%
}, Phys. Rept. \textbf{79}, 1 (1981).

\bibitem{BFGM1} S. Bellucci, S. Ferrara, M. G\"{u}naydin and A. Marrani,
\textit{Charge orbits of symmetric special geometries and attractors}, Int.
J. Mod. Phys. \textbf{A21} (2006) 5043-5098, \texttt{hep-th/0606209}.

\bibitem{CFM1} A. Ceresole, S. Ferrara and A. Marrani, \textit{4d/5d
Correspondence for the Black Hole Potential and its Critical Points}, Class.
Quant. Grav. \textbf{24} (2007) 5651-5666, \texttt{arXiv:0707.0964 [hep-th]}.

\bibitem{Small-Orbits} L. Borsten , M. J. Duff, S. Ferrara, A. Marrani and
W. Rubens, \textit{Small Orbits}, Phys. Rev. \textbf{D85} (2012) 086002,
\texttt{arXiv:1108.0424 [hep-th]}.

\bibitem{Patera} W. G. MacKay and J. Patera : \textit{\textquotedblleft
Tables of Dimensions, Indices, and Branching Rules for Representations of
Simple Lie Algebras"}, Marcel Dekker Inc. (New York and Basel, 1981). {}

\bibitem{Minchenko} A. N. Minchenko, \textit{The Semisimple Subalgebras of
Exceptional Lie Algebras}, Trans. Moscow Math. Soc. \textbf{67} (2006)
225--259.

\bibitem{Kephart} R. Feger and T. W. Kephart, \textit{LieART - A Mathematica
Application for Lie Algebras and Representation Theory}, Comput. Phys.
Commun. \textbf{192} (2015) 166-195, \texttt{arXiv:1206.6379 [math-ph]}.

\bibitem{GZ} M. G\"{u}naydin and M. Zagermann, \textit{Unified
Maxwell-Einstein and Yang-Mills-Einstein supergravity theories in
five-dimensions}, JHEP \textbf{0307} (2003) 023, \texttt{hep-th/0304109}. M.
G\"{u}naydin and M. Zagermann, \textit{Unified }$\mathit{N=2}$\textit{\
Maxwell-Einstein and Yang-Mills-Einstein supergravity theories in four
dimensions}, JHEP \textbf{0509} (2005) 026, \texttt{hep-th/0507227}.

\bibitem{FG-D=5} S. Ferrara and M. G\"{u}naydin, \textit{Orbits and
Attractors for }$\mathit{N=2}$\textit{\ Maxwell-Einstein Supergravity
Theories in Five Dimensions}, Nucl. Phys. \textbf{B759} (2006) 1-19, \texttt{%
hep-th/0606108}.

\bibitem{Vogel} P.Vogel, \textit{Algebraic structures on modules of diagrams}%
, J. Pure Appl. Algebra \textbf{215} (2011), no. 6, 1292-1339. P. Vogel,
\textit{The universal Lie algebra}, Preprint (1999),\texttt{\
www.math.jussieu.fr/ vogel/A299.ps.gz}.

\bibitem{LM-1} J. M. Landsberg and L. Manivel, \textit{Triality, exceptional
Lie algebras and Deligne dimension formulas}, Adv. Math. \textbf{171}
(2002), 59-85, \texttt{arXiv:math/0107032 [math.AG]}.

\bibitem{Brown} R. B. Brown, \textit{Groups of type }$\mathit{E}_{7}$, J.
Reine Angew. Math. \textbf{236} (1969) 79--102.

\bibitem{type-E7} S. Ferrara and R. Kallosh, \textit{Creation of Matter in
the Universe and Groups of Type }$\mathit{E}_{7}$, JHEP \textbf{1112} (2011)
096, \texttt{arXiv:1110.4048 [hep-th]}. S. Ferrara and A. Marrani, \textit{%
Black Holes and Groups of Type }$\mathit{E}_{7}$, Pramana \textbf{78} (2012)
893-905, \texttt{arXiv:1112.2664 [hep-th]}. S. Ferrara, R. Kallosh and A.
Marrani, \textit{Degeneration of Groups of Type E7 and Minimal Coupling in
Supergravity}, JHEP \textbf{1206} (2012) 074, \texttt{arXiv:1202.1290
[hep-th]}.

\bibitem{T^3} S. Bellucci, A. Marrani, E. Orazi and A. Shcherbakov, \textit{%
Attractors with Vanishing Central Charge}, Phys. Lett. \textbf{B655} (2007)
185-195, arXiv:0707.2730 [hep-th]. K. Saraikin and C. Vafa, \textit{%
Non-supersymmetric black holes and topological strings}, Class. Quant. Grav.
\textbf{25} (2008) 095007, \texttt{hep-th/0703214}.

\bibitem{DobMar2} V.K. Dobrev, A. Marrani, in preparation.
\end{thebibliography}
\end{document}